\documentclass[11pt]{article} 

\usepackage[utf8]{inputenc} 

\usepackage{geometry} 
\usepackage{graphicx} 

\usepackage{booktabs} 
\usepackage{array} 
\usepackage{paralist} 
\usepackage{verbatim} 
\usepackage{subfig} 

\usepackage[english]{babel}
\usepackage{csquotes}
\usepackage{graphicx} 
\usepackage{amssymb}
\usepackage{amsmath, amsfonts, bbm, dsfont}
\usepackage[all]{xy}

\usepackage{fancyhdr} 
\usepackage{sectsty}
\allsectionsfont{\sffamily\mdseries\upshape} 

\usepackage[nottoc,notlof,notlot]{tocbibind} 
\usepackage[titles,subfigure]{tocloft} 

\usepackage[
backend=biber,
style=alphabetic
]{biblatex}

\usepackage{amsthm}
\usepackage{tikz-cd}

\newtheorem{theo}{Theorem}[section]
\newtheorem{prop}[theo]{Proposition}
\newtheorem{lemma}[theo]{Lemma}
\newtheorem{coro}[theo]{Corollary}

\newtheorem{claim}[theo]{Claim}

\theoremstyle{definition}
\newtheorem{defi}[theo]{Definition}
\newtheorem{rem}[theo]{Remark}
\newtheorem{ex}[theo]{Example}
\newtheorem{conj}[theo]{Conjecture}

\newtheorem{constr}[theo]{Construction}

\title{Extending Hrushovski's groupoid-cover correspondence using simplicial groupoids}
\author{Paul Wang}
\date{}

\begin{document}
\setcounter{tocdepth}{1}
\maketitle

 ABSTRACT. Hrushovski's suggestion, given in [\textquotedblleft Groupoids, imaginaries and internal covers," \textit{Turkish Journal of Mathematics}, 2012], to capture the structure of the 1-analysable covers of a theory $T$ using simplicial groupoids definable in $T$ is realized here. The ideas of Haykazyan and Moosa, found in [\textquotedblleft Functoriality and uniformity in Hrushovski's groupoid-cover correspondence," \textit{Annals of Pure and Applied Logic}, 2018] are used, and extended, to define an equivalence of categories. Finally, a couple of examples are studied with these new tools. \footnote{\, \copyright \,  2021. This manuscript version is made available under the CC-BY-NC-ND 4.0 license http://creativecommons.org/licenses/by-nc-nd/4.0/}

  \tableofcontents

\begin{section}{Introduction}
When studying many-sorted first-order structures, one may pick a structure $\mathbb{U}$, and add a new sort $S$. If this does not create new definable sets in $\mathbb{U}$, the extended structure is called a cover of $\mathbb{U}$. 
Under an additional hypothesis of internality, it was shown by Hrushovski in \cite{hrushovski} that covers came with definable \textit{binding groupoids}, which are a generalized version of automorphism groups. Up to a suitable notion of equivalence, the internal covers of a fixed structure $\mathbb{U}$ are characterized by their binding groupoids.

This correspondence between internal covers and definable groupoids was then extended in \cite{haykazyan-moosa}, in two ways : First, it was shown to come from an equivalence of categories. Then, it was extended to the case of fibrations of internal covers, or 1-analysable covers, under a (rather strong) hypothesis of independence.

In this paper, we pursue the work done in \cite{haykazyan-moosa}. The aim is to extend the correspondence by dropping the hypothesis of independence found in \cite{haykazyan-moosa}. 

First, we recall some results regarding the internal and independent cases. We then write down the details of a proof left to the reader at the end of \cite{haykazyan-moosa}. The aim is to find ideas that can later on be extended to the more general case at hand.

Then, sections 4 and 5 are aimed at realizing Hrushovski's suggestion to deal with non-independent 1-analysable covers by using \textit{simplicial groupoids} instead of groupoids. Several new ideas appear, the most crucial one being that of \textit{coherence} of morphisms in a simplicial groupoid.

Then, two examples of non-internal and non-independent 1-analysable covers are studied, using the new tools given by simplicial groupoids. The aim is to find explicit descriptions of the binding simplicial groupoids of these covers, and deduce a few results.

Finally, in section 7, we drop yet another hypothesis, of finite generatedness of the language, and aim to generalize the correspondence again, using type-definable simplicial groupoids.

\end{section}

\section*{Acknowledgements}

Most of the work presented here was done under the supervision of Martin Bays. His help was invaluable in correcting and improving ideas, as well as for suggesting directions of research.

I would also like to thank Martin Hils, for pointing out a mistake which led to a better understanding of one of the examples presented here.

Other improvements, such as the section on type-definable groupoids, were suggested by Rahim Moosa, who is also one of the authors of a paper which inspired the present work.

Finally, thanks to the Ecole Normale Superieure for the great teachers, and classmates, I have, as well as for the housing and funding provided, which represent a much-appreciated help in my studies.

\begin{section}{Internal covers and 1-analysable covers with independent fibers}

In this section, our goal is to recall some results in the internal and independent cases. We also use this opportunity to give more detailed proofs of the fact that the structures defined in \cite{haykazyan-moosa}(Definition 3.9 and Remark 3.10) are indeed covers. The methods and ideas we use here will be generalized in the later sections.

\begin{subsection}{General definitions and context}

From now on, $\mathbb{U}$ denotes a saturated, i.e. $|\mathbb{U}|$-saturated , model of a complete theory $T$ that admits elimination of imaginaries. For more details on elimination of imaginaries, the reader can have a look at \cite{pillay} (from Lemma 1.2 to Remark 1.9).

    We will assume that all the theories we study have $\kappa$-saturated models in the same cardinality $\kappa$. If the theories are unstable, we may have to add axioms to ZFC to get such saturated models. Improving our proofs in order to get rid of this hypothesis might or might not be easy.

    \begin{defi}
    Let $T'$ be a complete extension of $T$, possibly with additional sorts. Let $\mathbb{U}'$ be a saturated model of $T'$, and $\mathbb{U}:= \mathbb{U'}|_T$. We say that $\mathbb{U}$ is \textit{stably embedded} in $\mathbb{U}'$ when the map $\mathrm{Aut}(\mathbb{U}')\rightarrow\mathrm{Aut}(\mathbb{U}) $ is surjective.

    In this case, we also say that $T'$ is a \textit{cover} of $T$, or that $\mathbb{U'}$ is a cover of $\mathbb{U}$.

    If there exist bijections between the new sorts and some definable sets of $\mathbb{U}$ that are definable, \textit{possibly with parameters}, the cover is said to be \textit{internal}.
    \end{defi}

    \begin{prop}
    Let $T'$ be a complete extension of $T$, possibly with additional sorts. Let $\mathbb{U}'$ be a saturated model of $T'$, and $\mathbb{U}:= \mathbb{U'}|_T$. The following are equivalent : 
    
    \begin{enumerate}
        \item $\mathbb{U}$ is stably embedded in $\mathbb{U}'$.
        \item Any subset of a product of sorts of $\mathbb{U}$ that is definable with parameters from $\mathbb{U}'$ is definable with parameters from $\mathbb{U}$.
        \item For all tuples $b \subseteq \mathbb{U}'$ of small length, there exists a small subset $A$ of $\mathbb{U}$ such that $tp(b / A) \models tp(b/ \mathbb{U})$.
    \end{enumerate}
    
    \end{prop}

    See \cite{chatz-hru} (Appendix, Lemma 1) for a proof of the equivalences.

    \begin{defi}
    A \textit{category} is a two-sorted structure, with a set of objects, a set of maps with source and target, an associative composition, and an identity map for each object.
    \end{defi}
    
    The following definitions are the same as those in \cite{haykazyan-moosa}. 

    \begin{defi}\label{defi_groupoid} \begin{enumerate}
        \item A \textit{groupoid} is a category where each morphism has an inverse. It is \textit{connected} if two arbitrary objects are isomorphic. It is \textit{canonical} if two isomorphic objects are equal.
    
    \item A definable concrete groupoid is given by a definable set of definable objects $(O_i)_{i \in I}$, and a definable set $M$ of definable bijections between these objects. In other words, there are definable sets $O, I, M, R$, a definable surjective map $f : O \rightarrow I$ such that $R \subseteq M \times O \times O$ and for all $m \in M$, the definable set $R_m \subseteq O \times O$ is the graph of a bijection from $O_i$ to $O_j$, for some $i,j \in I$. Here, $O_i$ denotes the fiber $f^{-1}(i) \subseteq O$.
    
    We require that the bijections define a category which is a groupoid. In particular, we require that the inverse of a bijection coded by some element in $M$ be also coded by an element in $M$. We also assume that the maps defined by the elements of $m$ are pairwise distinct. Note that the definition implies that the objects $O_i$ are pairwise disjoint. 
    
    The groupoid is \textit{definable over X} if the sets $O, I, M, R$ are definable over $X$. It is $0$-definable if it is definable over the empty set.

    \item A concrete groupoid $\mathcal{G}$ is \textit{finitely faithful} if, for each object $O$ of the groupoid $\mathcal{G}$, there exists a finite tuple $\overline{x}$ of elements of $O$ whose pointwise stabilizer in the group $\mathrm{Aut}_{\mathcal{G}}(O)$ is trivial.
    
    \item If $A$ is a definable set, a definable groupoid over $A$ is a uniformly definable family $(G_a)_{a \in A}$ of definable connected groupoids that are pairwise disjoint. See \cite{haykazyan-moosa}(Definition 5.1)
    \end{enumerate}
    \end{defi}

\begin{defi}
\begin{enumerate}
    \item If $A$ is a set, and $n$ is an integer, we let $A^{=n}$ denote the set of tuples $(a_1,...,a_n)$, where the $a_i$ are elements of $A$ that are pairwise distinct. We also let $A^{<\omega}$ denote the collection of finite subsets of $A$.
    \item If $\overline{c}$ is a finite tuple, we use the notation \textquotedblleft $a \in \overline{c}$" to denote \textquotedblleft $a$ is some coordinate of the tuple $\overline{c}$". We also use the notation \textquotedblleft$\overline{c} \subseteq \overline{d}$", where $\overline{c}, \overline{d}$ are tuples, to denote \textquotedblleft for all elements $a$, if $a \in \overline{c}$, then $a \in \overline{d}$".
    \item If $f : S \rightarrow A$ is a map, if $n$ is an integer and $\overline{c}$ is an element of $A^{=n}$, we define the set $S_{\overline{c}}$ as $S_{\overline{c}} := \bigcup\limits_{a \in \overline{c}} S_a  \subseteq S$.
\end{enumerate}
\end{defi}
    
        \begin{defi}\label{defi_covers} \begin{enumerate}
        \item Let $(\mathbb{U}, S)$ be a cover of $\mathbb{U}$, and $A \subseteq \mathbb{U}$ be a 0-definable set. The cover is \textit{1-analysable over A} if there exists a 0-definable surjective map $f : S \rightarrow A $ such that each fiber $S_a$ is internal to $\mathbb{U}$. In other words, $f$ and $S$ define a family of internal covers $(\mathbb{U}, S_a)$. The data of the map $f$ is considered as part of the data of the 1-analysable cover.
         
        \item If, additionally, we have $\mathrm{Aut}(S / \mathbb{U})=\prod_a \mathrm{Aut}(S_a / \mathbb{U})$, we say that the fibers are independent.
        
        \item Let $(\mathbb{U}, S)$ be a cover of $\mathbb{U}$, with a 0-definable relation $R \subseteq B \times S$, where $B$ is 0-definable in $\mathbb{U}$.
        If $b \in B$, we define $R_b \subseteq S$ as the set of elements $s \in S$ such that $\models R(b,s)$. We also define the structure $(\mathbb{U}, R_b)$ as that induced by the $b$-definable sets.
        
        If there exist finitely many relations $R_1(\overline{x}_1, y),..., R_k(\overline{x}_k, y)$ that are 0-definable in $(\mathbb{U},S)$ and such that, for each element $b$ in $B$, the relations $R_i(\overline{x}_i, b)$ generate the language of $(\mathbb{U}, R_b)$ over that of $\mathbb{U}$, then we say that the language of $(\mathbb{U}, R_b)$ is \textit{finitely generated} over that of $\mathbb{U}$, \textit{uniformly for $b$ in $ B$}.
        
        A special case of this definition is that where $B=A^{=n}$, and there exists a 0-definable map $f : S \twoheadrightarrow A$ such that $R(b,s)$ is defined as \textquotedblleft$f(s) $ is one of the coordinates of the tuple $b$". In that specific case, we say that the language of $(\mathbb{U}, S_b)$ is uniformly finitely generated over that of $\mathbb{U}$.
        
        \end{enumerate}
    \end{defi}
    
    \begin{rem}
    
        Note that the condition \textquotedblleft the language of $(\mathbb{U}, S_a)$ is finitely generated over that of $\mathbb{U}$, uniformly in $a \in A$" is not the same as finite generatedness of the language of the whole cover $(\mathbb{U},S)$ over that of $\mathbb{U}$. 
        
        The second example in section \ref{section_computations} can be rearranged to build a cover satisfying the latter condition, but not the former. See propositions \ref{no_uniform_finite_generatedness_DCF_0} for the main ideas.

    \end{rem}

    \begin{ex}
    If $A,K$ are groups that are 0-definable in $\mathbb{U}$, and $S$ is a group, with a 0-definable exact sequence $1 \rightarrow K \rightarrow S \rightarrow A \rightarrow 1$, then $(\mathbb{U}, S)$ is 1-analysable over $A$, with non-independent fibers in general.
    \end{ex}
\end{subsection}

\begin{subsection}{Connected groupoids}

We first recall the existence of the binding groupoid of an internal cover :

    \begin{theo}[\cite{hrushovski}, \cite{haykazyan-moosa}]\label{binding_groupoid_of_internal_cover}
    Let $\mathbb{U}' = (\mathbb{U}, S)$ be an internal cover of $\mathbb{U}$ whose language is finitely generated over that of $\mathbb{U}$. Then there exist 0-definable connected groupoids $\mathcal{G}$ in $\mathbb{U}$ and $\mathcal{G}'$ in $\mathbb{U}'^{eq}$, satisfying the following conditions :

    \begin{enumerate}
        \item The groupoid $\mathcal{G}'$ is an extension of $\mathcal{G}$, with $S$ as the only added object.
        \item The automorphism group of $S$ over $\mathbb{U}$ satisfies $\mathrm{Aut}_{\mathbb{U}'}(S / \mathbb{U}) = \mathrm{Aut}_{\mathcal{G}'}(S).$
    \end{enumerate}
    \end{theo}

We give the same proof as in \cite{haykazyan-moosa}(Theorem 3.5), because one idea that appears in said proof will be used later on.

\begin{proof}
    Let us first pick a definable bijection $f_b : O_a \rightarrow S$, where $a \in \mathbb{U}$ and $b \in \mathbb{U}'$. We wish to use $f$ and $O$ to define a groupoid. However, we need to \textit{restrict the set of parameters} that are allowed, in order to satisfy condition 2.

    Let $q:= tp(b / \mathbb{U}).$ By stable embeddedness, there exists a small $A \subseteq \mathbb{U}$ such that $p:=tp(b/A) \models q$. If $c,d \models p$, then $c \equiv_{\mathbb{U}}d$. In particular, $f_c, f_d$ are bijections from $O_a$ to $S$, thus the map $(f_d\circ f_c^{-1}) \cup id_{\mathbb{U}}$ is a well-defined automorphism of $\mathbb{U}'$ over $\mathbb{U}$.

    We have thus proved : $p(x) \cup p(y) \models $\textquotedblleft$(f_y\circ f_x^{-1})\cup id_{\mathbb{U}} \in \mathrm{Aut}_{}(\mathbb{U}' / \mathbb{U})$". 
    
    By compactness, we find a formula $\psi(z, \alpha) \in p$ such that $\psi(x,\alpha)\wedge \psi(y, \alpha) \models $\textquotedblleft$(f_y\circ f_x^{-1})\cup id_{\mathbb{U}} \in \mathrm{Aut}_{}(\mathbb{U}' / \mathbb{U})"$. We may assume $\alpha = a$. The following formula enables us to restrict the set of objects :  $\theta(t) := \forall x,y \, [\psi(x, t)\wedge \psi(y, t) \rightarrow $\textquotedblleft$(f_y\circ f_x^{-1})\cup id_{\mathbb{U}} \in \mathrm{Aut}_{}(\mathbb{U}' / \mathbb{U})"]$

    The objects of $\mathcal{G}'$ are the $O_{a'}$, where $\models \theta(a')$, along with $S$. The morphisms from $O_{a'}$ to $S$ are the $f_{b'}$, where $\models \psi(b', a')$. The other morphisms are defined using compositions and inverses.
    
\begin{center}
    
\begin{tikzcd}

O_{a_1} \arrow[rr, "{f_{b_2}^{-1}\circ f_{b_1}}"] \arrow[dr, "{f_{b_1}}"] && O_{a_2}     &&   S \arrow[rr, "{f_{b'}\circ f_{b}^{-1}}"] \arrow[dr, "{f_{b}^{-1}}"] && S    \\& S   \arrow[ru, "{f_{b_2}^{-1}}"]&&& & O_a   \arrow[ru , "f_{b'}"]
\end{tikzcd}

\end{center}

    Let's check that $Hom(O_{a'}, S) \circ Hom(S, O_{a'}) = \mathrm{Aut}_{}(S / \mathbb{U})$, for any $a'$ satisfying $\theta(t)$.

    The properties of $\theta$ and $\psi$ imply that $Hom(O_{a'}, S) \circ Hom(S, O_{a'}) \subseteq \mathrm{Aut}_{}(S / \mathbb{U})$. 
    
    Conversely, if $\sigma \in \mathrm{Aut}_{}(S / \mathbb{U}),$ then, for any $f_{b'} \in Hom(O_{a'}, S)$, we have $\sigma = f_{\sigma(b')}\circ f_{b'}^{-1} \in  Hom(O_{a'}, S) \circ Hom(S, O_{a'}) $.

    Finally, by stable embeddedness, the groupoid $\mathcal{G}:=\mathcal{G}'|_{\mathbb{U}}$ is 0-definable in $\mathbb{U}$.
    
    \end{proof}

\begin{rem}\label{types_of_bijections}
One of the main ideas in this theorem is to \textit{possibly add many objects} that classify the various \textquotedblleft kinds of bijections" between $S$ and the $O_a$, to ensure that two bijections $f_b$ and $f_{b'}$ of the \textquotedblleft same kind" yield a model-theoretic automorphism $f_{b'}\circ f_b^{-1}$.  
\end{rem}

\medskip

Let $\mathcal{G}$ be a 0-definable connected groupoid in $\mathbb{U}$. 
For the remainder of this subsection, we let $C(\mathcal{G}) = (\mathbb{U}, O_*, M_*, R_*)$ denote the extension of $\mathbb{U}$ built from $\mathcal{G}$ as defined in \cite{haykazyan-moosa}(3.9). Let us briefly recall the construction :

Let $O, I, M$ be sets defining the groupoid $\mathcal{G}$ as in Definition \ref{defi_groupoid}.
Let $O_i$ be some arbitrary object in $\mathcal{G}$. Let $O_*$ be a new sort, whose underlying set is a set-theoretic copy of the set defined by $O_i$. Namely, fix a bijection $f : O_i \rightarrow O_*$. Since $M$ is the set of morphisms in $\mathcal{G}$, let $M_*$ be a new sort, whose underlying set is $M\times \lbrace 0, 1 \rbrace \sqcup \mathrm{Aut}_{\mathcal{G}}(O_i) \times \lbrace 2 \rbrace$.

Let $R_* \subseteq (M_*) \times (O \cup O_*) \times (O \cup O_*) $ be a ternary relation, interpreting the elements of 
$M_*$ as bijections between the $O_j$ and $O_*$, or as permutations of $O_*$. More precisely, the elements of $M \times \lbrace 0 \rbrace$ are interpreted as bijections $O_j \rightarrow O_*$, corresponding to morphisms in $\mathrm{Hom}_{\mathcal{G}}(O_j, O_i)$, those in $M \times \lbrace 1 \rbrace$ are interpreted as bijections $O_* \rightarrow O_j$, corresponding to morphisms in $\mathrm{Hom}_{\mathcal{G}}(O_i, O_j)$. Finally, the elements of $\mathrm{Aut}_{\mathcal{G}}(O_i) \times \lbrace 2 \rbrace$ are interpreted as permutations of $O_*$.
This interpretation uses the bijection $f : O_i \rightarrow O_*$ we chose earlier.

Then, by construction, there is a 0-definable connected groupoid $\mathcal{G}'$ in $C(\mathcal{G})$, which is an extension of $\mathcal{G}$ with one extra object.

\begin{prop}
The structure $C(\mathcal{G})$ is interpretable in $\mathbb{U}$ and thus saturated.
\end{prop}
\begin{proof}
Using constants to create duplicates, we are able to follow the construction in \cite{haykazyan-moosa}(3.9) and interpret the whole structure of $C(\mathcal{G})$ in $\mathbb{U}$. For instance, an element $x$ in the sorts of $\mathbb{U}$ is interpreted by $x\times \lbrace0 \rbrace$, whereas an element in $O_*$ is interpreted by an element of $\mathbb{U}\times \lbrace 1 \rbrace$, the first coordinate itself being given by the construction in \cite{haykazyan-moosa}. Similarly, a morphism in $M_*$ is interpreted by an element of $\mathbb{U}\times \lbrace 2, 3 \rbrace$.
The new relation $R_*$ is interpreted the obvious way, again following the construction of $C(\mathcal{G})$.
\end{proof}

\begin{prop}
The structure $C(\mathcal{G})$ yields an internal cover of $\mathbb{U}$. 
\end{prop}
\begin{proof}
From the construction, it follows that the new sorts are definably isomorphic to definable sets in $\mathbb{U}$. In fact, all the morphisms in the extended groupoid $\mathcal{G}'$ with source $O_*$ yield definable bijections with $O_*$. As for the set of morphisms $M_*$, it is almost in bijection with the set of morphisms in $\mathcal{G}$. Using elimination of imaginaries, one can add the missing copies and get a definable subset $X$ of a product of sorts in $\mathbb{U}$, with a definable bijection $X \simeq M_*$. Thus, since $C(\mathcal{G})$ is saturated, it remains to show that no new structure is induced on the sorts of $\mathbb{U}$. To do this, we shall prove that automorphisms of $\mathbb{U}$ can be extended to automorphisms of $C(\mathcal{G})$.

Let $\sigma \in \mathrm{Aut}(\mathbb{U})$. Keeping the notations in \cite{haykazyan-moosa}, let $O_i$ be the object whose copy is $O_*$. Let $f_c : O_i \rightarrow O_{\sigma(i)}$ be a morphism in the original groupoid. We shall define an action on the elements of $O_*$ and on the morphisms in $M_*$.

If $x'$ is an element of $O_*$, let $x$ be the element of $O_i$ whose copy is $x'$. Then, we shall define $\sigma(x') := (f_c^{-1} \sigma(x))'$, i.e., it is the copy in $O_*$ of the image of $x$ under $f_c^{-1} \circ \sigma$.  

Now, if $m$ is a morphism in $M_*$, we shall distinguish different cases, depending on its domain and codomain :
\begin{itemize}
    \item If neither the domain nor the codomain of $m$ is $O_*$, then we let $\sigma$ act on $m$ as it did in $\mathbb{U}$.
    \item If $m : O_* \rightarrow O_j$ for any $j$ in $\mathbb{U}$, let $n : O_i \rightarrow O_j$ be the copy of $m$. Then $\sigma(n)\circ f_c : O_i \rightarrow O_{\sigma(j)}$. Therefore, let $\sigma(m) : O_* \rightarrow O_{\sigma(j)}$ be the copy of $\sigma(n)\circ f_c$.
    \item Similarly, if $m : O_j \rightarrow O_*$ for any $j$ in $\mathbb{U}$, let $\sigma(m) : O_{\sigma(j)} \rightarrow O_*$ be the copy of $f_c^{-1} \circ \sigma(n)$, using the same notations as before.
    \item Finally, if $m : O_* \rightarrow O_*$, let $\sigma(m) : O_* \rightarrow O_*$ be the copy of $f_c^{-1} \circ \sigma(n) \circ f_c$.
\end{itemize}

Now, studying the various cases, we can show that this extension of $\sigma$ preserves the relation $R$ that defines the behaviour of the bijections in $M_*$. For instance, using the same notations as before, let $m : O_* \rightarrow O_j$, $x' \in O_*$ and $y = m(x') \in O_j$. We compute 

$\sigma(m)(\sigma(x'))=\sigma(n)\circ f_c (\sigma(x'))=\sigma(n) \circ f_c (f_c^{-1}\circ \sigma (x)) = \sigma(n) (\sigma(x))=\sigma(n(x)),$ where $x \in O_i$ is the copy of $x'$.
Since we had assumed that $m(x') =  n(x) = y$, we have the equality $\sigma(m)(\sigma(x')) = \sigma(m(x'))$.

Finally, it is not hard to see what the inverse of $\sigma$ should be. For the action on $O_*$, instead of applying $\sigma$ and then $f_c^{-1}$, we should apply $f_c$ and then $\sigma^{-1}$. The action on the morphisms is defined in an analogous way, using $f_c^{-1}$ and $\sigma^{-1}$.
\end{proof}

\end{subsection}

\begin{subsection}{Non-connected groupoids}

\begin{defi}

\begin{enumerate}
    \item Let $A$ be a 0-definable set. We define the notion of \textit{0-definable families of connected groupoids over $A$} just as in \cite{haykazyan-moosa} (Definition 5.1).
    
    More formally : a \textit{0-definable family of connected groupoids over $A$} is a 0-definable groupoid $\mathcal{G}$, given by 0-definable sets $O, I, M, R$ just as in Definition \ref{defi_groupoid}, such that there exists a $0$-definable surjection $g : I \rightarrow A$ that satisfies the following property : for all $i,j \in I$, the objects $O_i$ and $O_j$ are in the same connected component of $\mathcal{G}$ if and only if $g(i)=g(j)$.
    
    In other words, the connected components of $\mathcal{G}$ are $0$-definably indexed by $A$.
    
    We shall also refer to these groupoids as \textit{groupoids over $A$}.
    \item If $\mathcal{G}$ is a groupoid over $A$, and $a$ is an element of $A$, we let $\mathcal{G}_a$ denote the connected component of $\mathcal{G}$ over $a$. It is an $a$-definable connected groupoid.
\end{enumerate}

\end{defi}

\begin{constr}
Let $\mathcal{G} \rightarrow A$ be a 0-definable family of connected groupoids in $\mathbb{U}$. Let $C(\mathcal{G})$ be the structure defined in \cite{haykazyan-moosa}(5.6). 

As in the case of a connected groupoid, it is defined by choosing one object in each isomorphism class of $\mathcal{G}$, and creating appropriate copies to yield a one-object extension of each isomorphism class.

Let us be a bit more explicit.

Let $O, I, M$ be sets defining the groupoid $\mathcal{G}$ as in Definition \ref{defi_groupoid}.
For each element $a \in A$, let $O_{i_a}$ be some arbitrary object in $\mathcal{G}_a$. Let $O_{*, a}$ be a set-theoretic copy of the set defined by $O_{i_a}$. Namely, fix a bijection $f_a : O_{*, a} \rightarrow O_{i_a}$. Let $O_*$ be a new sort, whose underlying set is the disjoint union of the $O_{*, a}$ for $a$ in $A$. Similarly, for each $a$ in $A$, let $M_{*, a}$ be $M_a\times \lbrace 0, 1 \rbrace \sqcup \mathrm{Aut}_{\mathcal{G_a}}(O_{i_a}) \times \lbrace 2 \rbrace$. Let $M_* = \bigcup_{a \in A} M_{*,a}$ be a new sort.
Finally, similarly to the connected case, let $R_* \subseteq (M_*) \times (O \cup O_*) \times (O \cup O_*) $ be a ternary relation, interpreting the elements of 
$M_*$ as bijections between the $O_j$ and the $O_{*, a}$, or as permutations of the $O_{*,a}$, for $a \in A$, $O_j$ an object of $\mathcal{G}_a$. This interpretation uses the bijections $f_a : O_{*, a} \rightarrow O_{i_a}$ we chose earlier.

We define $C(\mathcal{G})$ as the structure $(\mathbb{U}, O_*, M_*, R)$. 

Then, by construction, there is a 0-definable groupoid $\mathcal{G}'$ in $C(\mathcal{G})$, which is an extension of $\mathcal{G}$ with exactly one extra object in each isomorphism class. Note that $\mathcal{G}'$ is over $A$ as well, and that $C(\mathcal{G})$ comes with 0-definable maps $M_* \rightarrow A$ and $O_* \rightarrow A$, that yield the connected components of the new objects and morphisms.
\end{constr}

\begin{prop}\label{cover_in_independent_case}
The theory of $C(\mathcal{G})$ is a cover of the theory of $\mathbb{U}$.
\end{prop}
\begin{proof}
Let $\mathbb{V}' = (\mathbb{V}, S, M_*, M_* \rightarrow A, S \rightarrow A, R)$ be a saturated model of the theory of $C(\mathcal{G})$, such that $\mathbb{V}$ is elementarily equivalent to $\mathbb{U}$. We shall prove that $\mathbb{V}$ is stably embedded in $\mathbb{V}'$ using automorphisms.

Let $\sigma$ be an automorphism of $\mathbb{V}$. For each $a$ in $A$, we pick an index $i_a \in I_a$ and a morphism $f_a : S_a  \rightarrow O_{i_a}$.

We also pick a morphism $m_a : O_{i_{\sigma(a)}} \rightarrow O_{\sigma(i_a)}$. Such an $m_a$ exists because $\sigma(i_a)$ is in $I_{\sigma(a)}$, for $\sigma$ is an automorphism of $\mathbb{V}$, and because the connected components of $\mathcal{G}$ are given by the map $I \rightarrow A$. 

Now we shall define the action of $\sigma$ on $S$ and $M_*$.

On $S_a$, the function $\sigma$ shall act as the function $f_{\sigma(a)}^{-1} m_{a}^{-1}  \sigma f_a : S_a \rightarrow S_{\sigma(a)}$.

For the action on $M_*$, the point is to use the $f_a$ to get a morphism in $\mathcal{G}$, then apply $\sigma$ to it, then use morphisms again to get the appropriate domain and codomain.

\begin{itemize}
    \item If $m : O_j \rightarrow S_a$, we define $\sigma(m) = f_{\sigma(a)}^{-1} m_a^{-1}\sigma(f_a m) : O_{\sigma(j)} \rightarrow S_{\sigma(a)}$.
    \item If $m : S_a \rightarrow O_j$, we define $\sigma(m) = \sigma(m f_a^{-1}) m_a f_{\sigma(a)} : S_{\sigma(a)} \rightarrow O_{\sigma(j)}$.
    \item Finally, if $m : S_a \rightarrow S_a$, we define $\sigma(m) = f_{\sigma(a)}^{-1} m_a^{-1} \sigma(f_a m f_a^{-1}) m_a f_{\sigma(a)}  : S_{\sigma(a)} \rightarrow S_{\sigma(a)}$.
\end{itemize}

We can then check that $\sigma$ preserves the relation $R$, and commutes with the maps $M_*\rightarrow A$, $S \rightarrow A$.

\end{proof}

\begin{rem}
In the previous case, with only one isomorphism class, the morphism $m_a$ was given by the function defining $O_*$ as a copy of $O_i$. The morphism $n_a$ was given by the function $f_c^{-1}$. 
\end{rem}

\end{subsection}

\end{section}

\begin{section}{Functoriality in the independent case}

In this section, we write down the details of the proof of the equivalence of categories that were left to the reader at the end of section 5 of \cite{haykazyan-moosa}.

We find that all the ideas of the internal case extend naturally to the independent case. In fact, most of the time, dealing with independent fibers amounts to dealing with each fiber individually, without having to worry about coherence. 

\begin{subsection}{Categories of groupoids, categories of covers}

First, let us define what morphisms of groupoids are. Our definition is almost the same as that of \cite{haykazyan-moosa}(4.2), with only a small change to the condition (B) found there.

\begin{defi}
Let $\mathcal{G}_1$, $\mathcal{G}_2$ be 0-definable groupoids in $\mathbb{U}$. \begin{itemize}
    \item 
 A \textit{0-definable morphism} $H : \mathcal{G}_1 \rightarrow \mathcal{G}_2$ is a non-empty 0-definable family of definable functions between objects of $\mathcal{G}_1$ and objects of $\mathcal{G}_2$, as defined in \cite{haykazyan-moosa}(4.2), that satisfies the following conditions :

(A) : For each pair of maps $h_p : O_{1,i} \rightarrow O_{2,j}$, $h_m : O_{1,k} \rightarrow O_{2,l}$ in $H$, the sets of maps $h_p \circ \mathrm{Hom}_{\mathcal{G}_1}(O_{1,k}, O_{1,i})$ and $\mathrm{Hom}_{\mathcal{G}_2}(O_{2,l}, O_{2,j}) \circ h_m $ are equal.

(B') : The set of maps $H$ is stable under precomposition with morphisms of $\mathcal{G}_1$, and under postcomposition with morphisms of $\mathcal{G}_2$.

\smallskip

We recall the condition (B) given in \cite{haykazyan-moosa}(4.2) : 
(B) : The set of maps $H$ is stable under postcomposition with morphisms of $\mathcal{G}_2$.

Formally, a morphism of groupoids is given by a definable set of maps $H$ and a definable set $R \subseteq H \times O_1 \times O_2$, such that, for all $h \in H$, there exists a unique pair $(i_1, i_2) \in I_1 \times I_2$ such that the set $R_h \subset O_{1,i_1} \times O_{2,i_2}$ is the graph of a map from $(O_1)_{i_1}$ to $(O_2)_{i_2}$.

\item A morphism of groupoids is said to be \textit{injective} if all the maps it contains are injective.

\end{itemize}
\end{defi}

Note that the remark found at the end of the definition in \cite{haykazyan-moosa}(4.2) is incorrect : If we are given a groupoid morphism $H : \mathcal{G}_1 \rightarrow \mathcal{G}_2$ satisfying conditions (A) and (B'), and if we pick a nonempty proper 0-definable subset $X$ of the set of objects of $\mathcal{G}_1$, we can restrict $H$ to the set of functions whose domain is one of the objects in $X$. We get a nonempty set of functions satisfying (A) and (B), but it wouldn't satisfy (B').

We now prove that isomorphic connected groupoids are equivalent in the sense of \cite{hrushovski} (Section 3), which was not proved in full generality in \cite{haykazyan-moosa}.

\begin{prop}
Let $H : \mathcal{G}_1 \rightarrow \mathcal{G}_2$ be an isomorphism between 0-definable connected groupoids.
Then there exists a 0-definable connected groupoid $\mathcal{G}$, whose set of objects is the disjoint union of the objects of the $\mathcal{G}_i$, and 0-definable groupoid embeddings $f_i : \mathcal{G}_i \rightarrow \mathcal{G}$.
\end{prop}
\begin{proof}
The set of objects of $\mathcal{G}$ is defined as the disjoint union of $Ob(\mathcal{G}_1)$ and $Ob(\mathcal{G}_2)$. The morphisms are defined as follows : 

\begin{itemize}
    \item Between objects of $\mathcal{G}_1$, the morphisms are copies of morphisms of $\mathcal{G}_1$.
    \item Between objects of $\mathcal{G}_2$, the morphisms are copies of morphisms of $\mathcal{G}_2$.
    \item If $O_1$ is an object of $\mathcal{G}_1$ and $O_2$ is an object of $\mathcal{G}_2$, the set of morphisms from $O_1$ to $O_2$ is the set of maps in the morphism $H$ that are bijections from $O_1$ to $O_2$.
    \item Conversely, the set of maps from $O_2$ to $O_1$ is the set of inverses of maps in $H$ that are bijections from $O_1$ to $O_2$.
\end{itemize}

Composition is defined as expected, using composition of maps. 

Using condition (A), we check that if $h : O_2 \rightarrow O_1$ and $g: O'_1 \rightarrow O_2$ are morphisms in $\mathcal{G}$, then $h \circ g$ is a (copy of a) morphism in $\mathcal{G}_1$, and thus is a morphism of $\mathcal{G}$. Symmetrically, if $h : O_2 \rightarrow O_1$ and $g: O_1 \rightarrow O'_2$ are morphisms in $\mathcal{G}$, then $ g \circ h$ is a (copy of a) morphism in $\mathcal{G}_2$.

Now, condition (B) shows that we may postcompose morphisms between objects of $\mathcal{G}_1$ and objects of $\mathcal{G}_2$ with morphisms in $\mathcal{G}_2$. 

Using conditions (A) and (B'), we show that we can also precompose with morphisms in $\mathcal{G}_1$ : Let $h : O_1 \rightarrow O_2 $ be a function in $H$, let $\alpha : O'_1 \rightarrow O_1$ be a morphism in $\mathcal{G}_1$. Let $h' \in H$ be a map from $O'_1$ to $O_2$. Such a map exists by condition (B'). Using condition (A), there exists a morphism $\beta : O_2 \rightarrow O_2$ in the groupoid $\mathcal{G}_2$ such that $\beta \circ h' = h \circ \alpha$. By condition (B'), the function $\beta \circ h'$ is still in $H$, and thus so is $h \circ \alpha$.

We can also precompose inverses of maps in $H$ with maps in $\mathcal{G}_2$ and postcompose inverses of maps in $H$ with maps in $\mathcal{G}_1$.

Finally, this groupoid $\mathcal{G}$ is 0-definable, since both the $\mathcal{G}_i$ and $H$ are. The embeddings $f_i : \mathcal{G}_i \rightarrow \mathcal{G}$ are the obvious ones. 
\end{proof}

Now, if two groupoids are not connected, we want to specify which connected components are \textquotedblleft linked together" by the morphisms of groupoids between them.

\begin{defi}\label{morphisms_of_non_connected_groupoids}
 Let $A, B$ be 0-definable sets. Let $R \subseteq  A \times B$ be a definable relation. Let $\mathcal{G}_1, \mathcal{G}_2$ be groupoids over $A$ and $B$ respectively. A morphism of groupoids $H : \mathcal{G}_1 \rightarrow \mathcal{G}_2$ is said to be \textit{compatible with the relation $R$} if, for all $a$ in $A$, for all $b$ in $B$, there exists a map in $H$ from some object of $\mathcal{G}_{1,a}$ to some object of $\mathcal{G}_{2,b}$ if and only if $\models R(a,b)$.

If $A=B$, unless otherwise stated, we shall consider only morphisms of groupoids that are compatible with the equality relation. 
\end{defi}

\end{subsection}

In the rest of this section, we let $\mathrm{FCG_A}$ denote the category of finitely faithful 0-definable groupoids over $A$, whose morphisms are defined as above. It is indeed a category. For these groupoids, the set of connected components is (definably) indexed by $A$.

Let us now define morphisms of covers. For the various proofs, we refer the reader to \cite{haykazyan-moosa}.

\begin{defi}
Let $\mathbb{V}_1, \mathbb{V}_2$ be (models of) covers of $\mathbb{U}$, such that $\mathbb{V}_i$ has one extra sort $S_i$, compared to $\mathbb{U}$, for $i= 1,2$. A \textit{morphism of covers} $H : \mathbb{V}_1 \rightarrow \mathbb{V}_2$ is the theory of a structure on $S_1 \sqcup S_2 \cup \mathbb{U}$ which extends the structures on the $\mathbb{V}_i$, comes with a distinguished 0-definable map $h : S_1 \rightarrow S_2$, and defines a cover of both $Th(\mathbb{V}_1)$ and $Th(\mathbb{V}_2)$.
\end{defi}

Note that two distinct maps $h, h' : S_1 \rightarrow S_2$ may define the same theory, and thus the same morphism of covers.
Despite this, we may denote a morphism of covers with the same symbol as some 0-definable map which defines it.

If $h : \mathbb{V}_1 \rightarrow \mathbb{V}_2$, $g : \mathbb{V}_2 \rightarrow \mathbb{V}_3$ are morphisms of covers, we define the composite $g \circ h : \mathbb{V}_1 \rightarrow \mathbb{V}_3$ by taking the structure induced on $S_1 \cup S_3 \cup \mathbb{U}$ by the full structure on $S_1 \cup S_2 \cup S_3 \cup \mathbb{U}$, and by choosing $g \circ h : S_1 \rightarrow S_3$ as the distinguished 0-definable map.

\begin{prop}
Composition of morphisms of covers is well-defined and associative. This composition defines the category of the covers of $\mathbb{U}$ that only have one extra sort.
\end{prop}

\begin{proof}
See \cite{haykazyan-moosa}, Remark 4.1, Proposition 4.2, Definition 4.3. See also Lemma 2.4, and the paragraph preceding it, in the same paper, for more details on unions of covers. 
\end{proof}

Now, let $\mathrm{ACIF_A}$ be the category of 1-analysable covers over $A$ with independent fibers and finitely generated languages, defined just as in \cite{haykazyan-moosa}(between Proposition 5.7 and Example 5.8). The morphisms are morphisms of covers compatible with the maps $S \rightarrow A$. The same arguments as in \cite{haykazyan-moosa}(4.1) show that composition of morphisms is well-defined and yields a category. In fact, the proof there deals with the more general case of covers with one extra sort, so it can be applied to our specific case of 1-analysable ones.

\begin{subsection}{The functor $C : \mathrm{FCG_A} \rightarrow \mathrm{ACIF_A}$}
\begin{defi}
Let $h : \mathcal{G}_1 \rightarrow \mathcal{G}_2$ be a morphism of groupoids. Let $\mathcal{G}'_1$ and $\mathcal{G}'_2$ be the extensions of the $\mathcal{G}_i$ defined in \cite{haykazyan-moosa}(5.6). For each $a$ in $A$, let $h_a : O_{1, i_a} \rightarrow O_{2, j_a}$ be one of the maps in the morphism $h$, and let $f_a : S_{1,a} \rightarrow O_{1, i_a} $, \, $g_a : O_{2, j_a} \rightarrow S_{2, a}$ be morphisms in $\mathcal{G}'_1$ and $\mathcal{G}'_2$ respectively.

The morphism $\mathrm{C}(h)$ is defined as the theory of  $(\mathrm{C}(\mathcal{G}_1)\cup \mathrm{C}(\mathcal{G}_2), \bigcup\limits_a \, g_a \circ h_a\circ f_a)$.
\end{defi}

\begin{prop}\label{C(h)_free_choice}
The theory of $\mathrm{C}(h)$ depends only on the morphism $h$, and not on the choices made. More precisely, different choices of maps lead to structures that are isomorphic over $\mathbb{U}$. 
\end{prop}
\begin{proof}

Let $f'_a, g'_a, h'_a, i'_a, j'_a$ be other possible choices of maps for the index $a$. Using the connectedness of the groupoid, and preservation of the set of morphisms in $h$ under pre and post composition by groupoid morphisms, we may assume that only $h'_a$ is different from $h_a$, all else being equal. This is the same as assuming that the following diagram commutes, and replacing $h'_a$ with $h''_a.$

\begin{center}
\begin{tikzcd}S_{1,a} \arrow[r, "f_a"] \arrow[rd, "f'_a"] & O_{1, i_a}\arrow[r, "h''_a"] \arrow[d]& O_{2,j_a} \arrow[d]\arrow[r, "g_a"] & S_{2,a}\\
& O_{1, i'_a}\arrow[r, "h'_a"]& O_{2,j'_a} \arrow[ru, "g'_a"]  & \end{tikzcd}
\end{center}

Now, using property $(A)$ of the morphism $h$, we get : $h''_a \circ \sigma_a = h_a$, for some $\sigma_a \in \mathrm{Aut}_{\mathcal{G}_1}(O_{1, i_a})$.

Thus, we have $f_a^{-1} \circ \sigma_a \circ f_a \in \mathrm{Aut}_{\mathcal{G}_1}(S_{1, a}) = \mathrm{Aut}(S_{1, a}/\mathbb{U})$
.

By independence of the fibers, we can then construct an automorphism $\rho$ of the structure $\mathrm{C}(\mathcal{G}_1)\cup \mathrm{C}(\mathcal{G}_2)$ that fixes $\mathbb{U}$ pointwise, given by the union of the $f_a^{-1} \circ \sigma_a \circ f_a$ and $id_{\mathrm{C}(\mathcal{G}_2)}$. Finally, we check that $g_a h''_a f_a \rho = g_a h''_a f_a f_a^{-1} \sigma_a f_a = g_a h''_a \sigma_a f_a = g_a h_a f_a = \rho g_a h_a f_a $. Thus, the structures given by $h_a, f_a, g_a$ and $h''_a, f_a, g_a$ are isomorphic over $\mathbb{U}$.

\end{proof}

\begin{lemma}\label{C(h)_is_cover_of_U}
The theory of $\mathrm{C}(h)$ is a cover of the theory of $\mathbb{U}$.
\end{lemma}

\begin{proof}
We assume $\mathrm{C}(h)$ to be saturated. We shall prove that any automorphism of $\mathbb{U}$ can be lifted into an automorphism of $\mathrm{C}(h)$. Let $\sigma \in \mathrm{Aut}(\mathbb{U}).$ Let $\sigma' \in \mathrm{Aut}(\mathrm{C}(\mathcal{G}_1)\cup \mathrm{C}(\mathcal{G}_2))$ be an extension of $\sigma$. 

Now, we shall use the proof of the previous proposition : we notice that the maps $\sigma'(f_a), \sigma'(g_a), \sigma(h_a)$ represent another possible choice for the construction of $\mathrm{C}(h)$. Indeed, since $\sigma'$ is an element of $\mathrm{Aut}(\mathrm{C}(\mathcal{G}_1)\cup \mathrm{C}(\mathcal{G}_2))$, it preserves the structure of the 0-definable groupoids $\mathcal{G}'_1$ and $\mathcal{G}'_2$, and the 0-definable set of maps given by $h$.
Therefore, there exists some automorphism $\tau$ of $\mathrm{C}(\mathcal{G}_1)\cup \mathrm{C}(\mathcal{G}_2)$ fixing $\mathbb{U}$ pointwise such that, for all $a$ in $A$, we have 
$\tau \sigma'(g_a) \sigma(h_a) \sigma'(f_a) = g_{\sigma(a)} h_{\sigma(a)} f_{\sigma(a)} \tau$. Let $\rho = \tau \sigma'$. We check that this $\rho$ commutes with the map $\bigcup\limits_{a\in A} g_a h_a f_a$. Let $a$ be an element of $A$. We compute : $$\rho g_a h_a f_a = \tau \sigma' g_a h_a f_a = \tau \circ [\sigma'(g_a h_a f_a)] \circ \sigma' = g_{\sigma(a)} h_{\sigma(a)} f_{\sigma(a)} \tau \sigma' = g_{\sigma(a)} h_{\sigma(a)} f_{\sigma(a)} \rho.$$

Finally, since $\tau$ fixes $\mathbb{U}$ pointwise, the function $\rho$ extends $\sigma$.
\end{proof}

\begin{prop}\label{C(h)_is_morphism_of_covers}
The theory of $C(h)$ is a cover of both $\mathrm{C}(\mathcal{G}_1)$ and $\mathrm{C}(\mathcal{G}_2)$.
\end{prop}

\begin{proof}
Let $\sigma_1$ be an automorphism of $\mathrm{C}(\mathcal{G}_1)$. By the previous lemma, we may assume that $\sigma_1$ fixes $\mathbb{U}$ pointwise. We need to define the action of $\sigma_1$ on each of the $S_{2,a}$. Note that since $\sigma_1$ fixes $\mathbb{U}$ pointwise, an extension of $\sigma_1$ should fix each fiber $S_{2,a}$ setwise. Now, let $a$ be an element of $A$. We notice that the map $f_a \sigma_1|_{S_{1,a}} f_a^{-1}$ is a morphism in $\mathcal{G}_1$. Thus, by condition (A) of the morphism $h$, there exists some morphism $\gamma$ in $\mathcal{G}_2$ such that $h_a f_a \sigma_1|_{S_{1,a}} f_a^{-1} = \gamma h_a $. Now, since $\mathcal{G}'_2$, is a groupoid, there exists a morphism $\sigma_{2,a} \in \mathrm{Aut}_{\mathcal{G}'_2}(S_{2,a})$ such that $g_a\gamma = \sigma_{2,a} g_a$.

\begin{center}
\begin{tikzcd}
S_{1,a} \arrow[r, "f_a"] & O_{1,i} \arrow[r, "h_a"] & O_{2,j}\arrow[r, "g_a"] & S_{2,a} \\
S_{1,a} \arrow[u, "\sigma_1"]& O_{1,i}\arrow[l, "f_a^{-1}"]  \arrow[r, "h_a"]& O_{2,j}\arrow[u, dotted, "\exists \gamma"]  \arrow[r, "g_a"]& S_{2,a} \arrow[u, dotted, "\exists \sigma_{2,a}"] 

\end{tikzcd}
\end{center}

Finally, we have $g_a h_a f_a \sigma_1|_{S_{1,a}} = g_a\gamma h_a f_a = \sigma_{2,a} g_a h_a f_a.$ We then use the family of the $\sigma_{2,a}$ and the independence of the fibers in $\mathrm{C}(\mathcal{G}_2)$ to define an automorphism $\sigma_2$ of $\mathrm{C}(\mathcal{G}_2)$ that fixes $\mathbb{U}$ pointwise and which agrees with the $\sigma_{2,a}$.
Thus, the map $\sigma_1 \cup \sigma_2$ is an automorphism of $(\mathrm{C}(\mathcal{G}_1)\cup \mathrm{C}(\mathcal{G}_2), \bigcup\limits_a \, g_a \circ h_a\circ f_a)$ that extends $\sigma_1$.

Similarly, $\mathrm{C}(\mathcal{G}_2)$ is stably embedded in $(\mathrm{C}(\mathcal{G}_1)\cup \mathrm{C}(\mathcal{G}_2), \bigcup\limits_a \, g_a \circ h_a\circ f_a)$.
\end{proof}

\begin{prop}\label{C_is_functor_independent_case}
The action of $C$ on the morphisms defines a functor $C :\mathrm{FCG_A} \rightarrow \mathrm{ACIF_A}$.
\end{prop}

\begin{proof}
We need to show that $C$ sends identity morphisms of groupoids to identity morphisms of covers, and that it preserves composition.

For the first part, using proposition \ref{C(h)_free_choice}, we may assume that each $h_a$ is the identity morphism of some object $O_a$ in the groupoid $\mathcal{G}$. We may also assume that the maps $g_a$ and $f_a$ are inverses of one another. Thus, the map $\bigcup\limits_a \, g_a \circ h_a \circ f_a$ is the identity of the new sort $S$. Therefore, the cover $C(id_{\mathcal{G}})$ is indeed the identity morphism of $C(\mathcal{G})$.

Regarding composition, we use proposition \ref{C(h)_free_choice} again.
Let $h_{12} : \mathcal{G}_1 \rightarrow \mathcal{G}_2$, $h_{23} : \mathcal{G}_2 \rightarrow \mathcal{G}_3$ be morphisms. Let $h_{13} : \mathcal{G}_1 \rightarrow \mathcal{G}_3$ be the morphism $h_{23} \circ h_{12}$. Then we may assume that the morphisms picked to define $C(h_{12})$ and $C(h_{23})$ are consistent with the ones picked to define $C(h_{13})$. More precisely, we may assume that, for each $a \in A$, the morphisms $g_a$, $h_a$ and $f_a$ were picked to get the equality : $$g^{23}_a h^{23}_a f^{23}_a g^{12}_a h^{12}_a f^{12}_a = g^{13}_a h^{13}_a f^{13}_a$$

Indeed, the map $f_a^{23}g_a^{12} $ is a morphism in the groupoid $\mathcal{G}_2$. Therefore, by property (B) of the groupoid morphism $h_{12}$, the map $ f_a^{23}g_a^{12} h_a^{12}$ is again a map in the morphism $h_{12}$. Thus, the map $h_a^{23} f_a^{23} g_a^{12}  h_a^{12}$ is a map of the groupoid morphism $h_{23} \circ h_{12} = h_{13}$.

\begin{center}
\begin{tikzcd}S_{1,a} \arrow[r, "f^{12}_a"]\arrow[rd, "f^{13}_a"] & O_{1, i_a}\arrow[r, "h^{12}_a"]& O_{2,j_a} \arrow[r, "g^{12}_a"] & S_{2,a} \arrow[r, "f^{23}_a"] & O_{2, k_a}\arrow[r, "h^{23}_a"]& O_{3,l_a} \arrow[r, "g^{23}_a"] & S_{3,a} \\ & O_{1, i_a} \arrow[u, "id" description] \arrow[rrrr, "h^{13}_a"] &&&& O_{3, l_a}\arrow[u, "id" description]\arrow[ur, "g^{13}_a"] \end{tikzcd}
\end{center}

By proposition \ref{C(h)_free_choice}, we may assume that $h_a^{13}=h_a^{23} f_a^{23} g_a^{12}  h_a^{12}$ and $f_a^{12} = f_a^{13}$ and $g_a^{23} = g_a^{13}$. In such a case, the desired equality holds, i.e. the above diagram commutes.
Thus, from the definition of the composition of morphisms of covers, we get $$C(h_{23} \circ h_{12})=C(h_{23})\circ C(h_{12}).$$
\end{proof}

\end{subsection}

\begin{subsection}{The functor $G :\mathrm{ACIF_A} \rightarrow  \mathrm{FCG_A}$}
\begin{defi}
Let $(\mathbb{U}_1 \cup \mathbb{U}_2, h : S_1 \rightarrow S_2)$ be a morphism in the category $\mathrm{ACIF_A}$. We assume the structures to be saturated. Let $\mathcal{G}_1$, $\mathcal{G}_2$ be the liaison groupoids $\mathcal{G}(\mathbb{U}_1)$, $\mathcal{G}(\mathbb{U}_2)$. Let $\mathcal{G}'_1$, $\mathcal{G}'_2$ be the extensions of these groupoids that are 0-definable in $\mathbb{U}_1^{eq}$ and $\mathbb{U}_2^{eq}$ respectively. The groupoid morphism $G(h)$ is defined as the set of morphisms $H:=\lbrace g_a h_a f_a : O_{1, i_a} \rightarrow O_{2, j_a} | \, a \in A, \, i_a \in Ob(\mathcal{G}_{1,a}), j_a \in Ob(\mathcal{G}_{2,a}), f_a \in Hom_{\mathcal{G}'_1}(O_{1,i_a}, S_{1,a}) , g_a \in Hom_{\mathcal{G}'_2}(S_{2,a}, O_{2,j_a}) \rbrace$.

This set $H$ of morphisms is a 0-definable morphism of groupoids over $A$ in $(\mathbb{U}_1 \cup \mathbb{U}_2, h : S_1 \rightarrow S_2)^{eq}$. We wish to get a set of morphisms 0-definable in $\mathbb{U}$ with a 0-definable map to $A$. A compactness argument similar to the one in \cite{haykazyan-moosa}(5.6, §5), combined with the hypothesis of stable embeddedness, yields the result.

\end{defi}

\begin{prop}
The set of maps $G(h)$ depends only on the theory of $(\mathbb{U}_1 \cup \mathbb{U}_2, h)$, and is a morphism in $\mathrm{FCG}_A$.
\end{prop}
\begin{proof}
Let $h' : S_1 \rightarrow S_2$ such that $(\mathbb{U}_1 \cup \mathbb{U}_2, h) \equiv (\mathbb{U}_1 \cup \mathbb{U}_2, h')$. By saturation, there exists an isomorphism $\sigma : (\mathbb{U}_1 \cup \mathbb{U}_2, h) \simeq(\mathbb{U}_1 \cup \mathbb{U}_2, h')$. By stable embeddedness of $\mathbb{U}_1$ in $(\mathbb{U}_1 \cup \mathbb{U}_2, h)$, we may assume $\sigma_{\mathbb{U}_1} = id$. Thus, if $\tau = \sigma|_{\mathbb{U}_2} \in \mathrm{Aut}(\mathbb{U}_2/\mathbb{U}) = \prod_a \mathrm{Aut_{\mathcal{G}'_{2,a}}(S_{2,a})}$, we get $\tau h = h'$. Thus the families of maps induced by $h$ and $h'$ are equal.

We now need to check conditions (A) and (B'). The latter is easier to see, since postcomposing a map in $  Hom_{\mathcal{G}'_2}(S_{2,a}, O_{2,j_a}) \circ h_a \circ Hom_{\mathcal{G}'_1}(O_{1,i_a}, S_{1,a})$ by a morphism in $Hom_{\mathcal{G}_2} (O_{2,j_a}, O_{2, k_a})$ yields a map in $Hom_{\mathcal{G}'_2}(S_{2,a}, O_{2,k_a}) \circ h_a \circ Hom_{\mathcal{G}'_1}(O_{1,i_a}, S_{1,a}).$ Similar arguments work for precomposition.

Let us now check condition (A). Let $\alpha_1, \beta_1$ be morphisms in $\mathcal{G}'_1$ with target $S_{1,a}$, let $\alpha_2, \beta_2$ be morphisms in $\mathcal{G}'_2$ with source $S_{2,a}$.
We shall prove that, for each morphism $\alpha$ in $\mathcal{G}_1$ such that the composition $\alpha_2 h_a \alpha_1 \alpha$ is well-defined and has the same source as $\beta_2 h_a \beta_1$, there exists a morphism $\beta$ in $\mathcal{G}_2$ such that $\alpha_2 h_a \alpha_1 \alpha = \beta \beta_2 h_a \beta_1$. A symmetrical argument will then give the converse. Let such a morphism $\alpha$ be. Then, the morphism $\alpha_1 \alpha \beta_1^{-1}$ is in $\mathrm{Aut}_{\mathcal{G}'_1}(S_{1,a}) = \mathrm{Aut}(S_{1,a} / \mathbb{U})$. By stable embeddedness and independence of fibers, the morphism $\alpha_1 \alpha \beta_1^{-1}$ can be extended to an automorphism $\sigma \in \mathrm{Aut}(\mathbb{U}_1 / \mathbb{U})$. By stable embdeddness again, this automorphism can be extended to an automorphism $\tau$ of $\mathbb{U}_1 \cup\mathbb{U}_2, h$. We then have $\alpha_2 h_a \alpha_1 \alpha \beta_1^{-1} = \alpha_2 h_a \tau = \alpha_2 \tau h_a$. Here, $\alpha_2 \tau$ is a morphism of the connected groupoid $\mathcal{G}'_{2,a}$, and so is $\beta_2$. Moreover, these morphisms have the same source, which is $S_{2,a}$. Thus, there exists a morphism $\beta$ in $\mathcal{G}'_{2,a}$ such that $\beta \beta_2 = \alpha_2 \tau$. Thus $\alpha_2 h_a \alpha_1 \alpha = \beta \beta_2 h_a \beta_1$. Since $\beta$ is defined between objects of $\mathcal{G}_{2,a}$, it is in fact a morphism in $\mathcal{G}_{2,a}$.
\end{proof}

\begin{prop}
The action of $G$ on the morphisms of $ACIF_A$ defines a functor.
\end{prop}
\begin{proof}
The action on the identity morphisms is pretty clear : since we may pick any map $h : S_1 \rightarrow S_2 $, as long as it yields the desired theory, we may assume $h = id_{S_1}$. In that case, the set of maps in $G(h)$ is just the set of morphisms in $\mathcal{G}_1$, which is indeed the identity of the groupoid.

For composition, we use the same argument as in \cite{haykazyan-moosa}(4.11). Let $(\mathbb{U}_1 \cup \mathbb{U}_2, h_{12})$, $(\mathbb{U}_2 \cup \mathbb{U}_3, h_{23})$ be morphisms in $ACIF_A$.
The maps in $G(h_{23})\circ G(h_{12})$ are of the form $\alpha_3 h_{23, a} \alpha_2 \beta_2 h_{12, a} \alpha_1$, whereas those in $G(h_{23}\circ h_{12})$ are of the form $\beta_3 h_{23, a} h_{12, a} \beta_1$. Using the identity of $S_{2,a}$, which is indeed a morphism in $\mathcal{G}'_2$, we deduce that $G(h_{23}\circ h_{12}) \subseteq G(h_{23})\circ G(h_{12})$. For the other inclusion, let 
$\alpha_3 h_{23, a} \alpha_2 \beta_2 h_{12, a} \alpha_1$ be a map in $G(h_{23})\circ G(h_{12})$. Then, the map $\alpha_2 \beta_2$ is an automorphism of $S_{2,a}$ in $\mathcal{G}'_2$. Thus, ny stable embeddedness and independence of fibers, it can be extended by identity to an automorphism of $\mathbb{U}_2$ fixing $\mathbb{U}$ pointwise. By stable embeddedness again, this automorphism can be extended to an automorphism $\sigma$ of $(\mathbb{U}_1 \cup \mathbb{U}_2, h_{12})$. We then compute $\alpha_3 h_{23, a} \alpha_2 \beta_2 h_{12, a} \alpha_1= \alpha_3 h_{23, a} \sigma h_{12, a} \alpha_1=\alpha_3 h_{23, a} h_{12, a} \sigma\alpha_1$. Here, $\sigma\alpha_1$ is a morphism in $\mathcal{G}'_1$. This proves the equality $G(h_{23}\circ h_{12}) = G(h_{23})\circ G(h_{12})$.

\end{proof}

\end{subsection}

\begin{subsection}{The isomorphism $\eta : id_{ACIF_A} \simeq CG$}
The following proposition is a relativization of proposition 4.2 in \cite{haykazyan-moosa}, under an extra assumption.

\begin{prop}
Let $\mathbb{U}_1, \mathbb{U}_2$ be objects in the category $\mathrm{ACIF_A}$. Let $h : S_1 \rightarrow S_2$ be a map compatible with the maps $S_1 \rightarrow A$, $S_2 \rightarrow A$. Assume that $\mathbb{U}$ is stably embedded in $(\mathbb{U}_1 \cup \mathbb{U}_2, h)$. Then $(\mathbb{U}_1 \cup \mathbb{U}_2, h)$ is a morphism in $\mathrm{ACIF_A}$ if and only if, for all $a \in A$, the map $h_a \cup id_{\mathbb{U}} : (\mathbb{U}, S_{1,a}) \rightarrow (\mathbb{U}, S_{2,a})$ yields a morphism of covers from $(\mathbb{U}, S_{1,a})$ to $(\mathbb{U}, S_{2,a})$.
\end{prop}

\begin{proof}
We shall prove that $\mathbb{U}_1$ is stably embedded in $(\mathbb{U}_1 \cup \mathbb{U}_2, h)$. Let $\sigma$ be an automorphism of $\mathbb{U}_1$. We want to extend it to an automorphism of $(\mathbb{U}_1 \cup \mathbb{U}_2, h)$. We may assume that $\sigma$ fixes $\mathbb{U}$ pointwise. Thus, for each $a$ in $A$, $\sigma$ restricts into an automorphism $\sigma_a$ of $(\mathbb{U}, S_{1,a})$. Since $(\mathbb{U}, S_{1,a}, S_{2,a}, h_a)$ is a cover of $(\mathbb{U}, S_{1,a})$, there exists an automorphism $\tau_a \in \mathrm{Aut}((\mathbb{U}, S_{1,a}, S_{2,a}, h_a)/\mathbb{U})$ extending $\sigma_a$. Using independence of fibers and stable embeddedness of the $(\mathbb{U}, S_{2,a})$ in $(\mathbb{U}, S_2)$, the union of maps $\tau := \bigcup\limits_a \tau_a$ is an automorphism of $\mathbb{U}_1 \cup \mathbb{U}_2$ extending $\sigma$. By construction, it commutes with $h$. Thus $\mathbb{U}_1$ is stably embedded in $(\mathbb{U}_1 \cup \mathbb{U}_2, h)$, as desired.
\end{proof}

\begin{coro}
Let $\mathbb{U}_1, \mathbb{U}_2$ be objects in the category $\mathrm{ACIF_A}$. Let $h : S_1 \rightarrow S_2$ be a map compatible with the maps $S_1 \rightarrow A$, $S_2 \rightarrow A$. Assume that $\mathbb{U}$ is stably embedded in $(\mathbb{U}_1 \cup \mathbb{U}_2, h)$. 
Then $(\mathbb{U}_1 \cup \mathbb{U}_2, h)$ is an isomorphism of covers if and only if, for each $a$ in $A$, $(\mathbb{U}, S_{1,a}, S_{2,a}, h_a)$ is an isomorphism of covers.
\end{coro}

\begin{coro}
Let $\mathbb{U}_1, \mathbb{U}_2$ be objects in the category $\mathrm{ACIF_A}$. Let $h : S_1 \rightarrow S_2$ be a map compatible with the maps $S_1 \rightarrow A$, $S_2 \rightarrow A$. Assume that $\mathbb{U}$ is stably embedded in $(\mathbb{U}_1 \cup \mathbb{U}_2, h)$. Then $(\mathbb{U}_1 \cup \mathbb{U}_2, h)$ is an isomorphism in $\mathrm{ACIF_A}$ if and only if, for all $a \in A$, the map $h_a \cup id_{\mathbb{U}} : (\mathbb{U}, S_{1,a}) \rightarrow (\mathbb{U}, S_{2,a})$ is a bijective bi-interpretation.
\end{coro}

\begin{defi}
Let $\mathbb{V}=(\mathbb{U}, S, S \rightarrow A)$ be an object of $ACIF_A$. Let $\mathcal{G}'(\mathbb{V})$ be the extension of the liaison groupoid $\mathcal{G}(\mathbb{V})$ with one extra object $S_a$ in each isomorphism class, such that $\mathrm{Aut}_{\mathcal{G}'(\mathbb{V})}(S_a)=\mathrm{Aut}(S_a/\mathbb{U})$. Let $(\mathbb{U}, O_*)$ be the cover $C(\mathcal{G}(\mathbb{V}))$. Let $\mathcal{G}''(\mathbb{V})$ be the extension of $\mathcal{G}(\mathbb{V})$ defined in $C(\mathcal{G}(\mathbb{V}))^{eq}$ by formally creating copies of objects and morphisms of $\mathcal{G}(\mathbb{V})$.

Then $\eta_{\mathbb{V}}$ is defined as $(\mathbb{V}\cup C\mathcal{G}(\mathbb{V}), \bigcup\limits_a h_a g_a)$, where, for each $a$ in $A$, $i_a$ is some object in $\mathcal{G}(\mathbb{V})_a$, and the maps $h_a$ and $g_a$ are arbitrary morphisms in $Hom_{\mathcal{G}''(\mathbb{V})}(O_{1, i_a}, O_{*,a})$ and $Hom_{\mathcal{G}'(\mathbb{V})}(S_a, O_{1,i_a})$ respectively.

\begin{center}
\begin{tikzcd}S_{a} \arrow[r, "g_a"]  & O_{1, i_a}\arrow[r, "h_a"] & O_{*,a} \end{tikzcd}
\end{center}

\end{defi}

\begin{lemma}\label{eta}
The theory of $\eta_{\mathbb{V}}$ is a cover of $\mathbb{U}$.
\end{lemma}

\begin{proof}
Let $\sigma$ be an automorphism of $\mathbb{U}$. Let $\sigma'$ be an automorphism of $(\mathbb{V}\cup C\mathcal{G}(\mathbb{V}))^{eq}$ that extends $\sigma$. 
We shall proceed as in lemma \ref{C(h)_free_choice}, by exhibiting another possible choice of family of morphisms, and finding an automorphism connecting the two choices.
For each $a$ in $A$, let $h'_a=\sigma'(h_{\sigma^{-1}(a)})$ and $g'_a=\sigma'(g_{\sigma^{-1}(a)})$. Since $\sigma'$ is an automorphism, $h'_a$ is a morphism in the groupoid $\mathcal{G}''(\mathbb{V})$ with target $O_{*,a}$ and source in the objects of $\mathcal{G}(\mathbb{V})$. Similarly, $g'_a$ is a morphism in the groupoid $\mathcal{G}'(\mathbb{V})$ with source $S_{*,a}$ and target in the objects of $\mathcal{G}(\mathbb{V})$, and the composition $h'_a g'_a$ is well-defined. Now, there exist morphisms $\alpha_a$ and $\beta_a$ in $\mathcal{G}(\mathbb{V})$ such that $h'_a = h_a \alpha_a$ and $g'_a = \beta_a g_a $. Again, there exists a morphism $\tau_a$ in $\mathrm{Aut}_{\mathcal{G}'(\mathbb{V})}(S_a)=\mathrm{Aut}(S_a / \mathbb{U})$ such that $g_a \tau_a = \alpha_a \beta_a g_a$. So we have $h_a g_a \tau_a = h_a \alpha_a \beta_a g_a = h'_a g'_a$.

\begin{center}
\begin{tikzcd}S_{a} \arrow[r, "g_a"] \arrow[rd, "g'_a"] & O_{i} \arrow[d, "\beta_a"] &  & O_{i}\arrow[r, "h_a"]  & O_{*,a} & S_a\arrow[r, "g_a"] & O_i \\
& O_{j}  && O_{j}\arrow[u, "\alpha_a"]\arrow[ru, "h'_a"]  && S_a \arrow[u, "\tau_a"] \arrow[r, "g_a"]& O_i \arrow[u, "\alpha_a \beta_a"] \end{tikzcd}
\end{center}

Now, using independence of fibers, the morphisms $\tau_a$ define an automorphism of $\mathbb{V}$ fixing $\mathbb{U}$ pointwise, which we extend by identity on $O_*$.This way, we define an automorphism $\tau$ of $(\mathbb{V}\cup C\mathcal{G}(\mathbb{V}))$ such that, for all $a$ in $A$, we have $h_a g_a \tau \sigma'= \tau \sigma'(h_{\sigma^{-1}(a)}) \sigma'(g_{\sigma^{-1}(a)}) \sigma' = \tau \sigma' h_{\sigma^{-1}(a)} g_{\sigma^{-1}(a)}$. The automorphism $\tau \sigma'$ has the desired properties.

\end{proof}

\begin{prop}
The theory of $\eta_{\mathbb{V}}$ is an isomorphism between $\mathbb{V}$ and $C(\mathcal{G}(\mathbb{V}))$.
\end{prop}

\begin{proof}
We will prove that $\mathbb{V}$ is stably embedded in $\eta_{\mathbb{V}} = (\mathbb{V}\cup C(\mathcal{G}(\mathbb{V})), \bigcup\limits_a h_a g_a).$ Similar arguments show that $C(\mathcal{G}(\mathbb{V}))$ is stably embedded in $\eta_{\mathbb{V}} = (\mathbb{V}\cup C(\mathcal{G}(\mathbb{V}))$. Then, $\eta_{\mathbb{V}}$ will be a morphism of covers. The map $\bigcup\limits_a h_a g_a$ being bijective, the morphism will thus be an isomorphism.

Let $\sigma$ be an automorphism of $\mathbb{V}$. By the previous lemma, we may assume $\sigma$ to fix $\mathbb{U}$ pointwise. We only have to extend $\sigma$ to the sort $O_*$ in a way that make it commute with the map $\bigcup\limits_a h_a g_a$ and be an automorphism of $C(\mathcal{G}(\mathbb{V}))$. Luckily, since the map $\bigcup\limits_a h_a g_a$ is a bijection, there is only one way to do so. Since $\sigma$ fixes $\mathbb{U}$ pointwise, we may consider $\sigma_a=\sigma|_{S_a}$, for $a$ in $A$. The map $\sigma_a$ is an element of $\mathrm{Aut}_{\mathbb{V}}(S_a / \mathbb{U})= \mathrm{Aut}_{\mathcal{G}'(\mathbb{V})}(S_a)$. Therefore, the map $g_a \sigma_a g_a^{-1}$ is a morphism in $\mathcal{G}(\mathbb{V})$, so the map $h_a g_a \sigma_a g_a^{-1} h_a^{-1}$ is an automorphism of $O_{*,a}$ in $\mathcal{G}''(\mathbb{V})$. It can therefore be extended into an automorphism of $(\mathbb{U}, O_{*,a})$ fixing $\mathbb{U}$ pointwise. By stable embeddedness, this can be extended into an automorphism of $C(\mathcal{G}(\mathbb{V}))$. Finally, by independence of fibers, the map $\bigcup\limits_a h_a g_a \sigma_a g_a^{-1} h_a^{-1} \cup \sigma_S \cup id_{\mathbb{U}}$ is an automorphism of $\mathbb{V} \cup C(\mathcal{G}(\mathbb{V}))$ that commutes with $\bigcup\limits_a h_a g_a$.
\end{proof}

\begin{prop}
The family of maps $\eta$ defines a natural isomorphism between the identity functor of $ACIF_A$ and the functor $CG$.
\end{prop}
\begin{proof}
We have to prove that the naturality squares commute.

\begin{center}
\begin{tikzcd} \mathbb{U}_{1}\arrow[r, "{\eta_{\mathbb{U}_1}}"] \arrow[d, "h"]& C\mathcal{G}(\mathbb{U}_1) \arrow[d, "C\mathcal{G}(h)"] \\
 \mathbb{U}_{2}\arrow[r, "{\eta_{\mathbb{U}_2}}"]& C\mathcal{G}(\mathbb{U}_2)   \end{tikzcd}
\end{center}

In this situation, we recall that, by lemmas \ref{C(h)_free_choice} and \ref{eta} and lemma 2.4 of \cite{haykazyan-moosa}, we may freely pick the morphisms defining the $\eta_{\mathbb{U}_i}$ and $C\mathcal{G}(h)$ among all the possible choices, as long as they have appropriate domains and codomains. With this in mind, we may in fact assume that the diagram with the concrete maps defining the morphisms of covers commutes, not just the diagram with the theories.

\end{proof}
\end{subsection}

\begin{subsection}{The isomorphism $\varepsilon : id_{FCG_A} \simeq GC$}

\begin{defi}
Let $\mathcal{G} \rightarrow A$ be an object of $FCG_A$. Let $\mathcal{G}'$ be the extension of $\mathcal{G}$ constructed with the cover $C(\mathcal{G})$ by creating adequate copies of objects and morphisms in $\mathcal{G}$. Let $S$ denote the new sort in the cover $C(\mathcal{G})$. Let $\mathcal{G}''(C(\mathcal{G}))$ be the extension of $G(C(\mathcal{G}))$ with one extra object $S_a$ in each isomorphism class.

The set of maps of $\varepsilon_{\mathcal{G}}$ is defined as $\varepsilon_{\mathcal{G}}:=\lbrace g\circ f : O_{1, i_1} \rightarrow O_{2, i_2} \, |\, i_1 \in Ob(\mathcal{G}), i_2 \in Ob(GC(\mathcal{G})), \exists a \in A, \, g \in Hom_{\mathcal{G}''(C(\mathcal{G}))}(S_a, O_{2, i_2}) \wedge f \in Hom_{\mathcal{G}'}(O_{1, i_1}, S_a)\rbrace.$
\end{defi}

\begin{prop}\label{epsilon_independent_case}
The set of maps $\varepsilon_{\mathcal{G}}$ is an isomorphism in $FCG_A$.
\end{prop}
\begin{proof}
As each of the maps is a bijection, it suffices to check that $\varepsilon_{\mathcal{G}}$ is a morphism. 

First, regarding 0-definability in $\mathbb{U}$, stable embeddedness yields a 0-definable set of $\mathbb{U}$.

Then, as we have already seen a few times, condition (B) of being a morphism is easy to check, since each map is of the form $g\circ f$, where $g$ is a morphism in $\mathcal{G}''(C(\mathcal{G}))$.

Now, let's check condition (A). Let $a$ be an element of $A$. Let $g, g'$ be morphisms in $\mathcal{G}''(C(\mathcal{G}))$ with source $S_a$ and targets in $\mathcal{G}(C(\mathcal{G}))$. Let $f,f'$ be morphisms in $\mathcal{G}'$ with target $S_a$ and sources in $\mathcal{G}$. Let $\gamma$ be a morphism in $\mathcal{G}$ such that $gf\gamma$ has the same domain as $g'f'$. We seek some morphism $\beta$ in $\mathcal{G}(C(\mathcal{G}))$ such that $\beta g'f'=gf\gamma$. We notice that $h:=f\gamma f'^{-1}$ is an element of $\mathrm{Aut}_{\mathcal{G}'}(S_a)$. Using the properties of $\mathcal{G}'$ and $\mathcal{G}''(C(\mathcal{G}))$, we deduce that $f\gamma f'^{-1} \in \mathrm{Aut}_{\mathcal{G}'}(S_a) = \mathrm{Aut}_{C(\mathcal{G})}(S_a / \mathbb{U}) = \mathrm{Aut}_{\mathcal{G}''(C(\mathcal{G}))}(S_a)$.
Thus, the map $g f \gamma f'^{-1}$ is a morphism in $\mathcal{G}''(C(\mathcal{G}))$. Since it is in the same connected component as the morphism $g'$, there exists a morphism $\beta$ in $\mathcal{G}''(C(\mathcal{G}))$ such that $\beta g' = g f \gamma f'^{-1}$.

\begin{center}
\begin{tikzcd}S_{a} \arrow[r, "g"] & O_{2,j} & O_{2,j'} \arrow[l, dotted, "\exists ? \, \beta" description ] \\
O_{1,i} \arrow[u, "f"]& O_{1,i'}\arrow[l, "\gamma"]  \arrow[r, "f'" description]& S_{a}\arrow[llu, "h" description]\arrow[lu, "gh" description]\arrow[u,  "g'" description]  \end{tikzcd}
\end{center}

Since the source and target of $\beta$ are objects of $\mathcal{G}(C(\mathcal{G}))$, it belongs to the full subgroupoid $\mathcal{G}(C(\mathcal{G}))$.
The other inclusion in condition (A) is proved similarly.
\end{proof}

\begin{prop}
The family of maps $\varepsilon$ defines a natural isomorphism between the identity functor of $FCG_A$ and the functor $GC$.
\end{prop}
\begin{proof}
We have to prove that the naturality squares commute. Given a morphism $h : \mathcal{G}_1 \rightarrow \mathcal{G}_2$ in $FCG_A$, and an element $a$ of $A$, this amount to checking that the following square commutes :

\begin{center}
\begin{tikzcd} \mathcal{G}_{1, a}\arrow[r, "{\varepsilon_{\mathcal{G}_1}}_a"] \arrow[d, "h_a"]& GC(\mathcal{G}_1)_{a} \arrow[d, "GC(h)_a"] \\
 \mathcal{G}_{2, a}\arrow[r, "{\varepsilon_{\mathcal{G}_2}}_a"]& GC(\mathcal{G}_2)_{a}    \end{tikzcd}
\end{center}

Checking this consists in writing down the maps in $GC(\mathcal{G}_1)_{a}\circ {\varepsilon_{\mathcal{G}_1}}_a$ and in ${\varepsilon_{\mathcal{G}_2}}_a \circ h_a$ as compositions of maps in all the groupoids involved, finding subterms $\tau$ in the compositions that define automorphisms of the object $S_a$, using the binding groupoid statements to extend the $\tau$ to automorphisms of one of the internal covers above $a$, and using stable embeddedness to extend them again to automorphisms of $C(h)=(C(\mathcal{G}_1)\cup C(\mathcal{G}_2), \bigcup\limits_a g_a h_a f_a)$.
\end{proof}

\end{subsection}

\end{section}

\begin{section}{Analysable covers and simplicial groupoids}

We shall now try and generalise the constructions to the case of 1-analysable covers where the fibers are not assumed to be independent.
We follow the suggestion made in \cite{hrushovski}, in the remark following Lemma 3.1, to use simplicial groupoids.

\begin{defi}\label{defi_simplicial_groupoid}
Let $A$ be a definable set. \begin{itemize}
    \item A definable \textit{simplicial groupoid} $\mathcal{G}$ over $A$ is given by a sequence of 0-definable groupoids $\mathcal{G}_n$, for $n \geq 1$, along with morphisms of groupoids $\iota_{n,m} : \mathcal{G}_n \hookrightarrow \mathcal{G}_m$, for $n < m$, such that :

    \begin{enumerate}
        \item For each $n$, $\mathcal{G}_n$ is a groupoid over $A^{=n}/E_n$, where $E_n$ is the equivalence relation associated to the action of the symmetric group $\mathfrak{S}_n$ on the set $A^{=n}$. Informally, $\mathcal{G}_n$ is over the finite subsets of $A$ of size $n$.

        \item The operation $(n < m) \mapsto \iota_{n,m}$ behaves well with compositions, i.e. $\iota_{n,m} \circ \iota_{k,n} = \iota_{k,m}$, for all $k < n < m$.
        
        \item For each $n < m$, the morphism of groupoids $\iota_{n, m}$ is compatible with the relation $R(\overline{c}, \overline{d}) \leftrightarrow \overline{c} \subseteq \overline{d}$. See definition \ref{morphisms_of_non_connected_groupoids}. Moreover, all the maps comprising of this morphism of groupoids are injective. These maps may be referred to as \textit{inclusions}, and we may then qualify the morphism $\iota_{n, m}$ as being \textit{injective}.
    \end{enumerate}
    
 We may denote the groupoid in dimension $n$ using the subscript or superscript : $\mathcal{G}^n$ or $\mathcal{G}_n$. When dealing with the connected component over some finite set $\overline{c} \in A^{< \omega}$ of size $n$, we may write $\mathcal{G}^n_{\overline{c}}$ or $\mathcal{G}_{n,{\overline{c}}}$, or even $\mathcal{G}_{\overline{c}}$.

    \item If $\overline{c} \subset \overline{d} \in A^{<\omega}$ are finite subsets of $A$ of size $n$ and $m$ respectively, if $O_{\overline{c}}$ and $O_{\overline{d}}$ are objects in the groupoids $\mathcal{G}_{\overline{c}}$ and $\mathcal{G}_{\overline{d}}$ respectively, we let $\mathrm{Hom}_{\mathcal{G}}(O_{\overline{c}}, O_{\overline{d}})$ denote the set of maps in the morphism of groupoids $\iota_{n,m}$ from the object $O_{\overline{c}}$ to the object $O_{\overline{d}}$.

    \item We will mainly be considering simplicial groupoids over $A$ with the following property : for all $\overline{d}$ in $A^{<\omega}$, and $c_1, ..., c_k$ defining a partition of $\overline{d}$, if $O_{i_{\overline{d}}}, O_{i_{c_1}}, ..., O_{i_{c_k}}$ are objects in the corresponding groupoids, and $\iota_1 :   O_{i_{c_1}}\rightarrow O_{i_{\overline{d}}}$, ..., $\iota_k : O_{i_{c_k}} \rightarrow O_{i_{\overline{d}}}$ are inclusions between the adequate groupoids, then the map $\iota_1 \sqcup ... \sqcup \iota_k : O_{i_{c_1}} \sqcup ... \sqcup O_{i_{c_k}} \rightarrow  O_{i_{\overline{d}}}$ is a bijection that identifies $ (O_{i_{\overline{d}}}, \iota_1, ... , \iota_k)$ with the (set-theoretic) disjoint union $O_{i_{c_1}} \sqcup ... \sqcup O_{i_{c_k}}$.
    
    This property will be referred to as the \textit{disjoint union property}.
    
    \item A simplicial groupoid $\mathcal{G}$ is said to be finitely faithful if all the $\mathcal{G}_n$ are, for $n \geq 1$.
    
    \end{itemize}
\end{defi}

\begin{rem}
The morphisms of groupoids $\iota_{n,m}$ defining the simplicial structure only link the connected components one would expect them to.
\end{rem}

\begin{defi}
A morphism of simplicial groupoids is defined as a family of morphisms of groupoids that is compatible with the inclusion morphisms defining the simplicial structures.

More explicitly, if $(\mathcal{G}_1,(\iota^1_{n,m})), (\mathcal{G}_2,(\iota^2_{n,m}))$ are simplicial groupoids over $A$, then a morphism $ H : \mathcal{G}_1 \rightarrow \mathcal{G}_2$ is given by a family of morphisms $H_n : \mathcal{G}_{1,n} \rightarrow \mathcal{G}_{2,n}$ that are compatible with the equality relations $R_n \subseteq (A^{=n}/E_n)\times (A^{=n}/E_n) $, and such that, for all integers $n,m$ such that $1 \leq n < m$, the following diagram of morphisms of groupoids commutes : 
\begin{center}
\begin{tikzcd}
\mathcal{G}_{1,m} \arrow[rr, "H_m"]&& \mathcal{G}_{2,m}
\\
\mathcal{G}_{1,n} \arrow[u, "\iota^1_{n,m}"] \arrow[rr, "H_n"]&& \mathcal{G}_{2,n} \arrow[u, "\iota^2_{n,m}"]
\end{tikzcd}
    
\end{center}
\end{defi}

\begin{rem}
The reader might be wondering why we didn't use a definition from category theory to define simplicial groupoids as simplicial objects in the category of groupoids. 

The main reason is to avoid degeneracies, which are painful to deal with. In fact, most of the constructions that appear below are much easier within our framework. However, an algebraic topologist more technically gifted than the writer might disagree, and find that having to deal with the degeneracies is worth it.
\end{rem}

We now prove a general lemma that will be useful when studying 1-analysable covers.

\begin{lemma}\label{criterion_automorphism}
Let $(\mathbb{U}, T, f : T \rightarrow A)$ be a cover of $\mathbb{U}$.

Let $\tau : T \rightarrow T$ be a map compatible with $f : T \rightarrow A$. 

Then the map $\tau \cup id_{\mathbb{U}}$ is an automorphism of $(\mathbb{U}, T)$ if and only if, for all finite subsets $\overline{c} \subseteq A$, its restriction is an automorphism of $(\mathbb{U}, \bigcup\limits_{a \in \overline{c}}T_{a})$.
\end{lemma}

\begin{proof}
One direction is clear, as the structure on $(\mathbb{U}, \bigcup\limits_{a \in \overline{c}}T_{a})$ is the structure induced by the $\overline{c}$-definable sets.

Conversely, let's assume $\tau \cup id_{\mathbb{U}}$ satisfies the local condition. Let $\overline{b}$ be a finite tuple of elements of $(\mathbb{U}, T)$, and $R$ be a 0-definable relation in $(\mathbb{U}, T)$ such that $(\mathbb{U}, T) \models R(\overline{b})$. Let $\overline{c}$ be a tuple of elements of $A$ such that $\overline{b} \subseteq \mathbb{U}\cup \bigcup\limits_{a \in \overline{c}}T_{a}$. Let $R'$ be the restriction of the relation $R$ to $\mathbb{U}\cup \bigcup\limits_{a \in \overline{c}}T_{a}$. By definition of the structure of $\mathbb{U}\cup \bigcup\limits_{a \in \overline{c}}T_{a}$, the relation $R'$ is 0-definable in $(\mathbb{U}, \bigcup\limits_{a \in \overline{c}}T_{a})$.

Therefore, $(\mathbb{U}, \bigcup\limits_{a \in \overline{c}}T_{a}) \models R'((\tau\cup id_{\mathbb{U}})(\overline{b}))$.
So $(\mathbb{U}, T) \models R((\tau\cup id_{\mathbb{U}})(\overline{b}))$
\end{proof}

\begin{subsection}{Constructing a simplicial groupoid from a 1-analysable cover}\label{subsection_simplicial_from_cover}

In this subsection, we let $\mathbb{U}' = (\mathbb{U}, S, f : S \rightarrow A)$ denote a 1-analysable cover of $\mathbb{U}$ with the following properties : 
    \begin{itemize}
        \item (Uniform Finite Generatedness) For each integer $n \geq 1$, for each finite tuple $\overline{c} \in A^{=n}$, the language of $(\mathbb{U}, \bigcup\limits_{a \in \overline{c}}S_{a})$ is finitely generated over that of $\mathbb{U}$, and this happens uniformly for $\overline{c}$ in $A^{=n}$. See definition \ref{defi_covers} for a formal description of \textquotedblleft uniformly finitely generated languages".
        \item (Local Stable Embeddedness) For each finite subset $C \subseteq A$, the structure $(\mathbb{U}, \bigcup\limits_{a \in C}S_{a})$, with the structure induced by the $C$-definable sets of $(\mathbb{U}, S)$, is stably embedded in $(\mathbb{U}, S)$.
    \end{itemize}

\begin{theo}\label{theo_defi_simplicial_groupoid}
There exist simplicial groupoids $\mathcal{G}$ and $\mathcal{G}'$, over the set $A$, that are 0-definable in $\mathbb{U}$ and $\mathbb{U}'^{eq}$ respectively, with the following properties:

\begin{itemize}
    \item The simplicial groupoid $\mathcal{G}'$ is an extension of $\mathcal{G}$ with the extra object $S_{\overline{c}}$ in the connected component over $\overline{c}$. Moreover, the set-theoretic inclusions between the $S_{\overline{c}}$ belong to the simplicial groupoid.
    
    More explicitly : for each degree $n \geq 1$, for each element $\overline{c} \in A^{=n} / E_n$, the connected groupoid $\mathcal{G}'_{n, \overline{c}}$ is an extension of the connected groupoid $\mathcal{G}_{n, \overline{c}}$, with only one extra object, whose underlying set is the set $S_{\overline{c}} = \bigcup\limits_{a \in \overline{c}} S_a$. Moreover, if we have $n<m$, $\overline{c} \in A^{=n}/E_n$, $\overline{d} \in A^{=m}/ E_m$ and $\overline{c} \subseteq \overline{d}$, then the morphism of groupoids $\iota'_{n,m} : \mathcal{G}'_n \rightarrow \mathcal{G}'_m$ contains the set-theoretic inclusion map $S_{\overline{c}} \hookrightarrow S_{\overline{d}}$.

    \item The automorphism groups $\mathrm{Aut}_{\mathcal{G}'}(S_{\overline{c}})$, along with their action on $S_{\overline{c}}$, are canonically isomorphic to the groups $\mathrm{Aut}(S_{\overline{c}}/\mathbb{U})$.
\end{itemize}
\end{theo}
\begin{proof}
We shall define a simplicial groupoid which captures the whole structure of the cover $\mathbb{U}'$ and generalizes the case of independent fibers.

For each positive integer $n$, and each tuple $\overline{c} \in A^{=n}$, we notice that the structure $(S_{\overline{c}}, \mathbb{U})$ is an internal cover of $\mathbb{U}$. Using the hypothesis of uniformly finitely generated language, we can find a groupoid $\mathcal{G}^n$ over $A^{=n}$ which is 0-definable in $\mathbb{U}$, with an extension $\mathcal{G}'^n$ that is 0-definable in $(\mathbb{U}')^{eq}$. This extension has exactly one extra object $S_{a_1}\cup ... \cup S_{a_n}$ in the connected component over $(a_1,...,a_n)$ and satisfies the binding groupoid statement :  $\mathrm{Aut}_{\mathcal{G}'^n}(S_{a_1}\cup ... \cup S_{a_n}) = \mathrm{Aut}(S_{a_1}\cup ... \cup S_{a_n} / \mathbb{U})$.

Now, we want groupoids that are over $A^{=n}/ E_n$, and not over $A^{=n}$. To do this, it suffices to glue together the groupoids that correspond to a common orbit under the action of $\mathfrak{S}_n$. This is made possible by the existence of the common object $S_{\overline{c}}$, whose automorphism group is the same in all the connected components where it appears.

We define the simplicial structure using set-theoretic inclusions. For each inclusion of finite subsets of $A$, $\overline{c} \subseteq \overline{d}$, we find the definable map $S_{\overline{c}} \hookrightarrow S_{\overline{d}}$. We then create sets of maps between the groupoids $\mathcal{G}'_{\overline{c}}$ by precomposing and postcomposing these maps with morphisms in the groupoids $\mathcal{G}'^n$. 
To see that these sets of maps define groupoid morphisms, we need to use the second hypothesis on $(\mathbb{U}, S)$, that enables us to lift automorphisms and get surjective group homomorphisms $\mathrm{Aut}(S_{\overline{d}}/ \mathbb{U}) \twoheadrightarrow \mathrm{Aut}(S_{\overline{c}}/ \mathbb{U})$, for each inclusion $\overline{c} \subseteq \overline{d}$ :

\begin{prop}\label{condition_A_requires_lift}
Let $\iota$ be the sets of inclusion maps defined above. Then these sets of maps induce groupoid morphisms between the $\mathcal{G}'^n$, that define a simplicial groupoid.
\end{prop}

\begin{proof}
The main point is to check condition (A). We shall first prove the following partial result : If $\iota_0 : S_{\overline{c}} \rightarrow S_{\overline{d}}$ is a set-theoretic inclusion, and $\sigma, \tau$ are elements of $\mathrm{Aut}(S_{\overline{c}})$ and $\mathrm{Aut}(S_{\overline{d}})$ respectively, then there exist elements $\sigma'$ and $\tau'$ in $\mathrm{Aut}(S_{\overline{c}})$ and $\mathrm{Aut}(S_{\overline{d}})$ respectively making the following diagrams commute :

\begin{center}
\begin{tikzcd} 
S_{\overline{c}}\arrow[r, "\sigma"] \arrow[d, "\iota_0"]& S_{\overline{c}} \arrow[d, "\iota_0"]&& S_{\overline{c}}\arrow[r, dotted, "\exists \sigma'"] \arrow[d, "\iota_0"]&S_{\overline{c}} \arrow[d, "\iota_0"] \\

 S_{\overline{d}}\arrow[r, dotted, "\exists \tau'"]&S_{\overline{d}}  && S_{\overline{d}}\arrow[r, "\tau"]& S_{\overline{d}}
 \end{tikzcd}
\end{center}

\noindent We first use the binding groupoid statements. We define binary relations $R_i$ on $S_{\overline{c}}$ : $R_i(x, ) \leftrightarrow f(x) = c_i$ . These relations are  0-definable in $(\mathbb{U}, S_{\overline{c}})$, and thus preserved by $\sigma$. Therefore, we see that $\sigma$ is a union $\sigma_{c_1} \cup ... \cup \sigma_{c_n}$ of automorphisms of the $S_{c_i}$, for $c_i$ a coordinate of $\overline{c}$. We can also check that the restriction $\bigcup\limits_{c_i \in \overline{d}}\sigma_{c_i}$ is indeed an automorphism of $S_{\overline{d}}$. This restriction is the automorphism $\tau'$ we were looking for.

\noindent Now, to find $\sigma'$, we use the second hypothesis on the cover $(\mathbb{U}, S).$ We have assumed that $\tau$ induces an automorphism of the structure $(\mathbb{U}, \bigcup\limits_{a \in \overline{d}}S_a)$. Since $(\mathbb{U}, \bigcup\limits_{a \in \overline{c}}S_a)$ is a cover of $(\mathbb{U}, \bigcup\limits_{a \in \overline{d}}S_a)$ by hypothesis, such an automorphism exists.

Since the maps in $\iota$ are of the form $\alpha \iota_0 \beta$, where $\alpha$ and $\beta$ are morphisms in the $\mathcal{G}'_n$, the result proved above yields the full condition (A).

\end{proof}
\end{proof}

\begin{rem}
It is somewhat surprising that condition (A), which looks a bit formal, should translate as the model-theoretic stable embeddedness condition above. 
\end{rem}

\begin{prop}\label{remark_product}
The simplicial groupoid built here satisfies the disjoint union property, in addition to being finitely faithful : 
    Let $\overline{d}$ be a finite subset of $A$, and $c_1, ..., c_k$ define a partition of $\overline{d}$. Let $O_{i_{\overline{d}}}, O_{i_{c_1}}, ..., O_{i_{c_k}}$ be objects in the corresponding groupoids. Let $\iota_1 :   O_{i_{c_1}}\rightarrow O_{i_{\overline{d}}}$, ..., $\iota_k : O_{i_{c_k}} \rightarrow O_{i_{\overline{d}}}$ be maps belonging to the inclusion morphisms of groupoids. Then the map $\iota_1 \sqcup ... \sqcup \iota_k : O_{i_{c_1}} \sqcup ... \sqcup O_{i_{c_k}} \rightarrow  O_{i_{\overline{d}}}$ is a bijection that identifies $ (O_{i_{\overline{d}}}, \iota_1, ... , \iota_k)$ with the disjoint union $O_{i_{c_1}} \sqcup ... \sqcup O_{i_{c_k}}$.
\end{prop}

\begin{proof} The statement is true for the objects $S_{\overline{d}}$, $S_{c_j}$ and the set-theoretic inclusions. Now, condition (A) enables us to compare this setting to any other setting. Indeed, if $m_{ab} : O_{i_{ab}} \rightarrow S_{ab}$ is any morphism in $\mathcal{G}'$, condition (A) implies that there exist morphisms $m_a : O_{i_a} \rightarrow S_a$ and $m_b : O_{i_b} \rightarrow S_b$ in the degree 1 groupoid such that the following diagram commutes :

\begin{center}
    
    \begin{tikzcd}
    
     S_a \arrow[rd] && O_{i_a}\arrow[rd, "\iota_a" description] \arrow[ll, "{m_a}"]
    \\ & S_{ab} && O_{i_{ab}}\arrow[ll, "{m_{ab}}"]
    \\  S_b \arrow[ru] && O_{i_b} \arrow[ru, "\iota_b" description]\arrow[ll, "{m_b}"]
    \end{tikzcd}
    
\end{center}
    In a nutshell : $O_{i_a} \sqcup O_{i_b} \simeq S_a \sqcup S_b = S_{ab} \simeq O_{i_{ab}}$
\end{proof}

    We recall that this property will be referred to as the \textit{disjoint union property}.

\begin{rem}
These properties of the binding simplicial groupoids we build are preserved under isomorphisms of simplicial groupoids. Thus, our future assumptions on the simplicial groupoids we work with will be justified.
\end{rem}

\begin{prop}
Let $(\mathrm{Aut}(S_{\overline{d}}/ \mathbb{U}) \twoheadrightarrow \mathrm{Aut}(S_{\overline{c}}/ \mathbb{U}))_{\overline{c}\subseteq \overline{d}}$ be the projective system of groups mentioned at the end of Theorem \ref{theo_defi_simplicial_groupoid}. Then its projective limit in the category of groups is isomorphic to the automorphism group $\mathrm{Aut}(S / \mathbb{U}).$
\end{prop}

\begin{proof}
Using the second hypothesis on the cover $(\mathbb{U}, S)$, and lemma \ref{criterion_automorphism}, we see that each element of the projective limit yields an automorphism of $S$ over $\mathbb{U}$. Conversely, an automorphism of $S$ over $\mathbb{U}$ induces a family of automorphisms in the $\mathrm{Aut}(S_{\overline{c}} / \mathbb{U})$. Such a family is in the projective limit, since the projections act as restrictions.
\end{proof}

The next proposition is a generalization of proposition 3.6 in \cite{haykazyan-moosa}. It merely consists in making the compactness argument found there uniform over $A$.

\begin{prop}\label{criterion_uniformly_finite_language}
Let $\mathbb{V} = (\mathbb{U}, S, S \rightarrow A)$ be a 1-analysable cover of $\mathbb{U}$ over $A$. Assume that the automorphism groups $\mathrm{Aut}(S_a / \mathbb{U})$ along with their actions on the $S_a$ are uniformly interpretable in $\mathbb{V}$. Then the language of $(\mathbb{U}, S_a)$ is finitely generated over that of $\mathbb{U}$, uniformly in $a$.
\end{prop}

\begin{proof}
By compactness and 1-analysability, there exists a 0-definable set $B$ in $\mathbb{V}$ with a 0-definable surjective map $g : B \rightarrow A$, a 0-definable set $X$ in $\mathbb{U}$ and a 0-definable map $F : S \times_A B \rightarrow X$ such that for all $a$ in $A$, for all $b$ in $B_a$, the map $F_b$ is an injection $F_b : S_a \hookrightarrow X$.

Let $G(x,y,b,b')$ be the formula \textquotedblleft There exists some $a$ in $A$ such that the map $F_b^{-1} \circ F_{b'} : S_a \rightarrow S_a$ is a well-defined bijection sending $x$ to $y$".
Now, we shall use the uniform interpretability of $\mathrm{Aut}(S_a / \mathbb{U})$. Let $p(v,z)$ be the partial type $\lbrace z \in A \rbrace \cup \lbrace $\textquotedblleft $G_v $ is a permutation of $S_z"$ $ \rbrace \cup \lbrace \forall x \, \, (R((G_v\cup id_{\mathbb{U}})(x);z) \leftrightarrow R(x ; z)) \, | \, R : $ 0-definable relation in $\mathbb{U} \rbrace$.

We know that $p(v,z) \models G_v \in \mathrm{Aut}(S_z / \mathbb{U})$. Thus, by compactness, there exists formulas $R_1(x,y),...,R_n(x,y)$ in the language of $\mathbb{V}$ such that, for all $a$ in $A$, for all $v$ such that $G_v$ is a permutation of $S_a$, the map $G_v\cup id_{\mathbb{U}}$ is an automorphism of $(\mathbb{U}, S_a)$ if and only if it preserves the relations $R_1(x ,a), ..., R_n(x,a).$

Now, for $a$ in $A$, let $\Sigma_a$ be the language generated by the relations $R_i(x,a)$ along with the relation $R(x,y,b)$, the latter being defined as $F_b(x) = y$.
We thus get a family of finitely generated languages that vary uniformly over $A$. It remains to show that these languages generate the languages $L_a$ of $(\mathbb{U}, S_a)$. To simplify notations, let $\mathbb{V}_a$ and $\mathbb{V}_{0,a}$ be the structure $(\mathbb{U}, S_a)$ in the languages $L_a$ and $L \cup \Sigma_a$ respectively. By saturation, it suffices to prove that $\mathrm{Aut}(\mathbb{V}_a) = \mathrm{Aut}(\mathbb{V}_{0,a})$.
By stable embeddedness, it suffices to show that $\mathrm{Aut}(\mathbb{V}_a / \mathbb{U}) = \mathrm{Aut}(\mathbb{V}_{0,a}/ \mathbb{U})$. Let $\sigma$ be an automorphism of $\mathbb{V}_{0,a}$ fixing $\mathbb{U}$ pointwise. Let $b$ be a parameter in $B$ such that $F_b : S_a \hookrightarrow X$ is an injection. Since the formula defining $F$ is in the language $\Sigma_a$, we know that $F_{\sigma(b)}$ is also an injection from $S_a$ to $X$. In addition, its image is the same as that of $F_b$, for $X$ is a subset of $\mathbb{U}$ and $\sigma$ fixes $\mathbb{U}$ pointwise.
We now compute, for $x$ in $S_a$ : $F_b(x) = \sigma(F_b(x)) = F_{\sigma(b)}(\sigma(x))$. Thus the map $F_{\sigma(b)}^{-1} \circ F_b$ is equal to the map $\sigma|_{S_a}$. Finally, if $v=(\sigma(b), b)$, applying the property of the map $G_v$ described in the third paragraph of our proof, we deduce that $\sigma$ is an automorphism of $\mathbb{V}_a$, as desired.
\end{proof}

\begin{coro}
Let $\mathbb{V} = (\mathbb{U}, S, S \rightarrow A)$ be a 1-analysable cover of $\mathbb{U}$ over $A$. Let $n$ be an integer. Assume that the automorphism groups $\mathrm{Aut}(S_{\overline{c}} / \mathbb{U})$ along with their actions on the $S_{\overline{c}}$ are uniformly interpretable in $\mathbb{V}$, for $\overline{c} \in A^{=n}$. Then the language of $(\mathbb{U}, S_{\overline{c}})$ is finitely generated over that of $\mathbb{U}$, uniformly in $\overline{c}$.
\end{coro}

\begin{proof}
We first note that the automorphism groups $\mathrm{Aut}(S_{\overline{c}} / \mathbb{U})$ along with their actions on the $S_{\overline{c}}$ are uniformly interpretable in $\mathbb{V}$, for $\overline{c} \in A^{=n}$, if and only if their counterparts $\mathrm{Aut}(\prod\limits_{a \in \overline{c}}S_{\overline{a}} / \mathbb{U})$ are.

Moreover, the languages of the covers $(\mathbb{U}, S_{\overline{c}})$ are uniformly finitely generated over $\mathbb{U} $, for $\overline{c} \in A^{=n}$, if and only if the same holds true for the covers $(\mathbb{U}, \prod\limits_{a \in \overline{c}}S_{\overline{a}})$.

Thus, if $g : S^n \rightarrow A^n$ is the map $f\times ... \times f$, and $S'=g^{-1}(A^{=n})$, we can apply proposition \ref{criterion_uniformly_finite_language} to the 1-analysable cover $(\mathbb{U}, S', g : S' \rightarrow A^{=n})$. This yields the result.
\end{proof}

\begin{rem}
As in \cite{haykazyan-moosa} (Proposition 3.6), this result shows that the condition of uniformly finitely generated language is necessary in order to get 0-definable groupoids and not just type-definable ones. 
\end{rem}

\end{subsection}

\begin{subsection}{Constructing a 1-analysable cover from a simplicial groupoid}

In this subsection, we let $\mathcal{G}$ be a finitely faithful 0-definable simplicial groupoid over $A$ satisfying the properties of proposition \ref{remark_product}.

\begin{defi}
    A \textit{commuting system of inclusions} is a family of objects $(O_{i_c})_{c \subseteq A, \, |c| < \omega}$ with one object in each connected component, along with a family of maps $(\iota_{c,d} : O_{i_c} \hookrightarrow O_{i_d})_{c \subset d}$, such that, if $c \subset d \subset e$, then $\iota_{d,e} \circ \iota_{c,d} = \iota_{c,e}$.
    \end{defi}

    \begin{prop}\label{exists_system_of_inclusions}
    Let $\mathcal{G}$ be a definable simplicial groupoid over $A$. Then there exists a commuting system of inclusions in $\mathcal{G}$.
    \end{prop}
    We build such a system inductively. We will first need to prove weaker results :

\begin{lemma}\label{fill_square}
Let $\overline{a}$, $\overline{b}$, $\overline{c}$ be finite subsets of $A$, possibly not disjoint, such that $\overline{a} \subset \overline{ab}$, $\overline{a} \subset \overline{ac}$,  $\overline{ab}\subset \overline{abc}$, and $\overline{ac} \subset \overline{abc}$.
Let $O_{\overline{a}},O_{\overline{ab}}, O_{\overline{ac}}, O_{{\overline{abc}}}$ be objects in $\mathcal{G}$, in the appropriate connected components. Let $\iota_{1} : O_{\overline{a}} \hookrightarrow O_{\overline{abc}} , \iota_3 : O_{\overline{ac}} \hookrightarrow O_{\overline{abc}}$ be inclusions in the morphisms of groupoids defining the simplicial structure of $\mathcal{G}$. Then there exists a unique inclusion $\iota_4 : O_{\overline{a}} \hookrightarrow O_{\overline{ac}}$ making the following diagram commute :
\begin{center}

\begin{tikzcd}

 O_{\overline{abc}}
\\
& O_{\overline{ac}}\arrow[lu , "{\iota_3}"]
\\
O_{\overline{a}}\arrow[uu , "{\iota_1}"] \arrow[ru , "{\exists \, ! \, \iota_4}" description]
\end{tikzcd}

\end{center}
\end{lemma}
\begin{proof}
Uniqueness comes from injectivity of $\iota_3$. Let us now prove existence.

First, we recall the equality of sets of maps : $$ Hom(O_{\overline{a}}, O_{\overline{abc}}) = Hom(O_{\overline{ac}}, O_{\overline{abc}})\circ Hom(O_{\overline{a}}, O_{\overline{ac}})  $$

So there exist inclusion maps $j_3, j_4$ such that $j_3 \circ j_4 =  \iota_1$.
Then, by condition (A), there exists $\sigma \in \mathrm{Aut}(O_{\overline{ac}})$ such that $\iota_3 \sigma = j_3$. Then, the map $\iota_4 := \sigma j_4$ satsifies $\iota_3 \iota_4 = \iota_3 \sigma j_4 = j_3 j_4 =  \iota_1$, as desired.

\begin{center}
    
\begin{tikzcd}

O_{\overline{abc}}\arrow[rrr , "{id}"]&&&O_{\overline{abc}}
\\
& &&O_{\overline{ac}}\arrow[u , "{\exists \,j_3}"]\arrow[rr , "{\exists \, \sigma}" ]&& O_{\overline{ac}}\arrow[llu , "{\iota_3}" description]
\\
&&&O_{\overline{a}}\arrow[llluu , "{\iota_1}"] \arrow[rru , "{\sigma j_4}" description]\arrow[u, "\exists \, j_4"]

\end{tikzcd}
\end{center}

\end{proof}

\begin{lemma}\label{finite_systems_of_inclusions}
Let $k$ be an integer. Let $X \subset P(A)$ be a finite collection of finite subsets of $A$. Let's assume that $|X|=2^k-1$, and that the relation of inclusion on $X$ is isomorphic to that of the nonempty subsets of some $B \subseteq A$, where $|B|=k$. Let $(O_x)_{x \in X}$ be a family of objects of $\mathcal{G}$ in the appropriate connected components. Then there exists a commuting system of inclusions between the $O_x$. 
\end{lemma}
\begin{proof}
We prove the result by induction on $k$. For $k = 1$, there is nothing to prove. Now, let $k \geq 2$, let $X\subset P(A)$ and $(O_x)_{x \in X}$ be as in the lemma. Let $x_1, ..., x_k$ be the minimal elements of $X$. Using the induction hypothesis, there exists a commuting system of inclusions for the family $(O_x)_{x \in X, x_k \not\subseteq x }$. There also exists a commuting system of inclusions for the family $(O_x)_{x \in X, x_k \subset x }$. It remains to add $x_k$ to the latter system, and connect the two systems together.

From now on, to simplify notations, we shall assume that the $x_i$ are in fact elements of $A$, and that $X=P(B)$ for $B=\lbrace x_1,..., x_k \rbrace$. 

\noindent To add $x_k$ to the system of subsets of $B$ that strictly contain $x_k$, pick one inclusion $O_{x_k} \hookrightarrow O_B$. Then, if $i, j \neq k$, since there are already inclusions as in the following diagram, we can apply lemma \ref{fill_square} :

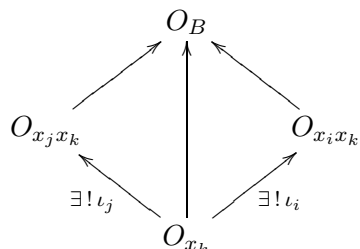
\xymatrix{&&&&&O_{B}
\\
&&&&O_{x_j x_k}\ar[ru] && O_{x_i x_k}\ar[lu]_{}
\\
&&&&&O_{x_k}\ar[uu]^{} \ar[ru]_{\exists \, ! \, \iota_i}\ar[lu]^{\exists \, ! \, \iota_j}}

Now, using injectivity of all the inclusions, commutativity of these squares, and the fact that the system between the $(O_x)_{x \in X, x_k \subset x }$ is commutative, we can prove that the all new squares created by our new maps $\iota_i$ commute. For instance, if $i\neq j$, we study the following diagram :

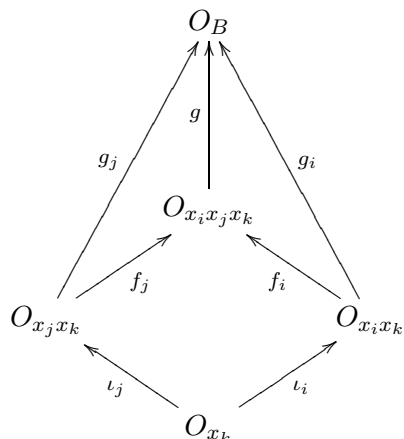
\xymatrix{&&&&&O_{B}
\\
\\
&&&&&O_{x_i x_j x_k}\ar[uu]^{g}
\\
&&&&O_{x_j x_k}\ar[ru]_{f_j} \ar[ruuu]^{g_j}&& O_{x_i x_k}\ar[lu]^{f_i}\ar[luuu]_{g_i}
\\
&&&&&O_{x_k} \ar[ru]_{ \iota_i}\ar[lu]^{  \iota_j}}

\noindent Here, all the arrows except $\iota_i$ and $\iota_j$ belong to the commuting system between the $(O_x)_{x \in X, x_k \subset x }$. So the two upper triangles commute. Moreover, by definition of $\iota_i$ and $\iota_j$, the outer square commutes, i.e. $g_j \iota_j = g_i \iota_i$. Thus, we find $g f_j \iota_j = g_j \iota_j = g_i \iota_i = g f_i \iota_i$. Since $g$ is injective, we deduce $f_j\iota_j=f_i \iota_i$.

Now, it remains to connect the objects that contain $x_k$ with those that do not contain it. Again, pick an arbitrary inclusion $f:O_{x_1...x_{k-1}} \hookrightarrow O_{B}$. 
Then, if $B'=B\setminus \lbrace x_i, x_k \rbrace $, apply lemma \ref{fill_square} in the following setting :

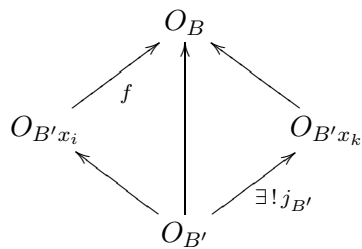
\xymatrix{&&&&&O_{B}
\\
&&&&O_{B' x_i}\ar[ru]_{f} && O_{B' x_k}\ar[lu]_{}
\\
&&&&&O_{B'}\ar[uu]^{} \ar[ru]_{\exists \, ! \, j_{B'}}\ar[lu]}

Here, we have $B'x_i = \lbrace x_1,..., x_{k-1} \rbrace$. The map in the middle is defined precisely as the composition of the maps on the left. The map $O_{B'x_k} \rightarrow O_B$ in the upper right corner belongs to the system $(O_x)_{x\in X, \, x_k \subset x}$. We find a uniquely determined inclusion $j_{B'} : O_{B'} \hookrightarrow O_{B' x_k}$ that makes the left triangle commute. By definition, this extended system commutes.

Now, deal with subsets $B' = B \setminus \lbrace x_i, x_j, x_k \rbrace$. There is already a unique inclusion map $O_{B'} \hookrightarrow O_B$, given by any path ending with $f$ :

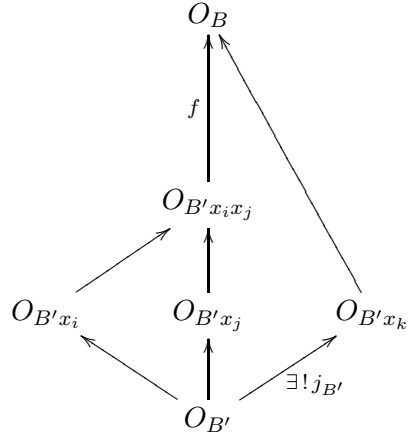
\xymatrix{&&&&&O_{B}
\\
\\
&&&&&O_{B'x_i x_j }\ar[uu]^{f}
\\
&&&&O_{B'x_i}\ar[ru]_{} & O_{B'x_j}\ar[u]^{}&O_{B'x_k}\ar[luuu]
\\
&&&&&O_{B'} \ar[u]_{}\ar[lu]^{} \ar[ru]_{\exists \, ! \, j_{B'}}  }

The map $O_{B'x_k} \rightarrow O_B$ in the upper right corner belongs to the system $(O_x)_{x\in X, \, x_k \subset x}$.

Again, using injectivity and commutativity of the previous smaller systems, we find that these new maps define a commuting system.

Then deal with the cases $B'= B \setminus\lbrace x_t, x_i, x_j, x_k \rbrace $ in a similar way, applying lemma \ref{fill_square} to find uniquely determined inclusions $j_{B'} : O_{B'} \hookrightarrow O_{B'x_k}$.

By decreasing induction, we finally find uniquely defined maps $j_{\lbrace x_i\rbrace} : O_{x_i} \hookrightarrow O_{x_i x_k}$ to complete our system.
Note that once $f$ was chosen, all the maps yet to be defined were uniquely determined. 

Using the same ideas as in the third paragraph of this proof, namely, injectivity and commutativity of the \textquotedblleft older" smaller systems, we check that the whole system commutes.

\end{proof}
    
We are now ready to prove proposition \ref{exists_system_of_inclusions} :

\begin{proof}[Proof of proposition \ref{exists_system_of_inclusions}]

Let $\kappa:= |A|$. Let us pick an enumeration of $A$ : $A=\lbrace a_{\alpha} \, | \, \alpha < \kappa \rbrace$. For $\alpha < \kappa$, let $A_{\alpha}$ denote the subset $\lbrace a_{\beta} \, | \, \beta < \alpha \rbrace$. We intend to define an increasing family of commuting systems of inclusions on the $A_{\alpha}$, and then take the increasing union of these partial systems.

This inductive construction uses both saturation and the ideas in the proof of lemma \ref{finite_systems_of_inclusions}.

\noindent Let $\gamma < \kappa$. Assume that there exists an increasing family of commuting systems of inclusions on the $A_{\beta}, $ for $\beta < \gamma$. We want to find a commuting system of inclusions on $A_{\gamma}$. We may assume that the ordinal $\gamma$ is a successor ordinal, otherwise just take the union of the previously defined systems. Let us write $\gamma = \alpha+1$. Pick any object $O_{i_{a_{\alpha}}}$ in the component of $\mathcal{G}$ over $a_{\alpha}$. We wish to add this object to the system. This consists in finding objects $O_{i_{Ba_{\alpha}}}$, for $B$ a finite subset of $A_{\alpha}$, along with inclusions $\iota_{B', B} : O_{i_{B'}} \hookrightarrow O_{i_B}$, for $B' \subset B$ finite subsets of $A_{\alpha}$, when $a_{\alpha} \in B$. 

 We shall view the codes of the inclusion maps we are looking for as \textit{variables in a partial type}. 
 
 \noindent In fact, commutativity of the finite diagrams we are considering is expressible in the language of $\mathbb{U}$. Moreover, the number of variables is the number of new inclusions we have to define, so its cardinal is bounded by $\alpha$. Similarly, the set of parameters of the partial type we define is given by the codes of the inclusions that are already defined. This set of parameters is thus also bounded in cardinality by $\alpha$.
 
 \noindent Thus, if we can show that the partial type expressing commutativity of all the diagrams involved is consistent, we can use saturation of $\mathbb{U}$ to extend our system by \textquotedblleft adding the point $a_{\alpha}$". Let us now prove the consistency of this partial type. 
 
 Let $A_1,...,A_n$ be finite subsets of $A_{\alpha}$. We wish to find commuting systems of inclusions for $A_1 \cup \lbrace a_{\alpha}\rbrace , ..., A_n \cup \lbrace a_{\alpha}\rbrace $ that extend the preexisting systems  defined on $A_1, ..., A_n$ respectively. 
 
 \noindent It suffices to consider the finite set $A_1 \cup ...\cup A_n \subseteq A_{\alpha}$, for which we have a commuting system of inclusions that contains the ones on the $A_i$. Indeed, managing to extend this system to $A_1 \cup ... A_n \cup \lbrace a_{\alpha} \rbrace$ is enough. 
 
 \noindent Now, we claim that the proof of the induction step for lemma \ref{finite_systems_of_inclusions} shows that, given a commuting system of inclusions on some finite set $A' \subset A$, one can \textquotedblleft add a point $a_{\alpha}$ to it" and get an extended commuting system on $A'\cup \lbrace a_{\alpha} \rbrace$.

\end{proof}

\begin{defi}
Let $\mathcal{G}_1$, $\mathcal{G}_2$ be 0-definable simplicial groupoids over $A$ with the disjoint union property, i.e. which satisfy the conclusion of Proposition \ref{remark_product}. Let $H : \mathcal{G}_1 \rightarrow \mathcal{G}_2$ be a morphism of simplicial groupoids.
Let $\iota_1= ( O_{i_{\overline{c}}} \rightarrow O_{i_{\overline{d}}})_{\overline{c} \subset \overline{d} \in A^{<\omega}}$ and $\iota_2 = (O_{j_{\overline{c}}} \rightarrow O_{j_{\overline{d}}})_{\overline{c} \subset \overline{d} \in A^{<\omega}} $ be commuting systems of inclusions in $\mathcal{G}_1$ and $\mathcal{G}_2$ respectively.

Let $B \subseteq A$. Let $(h_{a} : O_{i_{a}} \rightarrow O_{j_{a}})_{a \in B}$ be a family of maps belonging to the morphism of groupoids $H_1$. This family is called \textit{coherent}, with respect to $\iota_1$ and $\iota_2$, if, for each finite subset $\overline{c} \subset B$,  there exists a higher degree map $h_{\overline{c}}$ in $H$ such that, for each element $a \in \overline{c}$,  the following diagram commutes : 

\begin{center}
    \begin{tikzcd}
    O_{i_{\overline{c}}} \arrow[r, " h_{\overline{c}}"] & O_{j_{\overline{c}}}
    \\
    O_{i_{a}} \arrow[r, " h_{a}"]\arrow[u, "\iota_1"] & O_{j_{a}}\arrow[u, "\iota_2"]
    \end{tikzcd}
\end{center}
\end{defi}

We shall now prove that coherent families of maps always exist, for $B=A$.

\begin{lemma}
Let $\mathcal{G}_1$, $\mathcal{G}_2$ be 0-definable simplicial groupoids over $A$ with the disjoint union property. Let $H : \mathcal{G}_1 \rightarrow \mathcal{G}_2$ be a morphism of simplicial groupoids.
Let $\iota_1= ( O_{i_{\overline{c}}} \rightarrow O_{i_{\overline{d}}})_{\overline{c} \subset \overline{d} \in A^{<\omega}}$ and $\iota_2 = (O_{j_{\overline{c}}} \rightarrow O_{j_{\overline{d}}})_{\overline{c} \subset \overline{d} \in A^{<\omega}} $ be commuting systems of inclusions in $\mathcal{G}_1$ and $\mathcal{G}_2$ respectively. Let $n$ be an integer, let $\overline{c}=a_1...a_n$ be a set of elements of $A$. Let $m_{a_1},..., m_{a_n}$ be maps in $H_1$, and $f$ be a map in $H$ in higher dimension, such that the following squares commute, for $i= 1,..., n$ : 
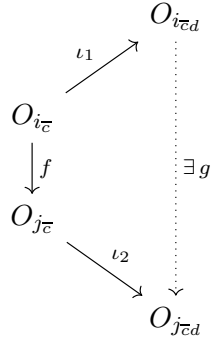
\begin{center}
\begin{tikzcd}
O_{i_{a_i}}\arrow[d, "m_{a_i}"]& & O_{i_{\overline{c}}}\arrow[d, "f"] \arrow[from=ll, "\iota_1"]  &

 \\
 O_{j_{a_i}}&   & O_{j_{\overline{c}}}  \arrow[from=ll, "\iota_2"] 
 \\
 \end{tikzcd}
\end{center}

Let $d$ be an element of $A$. Then there exists a degree 1 map $m_d : O_{i_{d}} \rightarrow O_{j_d}$ belonging to $H_1$, and a higher dimension map $g : O_{i_{\overline{c}d}} \rightarrow O_{j_{\overline{c}d}}$ in $H$ making the following diagram commute : 
\begin{center}
\begin{tikzcd}
 & O_{i_{{\overline{c}d}}}\arrow[ddd, "g"] \arrow[from=ld, "\iota_1"] \arrow[from=rd, "\iota_1"] \arrow[from=rrd, "\iota_1" description] 

 \\
 O_{i_{a_1}}\arrow[d, "m_{a_1}"]&&  O_{i_{a_n}}\arrow[d, "m_{a_n}"]   & O_{i_{d}}\arrow[d, "m_d"] &
 \\
 O_{j_{a_1}}&  &  O_{j_{a_n}} & O_{j_d} &
 \\
  & O_{j_{\overline{c}d}}  \arrow[from=lu, "\iota_2"] \arrow[from=ru, "\iota_2" description] \arrow[from=rru, "\iota_2" description]  &
 \end{tikzcd}
\end{center}

\end{lemma}

\begin{proof}
We recall that the inclusions $\iota_1$ and $\iota_2 $ belong to commuting systems. Thus, commutativity of the diagram above is equivalent to that of this one :

\begin{center}
\begin{tikzcd}
&&& O_{i_{{\overline{c}d}}} \arrow[ddddd,dotted, "\exists ? \, g"]\arrow[from=lld, "\iota_1"] \arrow[from=rdd, "\iota_1" description]

\\
 & O_{i_{{\overline{c}}}}  \arrow[ddd, "f"]\arrow[from=ld, "\iota_1"] \arrow[from=rd, "\iota_1"] & && & 

 \\
 O_{i_{a_1}}\arrow[d, "m_{a_1}"]&&  O_{i_{a_n}}\arrow[d, "m_{a_n}"]   && O_{i_{d}}\arrow[d, dotted, "\exists ? \, m_{d}"] & 
 \\
 O_{j_{a_1}}&  &  O_{j_{a_n}} && O_{j_{d}} & 
 \\
  & O_{j_{{\overline{c}}}} \arrow[from=lu, "\iota_2"] \arrow[from=ru, "\iota_2"] &&&&
  \\
   &&& O_{j_{\overline{c}d}} \arrow[from=ull, "\iota_2"] \arrow[from=ruu, "\iota_2" description]
 \end{tikzcd}
\end{center}
Now, since $H$ is a morphism of simplicial groupoids, by condition (A) (see the proof of lemma \ref{fill_square}), there exists a higher dimension morphism $g$ making the following square commute :

\begin{center}
\begin{tikzcd}
 & O_{i_{{\overline{c}d}}}\arrow[ddd, dotted, "\exists \, g"] \arrow[from=ld, "\iota_1"]  &

 \\
 O_{i_{{\overline{c}}}}\arrow[d, "f"]&&  
 \\
 O_{j_{{\overline{c}}}}&  & 
 \\
  & O_{j_{\overline{c}d}}  \arrow[from=lu, "\iota_2"]  &
 \end{tikzcd}
\end{center}

Applying condition (A) again, this time to the rightmost part of the big diagram, we find an appropriate morphism $m_{d}$. 

\end{proof}

\begin{prop}\label{coherent_families}
Let $\mathcal{G}_1$, $\mathcal{G}_2$ be 0-definable simplicial groupoids over $A$ with the disjoint union property. Let $H : \mathcal{G}_1 \rightarrow \mathcal{G}_2$ be a morphism of simplicial groupoids.
Let $\iota_1= ( O_{i_{\overline{c}}} \rightarrow O_{i_{\overline{d}}})_{\overline{c} \subset \overline{d} \in A^{<\omega}}$ and $\iota_2 = (O_{j_{\overline{c}}} \rightarrow O_{j_{\overline{d}}})_{\overline{c} \subset \overline{d} \in A^{<\omega}} $ be commuting systems of inclusions in $\mathcal{G}_1$ and $\mathcal{G}_2$ respectively. 

Then there exists a coherent family of maps $(m_a : O_{i_a} \rightarrow O_{j_a})_{a\in A}$ belonging to the morphism of simplicial groupoids $H$.

More generally, for all finite $B \subseteq A$, any \textit{coherent} family of maps $(m_a : O_{i_a} \rightarrow O_{j_a})_{a\in B}$ can be extended into a coherent family of maps $(m_a : O_{i_a} \rightarrow O_{j_a})_{a\in A}$.
\end{prop}

\begin{proof}

We shall build the additional morphisms inductively, using the previous lemma.

Let us take an enumeration of $A$ : $A = \lbrace a_{\alpha} \, | \, \alpha < \kappa \rbrace $ that begins with the elements of $B$. We will build inductively a family of degree 1 morphisms $m_{\alpha} : O_{i_{a_{\alpha}}} \rightarrow O_{j_{a_\alpha}}$ satisfying the following coherence property :

For all $\beta < \kappa$, for all finite tuples $\overline{c} = (a_1,...,a_n) \subseteq \lbrace a_{\alpha} \, | \, \alpha \leq  \beta \rbrace$, there exists a map $f$ belonging to $H$, in degree $|\overline{c}|$, making the following diagram commute :

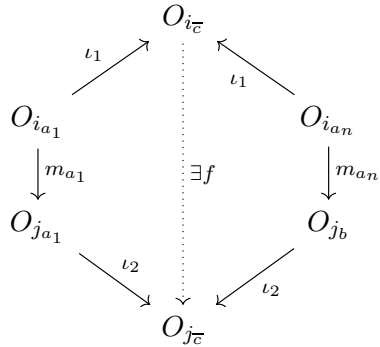
\begin{center}
\begin{tikzcd}
 & O_{i_{{\overline{c}}}}\arrow[ddd, dotted, "\exists f"] \arrow[from=ld, "\iota_1"] \arrow[from=rd, "\iota_1"] &

 \\
 O_{i_{a_1}}\arrow[d, "m_{a_1}"]&&  O_{i_{a_n}}\arrow[d, "m_{a_n}"]   
 \\
 O_{j_{a_1}}&  &  O_{j_b}
 \\
  & O_{j_{{\overline{c}}}}  \arrow[from=lu, "\iota_2"] \arrow[from=ru, "\iota_2"] &
 \end{tikzcd}
\end{center}

We start with the construction up to the finite ordinal $|B|$.

Let $\beta$ be an ordinal smaller than $\kappa$. Given the construction for all ordinals $ \alpha < \beta$, we shall extend it by realizing a partial type whose parameters are the morphisms already built.
Let $p_\beta (m)$ be the set of formulas containing the formula \textquotedblleft $m \in \mathrm{Hom}_{\mathcal{G}}(O_{i_{a_\beta}}, O_{j_{a_\beta}})"$ and, for each tuple $\overline{c} = (a,...,b) \subseteq \lbrace a_{\alpha} \, | \, \alpha \leq \beta  \rbrace$, the formula expressing :

\textquotedblleft There exists a higher dimension morphism $g : O_{i_{\overline{c}}} \rightarrow O_{j_{\overline{c}}}$ in $H$ making the following diagram commute : 
\begin{center}
\begin{tikzcd}
 & O_{i_{\overline{c}}}\arrow[ddd, "g"] \arrow[from=ld, "\iota_1"] \arrow[from=rd, "\iota_1"]  &

 \\
 O_{i_{a_1}}\arrow[d, "m_{a_1}"]&&  O_{i_{b}}\arrow[d, "m_{a_n}"]   & 
 \\
 O_{j_{a_n}}&  &  O_{j_b} & 
 \\
  & O_{j_{\overline{c}}}  \arrow[from=lu, "\iota_2" ] \arrow[from=ru, "\iota_2" description]  & && 
 \end{tikzcd} where $m_{a_{\beta}}$ is a notation for $m$."
\end{center}

Using the previous lemma and the induction hypothesis, we see that this set of formulas is consistent : if $\overline{c}_1,...,\overline{c}_m$ are finite tuples, by the commutativity of the systems of inclusions $\iota_1$ and $\iota_2$, and condition (A), it suffices to find a coherence morphism for the larger tuple $\overline{c}_1...\overline{c}_n$. The previous lemma shows that such a coherence morphism exists.

Moreover, the set of parameters of $p_\beta(m)$ is small enough in cardinality. Thus, by saturation, this partial type is realized in $\mathbb{U}$.
\end{proof}

Having introduced these two notions of coherence, we can now define the cover associated to a simplicial groupoid :

\begin{theo}\label{cover_from_simplicial_groupoid}
Let $\mathcal{G}$ be a 0-definable simplicial groupoid in $\mathbb{U}$ over $A$, with the disjoint union property. Then there exists a cover $(\mathbb{U}, O_*)$ of $\mathbb{U}$, along with an extension $\mathcal{G}'$ of $\mathcal{G}$ that is 0-definable in $(\mathbb{U}, O_*)$, such that :

\begin{itemize}

    \item The cover $(\mathbb{U}, O_*)$ is 1-analysable over $A$, and satisfies the two technical hypotheses given at the beginning of subsection \ref{subsection_simplicial_from_cover}.
    \item The simplicial groupoid $\mathcal{G}'$ is an extension of $\mathcal{G}$ with the extra object $O_{*,\overline{c}}$ in the connected component over $\overline{c}$, for every $\overline{c}$ in $A^{< \omega}$.
    \item The automorphism groups $\mathrm{Aut}_{\mathcal{G}'}(O_{*,\overline{c}})$, along with their action on $O_{*,\overline{c}}$, are isomorphic to the groups $\mathrm{Aut}(O_{*,\overline{c}}/\mathbb{U})$.
\end{itemize}

\end{theo}

\begin{proof}

Pick a commuting system of inclusions $\iota : O_{i_{\overline{c}}} \hookrightarrow O_{i_{\overline{d}}}$, in the groupoid morphisms defining the simplicial groupoid, where $\overline{c}\subseteq \overline{d}$ are finite tuples of elements of $A$.

Create a copy $O_{*,a}$ of $O_{i_a}$, with a map $f_a : O_{*,a} \rightarrow O_{i_a}$ identifying the two sets. The new objects in degree $n$ are the $O_{*, a_1}\cup ... \cup O_{*, a_n}$. Use the inclusions to build bijections ${f_{a_1...a_n}} : O_{*, a_1}\cup ... \cup O_{*, a_n} \simeq O_{i_{a_1...a_n}}$. These bijections are morphisms in the extended groupoid. They make the following diagrams commute : 

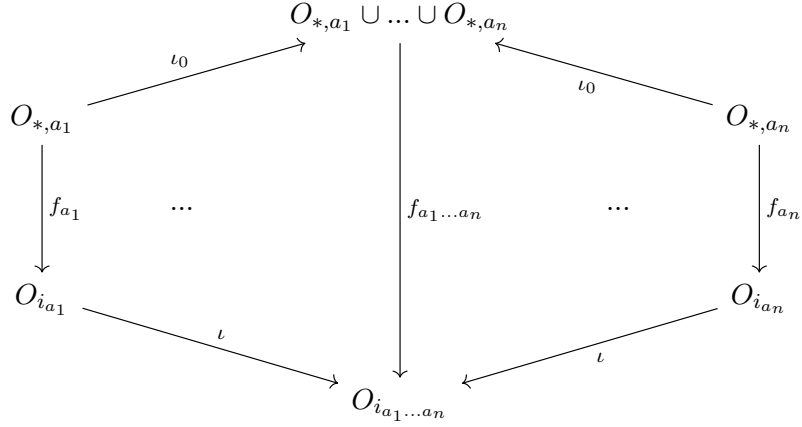
\begin{center}
\begin{tikzcd}
 && O_{*, a_1}\cup ... \cup O_{*, a_n}\arrow[dddd, "f_{a_1...a_n}"] \arrow[from=lld, "\iota_0"] \arrow[from=rrd, "\iota_0"] &

 \\
 O_{*,{a_1}}\arrow[dd, "f_{a_1}"]&&&&  O_{*,{a_n}}\arrow[dd, "f_{a_n}"]   
 \\
 & ... && ...
 \\
 O_{i_{a_1}}& && &  O_{i_{a_n}}
 \\
  && O_{i_{a_1...a_n}}  \arrow[from=llu, "\iota"] \arrow[from=rru, "\iota"] &
 \end{tikzcd} 
\end{center}
Here, the $\iota_0$ are the set-theoretic inclusions, and the $\iota$ are in the commuting system picked earlier. Note that we only consider the case of pairwise distinct $a_i$.

Extend the morphisms of groupoids, i.e. inclusions, in the only way that makes the maps $f_{a_1...a_n}$, $n \geq 1, a_1,...,a_n \in A$ isomorphisms of these groupoids. 

Define a multi-sorted extension of $\mathbb{U}$, using the whole simplicial groupoid structure : There are new sorts $O^{}_* = \bigcup\limits_{a \in A} O_{*, a}$ and $M^{(n)}_*$ for the groupoid in degree $n$, and new sorts $N^{(i,j)}_*$ for the morphisms from the degree $i$ groupoid to the degree $j$ groupoid. The new relations are the ones interpreting the morphisms in the new sorts as bijections.

We call $O_*$ the object sort.

We define the structure $\mathbb{U}'$ as the induced structure on $(\mathbb{U}, O_*)$.

We shall now check that this structure satisfies all the required properties. The final result will be Proposition \ref{local_embededness}.

\end{proof}
\begin{prop}\label{first_prop_cover}
The extension of $\mathcal{G}$ in the structure $(\mathbb{U}, O^{}_*, M^{(n)}_*, N^{(i,j)}_*)$ is a simplicial groupoid over $A$.
\end{prop}
\begin{proof}
First, in each degree, the extra objects $O_{*, \overline{c}}$ and the bijections $f_{\overline{c}} : O_{*, \overline{c}} \rightarrow O_{i_{\overline{c}}}$ define groupoids, just as in the internal or independent cases.

Then, we need to check that the extended inclusions satisfy condition (A). We know that the maps $f_{\overline{c}}$ connect the inclusions belonging to $\iota$ to the set-theoretic ones in $\iota_0$. Then, since the inclusion morphisms of groupoids in $\mathcal{G}$ satisfy condition (A), the extended ones also do.

Finally, let us show that all this data defines a simplicial groupoid, i.e. all the required relations between the inclusions hold. First, they hold for the set-theoretic maps. Thus, by condition (A), they hold globally for the sets of maps that define the morphisms of groupoids.
\end{proof}

\begin{prop}
The structure $(\mathbb{U}, O_*)$ is stably embedded in the simplicial groupoid structure $(\mathbb{U}, O^{}_*, M^{(n)}_*, N^{(i,j)}_*)$.
\end{prop}

\begin{proof}
Using finite faithfulness, we see that the morphism sorts are included in the definable closure, in $(\mathbb{U}, O^{(n)}_*, M^{(n)}_*, N^{(i,j)}_*)$, of the object sorts.
\end{proof}

Now, to prove stable embeddedness of $\mathbb{U}$ in this new structure, we shall proceed more carefully than in the case of independent fibers. Just as in proposition \ref{cover_in_independent_case}, if $\sigma$ is an automorphism of $\mathbb{U}$, in order to define its action on the morphisms in the extended groupoid, we need to find maps $m_a : O_{ i_{\sigma(a)}} \rightarrow O_{\sigma(i_a)} $ in the degree 1 groupoid. Then, we will be able to extend $\sigma$ as follows: if, for instance, $m : O_{*, a} \rightarrow O_j$ is a degree 1 morphism in the extended groupoid, we define $\sigma(m)$ as the morphism $\sigma(m f_a^{-1}) m_a f_{\sigma(a)} : O_{*, \sigma(a)} \rightarrow O_{\sigma(j)}$. 

However, this time we need to make sure our choices of maps come from morphisms in the higher degree groupoids.
This is necessary in order to define the action of $\sigma$ on higher degree morphisms in a way that is compatible with the inclusions.

\begin{prop}\label{cover_general_case}
The structure $(\mathbb{U}, O_*)$ is a 1-analysable cover of $\mathbb{U}$.
\end{prop}

\begin{proof}
We wish to extend automorphisms of $\mathbb{U}$ using the isomorphisms $f_a^{-1} : O_{i_a} \simeq O_{*,a}$.

\medskip

Using proposition \ref{coherent_families}, we find morphisms $m_a :  O_{i_{\sigma(a)}} \rightarrow O_{\sigma(i_a)}$ in the original groupoid, with the following coherence property : 
For all finite tuples $\overline{c} = (a_1,...,a_n) \subseteq A$, there exists a morphism $f$ in the groupoid of degree $|\overline{c}|$ making the following diagram commute :

\begin{center}
\begin{tikzcd}
 & O_{i_{\sigma({\overline{c}})}}\arrow[ddd, dotted, "\exists f"] \arrow[from=ld, "\iota"] \arrow[from=rd, "\iota"] &

 \\
 O_{i_{\sigma(a)}}\arrow[d, "m_{a_1}"]&&  O_{i_{\sigma(b)}}\arrow[d, "m_{a_n}"]   
 \\
 O_{\sigma(i_a)}&  &  O_{\sigma(i_b)}
 \\
  & O_{\sigma(i_{{\overline{c}}})}  \arrow[from=lu, "\sigma(\iota)"] \arrow[from=ru, "\sigma(\iota)"] &
 \end{tikzcd}
\end{center}

We can thus extend the action of $\sigma$ to the degree 1 groupoid : if, for instance, $m : O_{*, a} \rightarrow O_j$ is a degree 1 morphism in the extended groupoid, we define $\sigma(m)$ as the morphism $\sigma(m f_a^{-1}) m_a f_{\sigma(a)} : O_{*, \sigma(a)} \rightarrow O_{\sigma(j)}$.

\noindent Then, using the coherence property of the family of maps we built, we find uniquely determined maps $m_{\overline{c}} : O_{i_{\sigma(\overline{c})}} \rightarrow O_{\sigma(i_{\overline{c}})}$, $c \in A^{< \omega}$ lifting the degree 1 morphisms $m_a$. In fact, uniqueness comes from the disjoint union property proved in \ref{remark_product}, and forces commutativity of the coherence diagrams between these maps. Using these higher degree bijections, we define $\sigma$ on the higher degree groupoids in a similar way. See proposition \ref{cover_in_independent_case} for the full definition.

\noindent It remains to define the action of $\sigma$ on the inclusions. The distinguished \textquotedblleft set-theoretic" inclusions $O_{*,a_1} \cup ... \cup O_{*, a_{n}} \rightarrow O_{*,a_1} \cup ... \cup O_{*, a_{n+1}}$ are sent to their counterparts $O_{*, \sigma(a_1)} \cup ... \cup O_{*, \sigma(a_{n})} \rightarrow O_{*,\sigma(a_1)} \cup ... \cup O_{*, \sigma(a_{n+1})}$.

To check that this defines an automorphism of the extended simplicial groupoid, viewed as a first-order multi-sorted structure, we use  the 0-definability in $\mathbb{U}$ of the smaller simplicial groupoid structure, the fact $\sigma$ is an automorphism of $\mathbb{U}$, and the coherence property described above :

\noindent As a first step, let's check that, if the left square of the following diagram commutes, then so does the right square :

\begin{center}
\begin{tikzcd} 
O_{2}\arrow[r, "n"] \arrow[from=d, "\iota'"]& O_{*,ab} \arrow[from=d, "\iota_0"]&& \sigma(O_{2})\arrow[r, "\sigma(n)"] \arrow[from=d, "\sigma(\iota')"]&O_{*,\sigma(ab)} \arrow[from=d, "\sigma(\iota_0)"] \\

 O_1\arrow[r, "n_1"]&O_{*,a}  && \sigma(O_1)\arrow[r, "\sigma(n_1)"]& O_{*,\sigma(a)}
 \end{tikzcd}
\end{center}

Here, $\iota_0$ is the distinguished set-theoretic inclusion, whereas $\iota'$ is some inclusion in the original simplicial groupoid. 
We recall that $\sigma(n) = f_{\sigma(ab)}^{-1}m_{ab}^{-1}\sigma(f_{ab}n)$ and $\sigma(n_1) = f_{\sigma(a)}^{-1}m_{a}^{-1}\sigma(f_{a}n_1)$. Thus we can compute : $$\sigma(\iota_0)\sigma(n)
=\sigma(\iota_0) f_{\sigma(ab)}^{-1}m_{ab}^{-1}\sigma(f_{ab}n)$$

$$ = \sigma(\iota_{0}) f_{\sigma(ab)}^{-1}m_{ab}^{-1}\sigma(f_{ab}n)$$

$$ = f_{\sigma(a)}^{-1} \iota m_{ab}^{-1}\sigma(f_{ab}n)$$

$$ = f_{\sigma(a)}^{-1} m_a^{-1} \sigma(\iota) \sigma(f_{ab}n) $$

$$ = f_{\sigma(a)}^{-1} m_a^{-1} \sigma(\iota f_{ab} n )$$

$$ = f_{\sigma(a)}^{-1} m_a^{-1} \sigma(f_{a}\iota_0 n ) $$
$$ = f_{\sigma(a)}^{-1} m_a^{-1} \sigma(f_{a} n_1 \iota' )$$

$$ = f_{\sigma(a)}^{-1} m_a^{-1} \sigma(f_{a} n_1) \sigma(\iota') $$

$$ = \sigma(n_1) \sigma(\iota').$$

In other words, the following diagram commutes :

\begin{center}
    \begin{tikzcd}
    && O_{i_{\sigma(a)}} \arrow[r, "f_{\sigma(a)}^{-1}"] & O_{*, \sigma(a)}
    \\
    && O_{i_{\sigma(ab)}} \arrow[r, "f_{\sigma(ab)}^{-1}"] \arrow[from=u, "\iota"] & O_{*, \sigma(ab)} \arrow[from=u, "\sigma(\iota_0)"]
    \\
    O_{\sigma(i_a)} \arrow[rruu, "m_a^{-1}"] && O_{\sigma(i_{ab})} \arrow[u, "m_{ab}^{-1}"]\arrow[from=ll, "\sigma(\iota)" description]
    \\
    && 
    \\
     \sigma(O_1)\arrow[uu, "\sigma(f_a n_1)"]\arrow[uurr, "\sigma(f_{ab}\iota_0 n_1)" description] && \arrow[from=ll, "\sigma(\iota')"] \sigma(O_2) \arrow[uu, "\sigma(f_{ab} n)" description]

    \end{tikzcd}
\end{center}

The bottom part of the diagram commutes because $\sigma$ is both an automorphism of the disjoint union of the extended groupoids $\sqcup_n\mathcal{G'}_n$ and an automorphism of the simplicial groupoid $\mathcal{G}$, and because the following diagram commutes : 

\begin{center}
    \begin{tikzcd}

    O_{i_a}&& O_{i_{ab}} \arrow[from=ll, "\iota" description]
    \\
    O_{*,a}\arrow[u, "f_a"] && \arrow[from=ll, "\iota_0"] O_{*, ab}\arrow[u, "f_{ab}" description]
    \\
     O_1\arrow[u, "n_1"] && \arrow[from=ll, "\iota'"] O_2 \arrow[u, "n" description]

    \end{tikzcd}
\end{center}

Similar computations deal with the cases where the objects $O_1$, $O_2$ may be equal to the new objects $O_{*,a}$ or $O_{*, ab}$. In these cases, we use the description of the automorphisms of the $O_{*, \overline{c}}$ as the conjugates of the automorphisms of $O_{i_{\overline{c}}}$ by the maps $f_{\overline{c}}$.

\noindent Therefore, defining $\sigma$ on the new inclusions with $\sigma(\alpha\iota'\beta) := \sigma(\alpha) \sigma(\iota') \sigma(\beta)$ yields a well-defined automorphism of the structure composed of the family of the extended groupoids $\mathcal{G}'_n$ along with the extended inclusion morphisms between them. 

\end{proof}

\begin{rem}
A cautious reader might object that we need to take saturated models of the theories we define, before we can use automorphism arguments to prove that we have defined a cover. This issue can be addressed by noticing that the inductive constructions that constitute the core of the proofs only use condition (A) for the various morphisms of groupoids involved, which is preserved under elementary equivalence.
\end{rem}

\begin{prop}
The cover built from the simplicial groupoid $\mathcal{G}$ does not depend on the choices made.
\end{prop}

\begin{proof}
Let us assume we had picked another commuting system of inclusions $\iota'$, with objects $O_{j_{\overline{c}}}$, and copying maps $f'_{\overline{c}} : O_{j_{\overline{c}}} \rightarrow O'_{*, \overline{c}}$, for $\overline{c}$ in $A^{< \omega}$. Using proposition \ref{coherent_families}, we find bijections $g_{\overline{c}} : O_{i_{\overline{c}}} \rightarrow O_{j_{\overline{c}}}$ in the simplicial groupoid $\mathcal{G}$ that make the following diagrams commute : 

\begin{center}
    \begin{tikzcd}
    O_{i_{\overline{c}}} \arrow[from=d, "\iota"]\arrow[r, "g_{\overline{c}}"]& O_{j_{\overline{c}}}\arrow[from=d, "\iota'"]
    \\
    O_{i_a} \arrow[r, "g_a"]& O_{j_a}
    \end{tikzcd}
\end{center}

Now, we define an isomorphism $h$ between the two extended simplicial groupoids, using the only maps $h_{\overline{c}}$ that make the following diagrams commute :

\begin{center}
    \begin{tikzcd}
    O_{i_{\overline{c}}} \arrow[d, "f_{\overline{c}}"]\arrow[r, "g_{\overline{c}}"]& O_{j_{\overline{c}}}\arrow[d, "f'_{\overline{c}}"]
    \\
    O_{*, \overline{c}} \arrow[r, "h_{\overline{c}}"]& O'_{*, \overline{c}}
    \end{tikzcd}
\end{center}

This isomorphism of simplicial groupoids defines an isomorphism between the structures $(\mathbb{U}, O_*)$ and $(\mathbb{U}, O'_*)$, that fixes $\mathbb{U}$ pointwise.

\end{proof}

\begin{prop}\label{binding_groupoid_statement}
Let $\overline{c} \in A^{=n}$ be a finite tuple. Then the group $\mathrm{Aut}(O_{*, \overline{c}} / \mathbb{U})$ along with its action on $O_{*, \overline{c}}$, is isomorphic to the group $\mathrm{Aut}_{\mathcal{G}'}(O_{*, \overline{c}}).$
\end{prop}

\begin{proof}
Let $\sigma$ be an automorphism in $\mathrm{Aut}(O_{*, \overline{c}} / \mathbb{U})$. If $m \in \mathrm{Hom}_{\mathcal{G}'}(O_{*, \overline{c}}, O_{i_{\overline{c}}})$ and $x \in O_{*, \overline{c}}$, we compute : $m(x) = \sigma(m(x)) = \sigma(m) (\sigma(x))$. The first equality comes from the fact $m(x) \in \mathbb{U}$. The second one comes from the 0-definability of the action of the morphisms, which is therefore preserved by $\sigma$. Thus $\sigma|_{O_{*, \overline{c}}} = \sigma(m)^{-1} \circ m \in \mathrm{Aut}_{\mathcal{G}'}(O_{*, \overline{c}}).$

For the converse, we use proposition \ref{coherent_families}. We are given an element $\sigma_{\overline{c}}$ of $\mathrm{Aut}_{\mathcal{G}'}(O_{*, \overline{c}})$, and wish to extend it into an automorphism of the whole structure $(\mathbb{U}, O^{(n)}_*, M^{(n)}_*, N^{(i,j)}_*)$ fixing $\mathbb{U}$ pointwise. We shall find morphisms $\sigma_{\overline{d}} \in \mathrm{Aut}_{\mathcal{G}'}(O_{*, \overline{d}})$ for each finite tuple $\overline{d}$ in $A$, with the following coherence condition : if $\overline{a}$ is a subtuple of $\overline{b}$, then the following diagram commutes :

\begin{center}
\begin{tikzcd} O_{*, \overline{b}}\arrow[r, "\sigma_{\overline{b}}"] \arrow[from=d, "\iota"]& O_{*, \overline{b}} \arrow[from=d, "\iota"] \\
 O_{*, \overline{a}}\arrow[r, "\sigma_{\overline{a}}"]& O_{*, \overline{a}}   \end{tikzcd}
\end{center}

\noindent For instance, let's assume that we start with an automorphism $\sigma_{abc}$ of the object $O_{*,abc}$. The degree 1 morphisms $\sigma_a, \sigma_b, \sigma_c$ are the ones uniquely determined by the commutativity of the following diagrams :

\begin{center}
\begin{tikzcd} 
O_{*, abc}\arrow[r, "\sigma_{abc}"] \arrow[from=d, "\iota"]& O_{*, abc} \arrow[from=d, "\iota"] && O_{*, abc}\arrow[r, "\sigma_{abc}"] \arrow[from=d, "\iota"]& O_{*, abc} \arrow[from=d, "\iota"]  && O_{*, abc}\arrow[r, "\sigma_{abc}"] \arrow[from=d, "\iota"]& O_{*, abc} \arrow[from=d, "\iota"]

\\
 O_{*, a}\arrow[r, "\sigma_{a}"]& O_{*, a}   && O_{*, b}\arrow[r, "\sigma_{b}"]& O_{*, b} && O_{*, c}\arrow[r, "\sigma_{c}"]& O_{*, c} 
 
 \end{tikzcd}
\end{center}
By proposition \ref{coherent_families}, this finite coherent family can be extended into a full coherent family $(\sigma_a)_{a \in A}$ of degree 1 automorphisms of the $O_{*,a}$. Because of the disjoint union property, the higher degree automorphisms $\sigma_{\overline{c}}$ are uniquely determined by the $\sigma_a$.

Finally, the action on morphisms is defined either as precomposition with the appropriate $\sigma_{\overline{c}}$, postcomposition with the inverse of such a map, conjugation by it, or identity, depending on the domain and codomain of the morphism.
\end{proof}

\begin{coro}
Let $n$ be a positive integer. Then the language of $(\mathbb{U}, \bigcup\limits_{a \in \overline{c}} O_{*,a})$ is finitely generated over that of $\mathbb{U}$, uniformly for $\overline{c}$ in $A^{=n}$.
\end{coro}

\begin{proof}
We use proposition \ref{criterion_uniformly_finite_language}. The automorphism groups $\mathrm{Aut}(\bigcup\limits_{a \in \overline{c}} O_{*,a} / \mathbb{U})$ are uniformly interpretable in $\mathbb{U}'$, since they are defined in the groupoid of degree $n$. 
\end{proof}

\begin{prop}\label{local_embededness}
Let $C \subseteq A$ be a finite subset. Then the structure $(\mathbb{U}, \bigcup\limits_{a \in C}O_{*,a})$ is stably embedded in  $(\mathbb{U}, O_*)$.
\end{prop}
\begin{proof}
The proof is similar to that of propositions \ref{cover_general_case} and \ref{binding_groupoid_statement}. The automorphisms $\sigma$ we begin with are already defined on some of the new objects and morphisms. For instance, we are already given morphisms $m_a = \sigma|_{O_{*,a}} : O_{*,a} \rightarrow O_{*, a}$ for $a$ in $C$. In fact, by finite faithfulness, we can extend $\sigma$ in a unique way to the simplicial groupoid restricted to $C$. It suffices to check that these maps define a coherent family, with respect to the set-theoretic inclusions, and to apply proposition \ref{coherent_families}. 
\end{proof}

\end{subsection}

\begin{subsection}{Comparing covers}

From now on, we let $\mathbb{V} = (\mathbb{U}, S, p : S\rightarrow A)$ be a 1-analysable cover satisfying the  conditions given at the beginning of subsection \ref{subsection_simplicial_from_cover}, i.e., (Local Stable Embeddedness) and (Uniform Finite Generatedness). The simplicial groupoids $\mathcal{G}(\mathbb{V})$ and $\mathcal{G}'(\mathbb{V})$ are the binding simplicial groupoids definable in $\mathbb{U}$ and $\mathbb{V}^{eq}$ respectively.
As in the case of independent fibers, we shall also work with the cover $C\mathcal{G}(\mathbb{V}) = (\mathbb{U}, O_*)$ and the simplicial groupoids $\mathcal{G}(C\mathcal{G}(\mathbb{V}))$ and $\mathcal{G}''(C\mathcal{G}(\mathbb{V}))$.

\begin{theo}\label{theo_comparing_covers}
There exists an isomorphism of covers $\eta_{\mathbb{V}} : \mathbb{V} \simeq C\mathcal{G}(\mathbb{V}).$
\end{theo}

\begin{proof}
Pick a coherent systems of inclusions $\iota = (O_{i_{\overline{c}}} \rightarrow O_{i_{\overline{d}}})$ in $\mathcal{G}(\mathbb{V})$. Then pick families of morphisms $(g_a : S_a \rightarrow O_{i_a})_a$ and $(h_a : O_{i_a} \rightarrow O_{*,a})_a$ that belong to $\mathcal{G}'(\mathbb{V})$ and $\mathcal{G}''(C\mathcal{G}(\mathbb{V}))$ respectively, and which are \textit{coherent}, with respect to $\iota$ and the set-theoretic inclusions. By proposition \ref{coherent_families}, such families exist. We shall prove that the structure $(\mathbb{V}\cup C\mathcal{G}(\mathbb{V}), \bigcup\limits_{a \in A} h_a g_a)$ defines an isomorphism of covers.

First, let us show that $\mathbb{U}$ is stably embedded in this structure. Let $\sigma$ be an automorphism of $\mathbb{U}$. We extend it \textit{arbitrarily} into an automorphism of $\mathbb{V}\cup C\mathcal{G}(\mathbb{V})$. Then, the maps $\sigma(h_a)$ and $\sigma(g_a)$ are morphisms in $\mathcal{G}''(C\mathcal{G}(\mathbb{V}))$ and $\mathcal{G}'(\mathbb{V})$ respectively. Since the set theoretic inclusions are 0-definable, they are preserved by $\sigma.$ Thus, the family of maps $(\sigma(h_a) \sigma(g_a))$ is coherent with respect to the set-theoretic inclusions in $S$ and $O_*$, just as $(h_a g_a)$ is. As a consequence, the family of maps $(\sigma(g_a^{-1})\sigma(h_a^{-1}) h_{\sigma(a)} g_{\sigma(a)})$ is also coherent. To sum up, for all finite tuples $\overline{c} = (a,...,b) \subseteq A$ made of pairwise distinct elements, there exist higher degree morphisms $g,g'$ in $\mathcal{G}'(\mathbb{V})$ and $h,h'$ in $\mathcal{G}''(C\mathcal{G}(\mathbb{V}))$ making the following diagram commute :

\begin{center}
    \begin{tikzcd}
    S_{\sigma(a)} \arrow[r,  " g_{\sigma(a)}"] & O_{i_{\sigma(a)}} \arrow[r,  " h_{\sigma(a)}"]& O_{*, \sigma(a)}\arrow[r,  "\sigma(h_a^{-1}) "]& O_{\sigma(i_a)} \arrow[r,  "\sigma(g_a^{-1}) "]& S_{ \sigma(a)}
    \\
    S_{\sigma(\overline{c})} \arrow[u, hook] \arrow[d, hook]\arrow[r, dotted, "\exists g"]& O_{i_{\sigma(\overline{c})}} \arrow[d, hook, "\iota"]\arrow[u, hook, "\iota"] \arrow[r, dotted, "\exists h"]& O_{*, \sigma(\overline{c})}\arrow[u, hook]\arrow[d, hook] \arrow[r, dotted, "\exists h'^{-1}"]& O_{\sigma(i_{\overline{c}})} \arrow[d, hook, "\sigma(\iota)"]\arrow[u, hook, "\sigma(\iota)"] \arrow[r, dotted, "\exists g'^{-1}"]& S_{\sigma(\overline{c})}\arrow[u, hook]\arrow[d, hook]
    \\
    S_{\sigma(b)} \arrow[r, "g_{\sigma(b)}"] & O_{i_{\sigma(b)}} \arrow[r, "h_{\sigma(b)}"]& O_{*, \sigma(b)} \arrow[r,  "\sigma(h_b^{-1}) "]& O_{\sigma(i_b)} \arrow[r, "\sigma(g_b^{-1})"]& S_{\sigma(b)}

    \end{tikzcd}
    
\end{center}

\noindent In fact, the maps $g$ and $h$ come from the coherence of the families $(g_a)$ and $(h_a)$. Besides, the families $(h_a^{-1})$ and $(g_a^{-1})$ are also coherent, with respect to the set-theoretic inclusions and $\iota$. Applying the model-theoretic automorphism $\sigma$ to the commutative diagrams shows that the families $(\sigma(h_a^{-1}))$ and $(\sigma(g_a^{-1}))$ are coherent with respect to the set-theoretic inclusions and $\sigma(\iota)$. This implies the existence of the maps $g'$ and $h'$.

\noindent Now, let $\tau$ be the map $\bigcup\limits_{a \in A} \sigma(g_a^{-1})\sigma(h_a^{-1}) h_{\sigma(a)} g_{\sigma(a)}$. Using lemma \ref{criterion_automorphism}, we see that $\tau \cup id_{\mathbb{U}}$ is an automorphism of $\mathbb{V}$ that fixes $\mathbb{U}$ pointwise, and sends the family of maps $(h_{\sigma(a)} g_{\sigma(a)})$ to the family of maps $(\sigma(h_a)\sigma(g_a))$. Finally, the map $(\tau \cup id_{\mathbb{U}})^{-1}\circ \sigma$ is an automorphism of $(\mathbb{V}\cup C\mathcal{G}(\mathbb{V}), \bigcup\limits_{a \in A} h_a g_a)$ extending $\sigma$. 

\noindent Thus, $\mathbb{U}$ is stably embedded in the structure $(\mathbb{V}\cup C\mathcal{G}(\mathbb{V}), \bigcup\limits_{a \in A} h_a g_a)$.

\medskip

Finally, we prove that $\mathbb{V}$ is stably embedded in $(\mathbb{V}\cup C\mathcal{G}(\mathbb{V}), \bigcup\limits_{a \in A} h_a g_a)$, the other stable embeddedness statement being similar. Let $\sigma$ be an automorphism of $\mathbb{V}$, which we may assume to fix $\mathbb{U}$ pointwise. We define the action of $\sigma$ on $O_{*}$ by conjugating $\sigma$ with the map $\bigcup\limits_{a \in A} h_a g_a$.
Let $\tau$ be the map $(\bigcup\limits_{a \in A} h_a g_a) \circ \sigma|_S \circ (\bigcup\limits_{a \in A} h_a g_a)^{-1}.$

\noindent From the binding groupoid statements in $\mathbb{V}$ and $C\mathcal{G}(\mathbb{V})$, and lemma \ref{criterion_automorphism}, this map $\tau$, when extended by the identity on $\mathbb{U}$, is an automorphism of the structure $C\mathcal{G}(\mathbb{V}).$ Thus, the map $\tau \cup \sigma$ is an automorphism of $(\mathbb{V}\cup C\mathcal{G}(\mathbb{V}), \bigcup\limits_{a \in A} h_a g_a)$ extending $\sigma$. 

\noindent So $\mathbb{V}$ is stably embedded in $(\mathbb{V}\cup C\mathcal{G}(\mathbb{V}), \bigcup\limits_{a \in A} h_a g_a)$.

\end{proof}

\begin{rem}
This theorem implies that any interpretation of the extended simplicial groupoid $\mathcal{G}''(C\mathcal{G}(\mathbb{V}))$ yields an interpretation of the cover $\mathbb{V}$ itself.

A more precise study of this problem could hypothetically yield criteria for the interpretability of the cover $\mathbb{V}$ in the original structure $\mathbb{U}$.

For instance, if one could, inside the groupoid $\mathcal{G}$, find canonical groupoids in each degree, and a 0-definable commuting system of inclusions between the unique objects of these canonical groupoids, one would be able to reconstruct the cover $C\mathcal{G}(\mathbb{V})$ inside $\mathbb{U}$, and thus interpret the cover $\mathbb{V}$ in the structure $\mathbb{U}$ itself.
\end{rem}

\end{subsection}

\begin{section}{Functoriality}\label{section_functoriality_finite_generatedness}

From now on, we let $AC_A$ denote the category of 1-analysable covers over $A$ that satisfy the conditions of (Uniform Finite Generatedness) and (Local Stable Embeddedness) given at the beginning of subsection \ref{subsection_simplicial_from_cover}. Note that an object in $AC_A$ is given by a cover of $\mathbb{U}$ with exactly one extra sort $S$ and a definable surjective map $S \rightarrow A$, that satisfy some properties.

On the other hand, let $SFCG_A$ denote the category of the finitely faithful simplicial groupoids over $A$ that satisfy the disjoint union property. See Definition \ref{defi_simplicial_groupoid}.

If $\mathcal{C}$ is a category, we let $Iso(\mathcal{C})$ denote the subcategory of $\mathcal{C}$ with the same objects, but where the morphisms are the isomorphisms in $\mathcal{C}$.

We recall that morphisms of simplicial groupoids are sequences of morphisms of groupoids that are compatible with the inclusion morphisms defining the simplicial structures. They are required to satisfy condition $(C)$ in each connected component. On the other hands, morphisms in $AC_A$ are morphisms of covers compatible with the surjective maps $S \rightarrow A$.

\begin{subsection}{The functor $G : AC_A \rightarrow SFCG_A$}

\begin{defi}
Let $h : \mathbb{U}_1 \rightarrow \mathbb{U}_2$ be a morphism of 1-analysable covers over $A$. We define the simplicial groupoid morphism $G(h)$ as follows : 
In degree $n$, the set of maps defining the groupoid morphism between the degree $n$ groupoids $\mathcal{G}(\mathbb{U}_1)_n$ and $\mathcal{G}(\mathbb{U}_2)_n$ is given by $G(h)_n:=\lbrace \alpha_2 \circ h_{\overline{c}} \circ \alpha_1 : O_1 \rightarrow O_2\, | \, \overline{c} \in A^{=n}, \, O_2 \in Ob(\mathcal{G}(\mathbb{U}_2)_n),\, O_1 \in Ob(\mathcal{G}(\mathbb{U}_1)_n), \, \alpha_1 \in \mathrm{Hom}_{\mathcal{G}'(\mathbb{U}_1)_n}(O_1, S_{1, \overline{c}}), \, \alpha_2 \in \mathrm{Hom}_{\mathcal{G}'(\mathbb{U}_2)_n}(S_{2, \overline{c}}, O_2 ) \rbrace$.
\end{defi}

\begin{prop}
For each integer $n$, the set of maps $G(h)_n$ is a 0-definable morphism of groupoids from $\mathcal{G}(\mathbb{U}_1)_n$ to $\mathcal{G}(\mathbb{U}_2)_n$, that only depends on the theory of $(\mathbb{U}_1 \cup \mathbb{U}_2, h)$.
\end{prop}

\begin{proof}
As in the independent case, if we pick another map $h' : \mathbb{U}_1 \rightarrow \mathbb{U}_2$ such that $(\mathbb{U}_1 \cup \mathbb{U}_2, h) \equiv (\mathbb{U}_1 \cup \mathbb{U}_2, h')$, then, taking saturated models, the structures are isomorphic. So there exists an automorphism $\sigma$ of $\mathbb{U}_1 \cup\mathbb{U}_2$ such that $\sigma h = h' \sigma$. By stable embeddedness of $\mathbb{U}_1$ in $(\mathbb{U}_1 \cup \mathbb{U}_2, h)$, we may assume that $\sigma$ is the identity on $\mathbb{U}_1$. So we get $\sigma h = h'$. The binding groupoid statements in $\mathbb{U}_2$ enable us to conclude that, in each connected component, the sets of maps are equal.

Now, we need to check that the sets of maps satisfy the morphism conditions. We shall use the results proved in the simpler cases. We notice that, for each tuple $\overline{c}$ in $A^{=n}$, the map $h_{\overline{c}} : S_{1, \overline{c}} \rightarrow S_{2, \overline{c}}$ defines a morphism of covers. Moreover, the connected components over $\overline{c}$ of the binding groupoids $\mathcal{G}'(\mathbb{U}_1)$ and $\mathcal{G}'(\mathbb{U}_2)$ are the binding groupoids of the covers $(\mathbb{U}, S_{1, \overline{c}}) $ and $(\mathbb{U}, S_{2, \overline{c}}) $. Now, these covers are internal, so we already know that the set of maps $\lbrace \alpha_2 \circ h_{\overline{c}} \circ \alpha_1 : O_1 \rightarrow O_2\, | \, O_2 \in Ob(\mathcal{G}(\mathbb{U}_2)_n),\, O_1 \in Ob(\mathcal{G}(\mathbb{U}_1)_n), \, \alpha_1 \in \mathrm{Hom}_{\mathcal{G}'(\mathbb{U}_1)_n}(O_1, S_{1, \overline{c}}), \, \alpha_2 \in \mathrm{Hom}_{\mathcal{G}'(\mathbb{U}_2)_n}(S_{2, \overline{c}}, O_2 ) \rbrace$ defines a morphism of groupoids.

\end{proof}

\begin{prop}\label{G(h)_is_a_simplicial_morphism}
The family of morphisms of groupoids $G(h)_n$ defines a morphism of simplicial groupoids.
\end{prop}

\begin{proof}
We need to check that the following diagrams commute, for each choice of inclusion $\iota$ :

\begin{center}
\begin{tikzcd} \mathcal{G}(\mathbb{U}_1)_{n+1}\arrow[r, "G(h)_{n+1}"] \arrow[from=d, "\iota"]& \mathcal{G}(\mathbb{U}_2)_{n+1} \arrow[from=d, "\iota"] && \\
 \mathcal{G}(\mathbb{U}_1)_{n}\arrow[r, "G(h)_{n}"]& \mathcal{G}(\mathbb{U}_2)_{n} && 
 \end{tikzcd}
\end{center}

It is easier to prove that the diagrams involving the extended groupoids commute :

\begin{center}
\begin{tikzcd} \mathcal{G}'(\mathbb{U}_1)_{n+1}\arrow[r, "G(h)_{n+1}"] \arrow[from=d, "\iota"]& \mathcal{G}'(\mathbb{U}_2)_{n+1} \arrow[from=d, "\iota"] && \\
 \mathcal{G}'(\mathbb{U}_1)_{n}\arrow[r, "G(h)_{n}"]& \mathcal{G}'(\mathbb{U}_2)_{n} && 
 
 \end{tikzcd}
\end{center}

Since all the sets of maps involved are morphisms of groupoids, we only need to check that, in each connected component, there exists a common map between the composite morphisms of groupoids. To simplify notations, we assume that the inclusion maps that define the morphism of groupoids $\iota$ are the set-theoretic inclusions.
Let $\overline{c}=(a_1,...,a_{n})$ be a tuple in $ A^{n}$ and $\overline{d}=(a_1,...,a_{n+1})$ be a tuple in $ A^{n+1}$ containing $\overline{c}$. We claim that the map $x \in S_{1, \overline{c}} \mapsto h(x) \in S_{2, \overline{d}} $ belongs to both $\iota \circ G(h)_{n}$ and $G(h)_{n+1} \circ \iota$.

\end{proof}

\begin{prop}
The action of $G$ on morphisms defines a functor $G : AC_A \rightarrow SFCG_A$.
\end{prop}
\begin{proof}
The fact that $G$ sends identities of $AC_A$ to identities of $SFCG_A$ is proved in the same way as in the independent case.

Now, it remains to check that finite composition is preserved by $G$. Working fiberwise, this follows from the case of internal covers.  
\end{proof}

\end{subsection}

\begin{subsection}{The functor $C : Iso(SFCG_A) \rightarrow AC_A$}

Let $H : \mathcal{G}_1 \rightarrow \mathcal{G}_2$ be an isomorphism in $SFCG_A$. We wish to define a morphism of covers $C(H) : C\mathcal{G}_1 \rightarrow C\mathcal{G}_2.$ As in the independent case, we shall use the extended groupoids $\mathcal{G}'_1$ and $\mathcal{G}'_2$ defined in the covers $C\mathcal{G}_1$ and $C\mathcal{G}_2$ respectively.

\begin{defi}
Let $\iota_1$, $\iota_2$ be commuting systems of inclusions in $\mathcal{G}_1$ and $\mathcal{G}_2$ respectively.
Let $(f_{\overline{c}})$, $(h_{\overline{c}})$, $(g_{\overline{c}})$ be families of maps that are coherent with respect to $\iota_1$, $\iota_2$ and the set-theoretic inclusions. In other words, the following diagrams commute : 
\begin{center}
    \begin{tikzcd}
    S_{1,a} \arrow[r, " f_a"] & O_{1,a} \arrow[r, "h_a"]& O_{2, a} \arrow[r, "g_a"] & S_{2,a}
    \\
    S_{1, \overline{c}} \arrow[from=u, hook] \arrow[from=d, hook]\arrow[r, "f_{\overline{c}}"]& O_{1, \overline{c}} \arrow[from=u, "\iota_1" description] \arrow[from=d,"\iota_1" description] \arrow[r, "h_{\overline{c}}"]& O_{2, \overline{c}}\arrow[from=u,"\iota_2" description]\arrow[from=d, "\iota_2" description] \arrow[r, "g_{\overline{c}}"] & S_{2, \overline{c}} \arrow[from=u, hook]\arrow[from=d, hook]
    \\
    S_{1,b} \arrow[r, "f_b"] & O_{1, b} \arrow[r, "h_b"]& O_{2, b} \arrow[r, "g_b"] & S_{2,b}

    \end{tikzcd}
    
\end{center} We define $C(H)$ as the theory of $(C\mathcal{G}_1 \cup C\mathcal{G}_2, \bigcup\limits_{a \in A} g_a h_a f_a).$
\end{defi}

First, as in the independent case, we need to make sure the choices made do not change the theory of $C(H)$.

\begin{lemma}
Different choices for the families of maps $(g_{\overline{c}}), (h_{\overline{c}}), (f_{\overline{c}}), \iota_1, \iota_2$ yield structures $(C\mathcal{G}_1 \cup C\mathcal{G}_2, \bigcup\limits_{a \in A} g_a h_a f_a)$ that are isomorphic over $\mathbb{U}$.
\end{lemma}

\begin{proof} We shall combine the ideas of the proof of proposition \ref{C(h)_free_choice} and the coherence conditions found in this section. Let $\iota'_1$ and $\iota'_2$ be different choices of commuting systems of inclusions, and $(f'_{\overline{c}}), (g'_{\overline{c}})$ and $(h'_{\overline{c}})$ be different choices of coherent families. We want to reduce to the case where only the $f'_{\overline{c}}$ differ. To do this, we need to look at the following commutative diagrams : 
\begin{center}
    
\begin{tikzcd}[row sep=scriptsize, column sep=scriptsize]
S_{1,\overline{c}}\arrow[from=dd]\arrow[rrr, "f_{\overline{c}}"]&&& O_{1, \overline{c}}   \arrow[from=dd, dotted] \arrow[rrrr, "h_{\overline{c}}"]& &  & & O_{2, \overline{c}} \arrow[drr, "g_{\overline{c}}"]  \arrow[from=dd, dotted] & &\\

&S_{1,\overline{c}} \arrow[rrr, crossing over, "f'_{\overline{c}}"]  && & O'_{1, \overline{c}} \arrow[rr, "h'_{\overline{c}}"] &  &O'_{2, \overline{c}}\arrow[ur] \arrow[rrr, crossing over, "g'_{\overline{c}}"]  & && S_{2, \overline{c}} \\

S_{1,a}\arrow[rrr, "f_{a}" description]&&& O_{1, a} \arrow[rrrr, "h_a"] & &  && O_{2, a}  \arrow[drr, "g_a"]  & &\\

&S_{1,a}\arrow[uu, crossing over] \arrow[rrr, "f'_a"] && & O'_{1, a} \arrow[rr, "h'_a"] \arrow[uu, dotted, crossing over]&  &O'_{2,a}\arrow[uu, dotted, crossing over]\arrow[ur] \arrow[rrr, "g'_a"] && & S_{2,a} \arrow[uu, crossing over]\\
\end{tikzcd}

\end{center}
\noindent Here, the vertical lines are inclusions, the non-dotted ones being the set-theoretic inclusions, the others being the $\iota_i$ or $\iota'_i$, $i=1,2$. 

For each integer $n$ and each tuple $\overline{c}$ in $A^{=n}$, condition (A) of the morphism $H_n$, and bijectivity of the maps in $H$, implies the existence of a unique morphism $\sigma_{\overline{c}} \in \mathrm{Aut}_{\mathcal{G}'_1}(S_{1,c}) = \mathrm{Aut}(S_{1, \overline{c}} / \mathbb{U})$ making the following diagram commute

\begin{center}
    
\begin{tikzcd}[row sep=scriptsize, column sep=scriptsize]
S_{1,\overline{c}}\arrow[from=dd]\arrow[rrr, "f_{\overline{c}}"]\arrow[dr, dotted, "\exists \, ! \,\sigma_{\overline{c}}"]&&& O_{1, \overline{c}} \arrow[dr, dotted, "\exists \, !"]   \arrow[rrrr, "h_{\overline{c}}"]& &  & & O_{2, \overline{c}} \arrow[drr, "g_{\overline{c}}"]   & &\\

&S_{1,\overline{c}} \arrow[rrr, crossing over, "f'_{\overline{c}}"] && & O'_{1, \overline{c}} \arrow[rr, "h'_{\overline{c}}"] &  &O'_{2, \overline{c}}\arrow[ur] \arrow[rrr, crossing over, "g'_{\overline{c}}"]  & && S_{2, \overline{c}} \\

S_{1,a}\arrow[rrr,  "f_a" description] &&& O_{1, a}  \arrow[dr, dotted, "\exists \, !"] \arrow[rrrr, "h_a"] & &  && O_{2, a}  \arrow[drr, "g_a"]  & &\\

&S_{1,a} \arrow[uu, crossing over]\arrow[rrr, "f'_a"] && & O'_{1, a} \arrow[rr, "h'_a"] &  &O'_{2,a}\arrow[ur] \arrow[rrr, "g'_a"] && & S_{2,a} \arrow[uu, crossing over]\\
\end{tikzcd}

\end{center}

Now, by uniqueness, the family $(\sigma_{\overline{c}})_{n \in \mathbb{N}, \, \overline{c} \in A^{n}}$ makes the following diagram commute, for all finite tuples $\overline{c}$ and for all $a$ in $\overline{c}$ :

\begin{center}
    
\begin{tikzcd}[row sep=scriptsize, column sep=scriptsize]
S_{1,\overline{c}}\arrow[rrr, "f_{\overline{c}}"]\arrow[dr, "\sigma_{\overline{c}}"]&&& O_{1, \overline{c}} \arrow[dr]   \arrow[rrrr, "h_{\overline{c}}"]& &  & & O_{2, \overline{c}} \arrow[drr, "g_{\overline{c}}"]   & &\\

&S_{1,\overline{c}} \arrow[rrr, crossing over, "f'_{\overline{c}}"] && & O'_{1, \overline{c}} \arrow[rr, "h'_{\overline{c}}"] &  &O'_{2, \overline{c}}\arrow[ur] \arrow[rrr, crossing over, "g'_{\overline{c}}"]  & && S_{2, \overline{c}} \\

S_{1,a}\arrow[uu, crossing over]\arrow[rrr,  "f_a" description] \arrow[dr, "\sigma_a"]&&& O_{1, a}  \arrow[dr] \arrow[rrrr, "h_a"] & &  && O_{2, a}  \arrow[drr, "g_a"]  & &\\

&S_{1,a} \arrow[uu, crossing over]\arrow[rrr, "f'_a"] && & O'_{1, a} \arrow[rr, "h'_a"] &  &O'_{2,a}\arrow[ur] \arrow[rrr, "g'_a"] && & S_{2,a} \arrow[uu, crossing over]\\
\end{tikzcd}

\end{center}

Finally, using the criterion given in lemma \ref{criterion_automorphism}, we find that the map $\tau := \bigcup\limits_{a}\sigma_a \cup id_{\mathbb{U}}$ is an automorphism of $C\mathcal{G}_1$. Thus the map $\rho := \tau \cup id_{C\mathcal{G}_2}$ is an automorphism of $C\mathcal{G}_1 \cup C\mathcal{G}_2$ that fixes $\mathbb{U}$ pointwise and satisfies $\rho \circ (\bigcup\limits_{a \in A} g_a h_a f_a) = (\bigcup\limits_{a \in A} g'_a h'_a f'_a)\circ \rho$.

\end{proof}

We shall now prove that this $C(H)$ defines a morphism of covers. We will proceed as in lemma \ref{C(h)_is_cover_of_U} and proposition \ref{C(h)_is_morphism_of_covers}.

\begin{lemma}
The theory of $C(H)$ is a cover of the theory of $\mathbb{U}$.

\end{lemma}

\begin{proof}
The same proof as in lemma \ref{C(h)_is_cover_of_U} works.
\end{proof}

\begin{prop}
The theory of $C(H)$ is a cover of both $C\mathcal{G}_1$ and $C\mathcal{G}_2$.
\end{prop}

\begin{proof}
Adapting the proof of proposition \ref{C(h)_is_morphism_of_covers} to the need for higher degree coherence morphisms, we find a way of extending automorphisms of $C\mathcal{G}_1$ into automorphisms of $(C\mathcal{G}_1 \cup C\mathcal{G}_2, \bigcup\limits_{a \in A} g_a h_a f_a)$ :

Let $\sigma_1$ be an automorphism of $C\mathcal{G}_1$. We may assume that it fixes $\mathbb{U}$ pointwise. Using the binding groupoid statement for $\mathcal{G}'_1$ in $C\mathcal{G}_1$ and condition (A) of the morphism of groupoids $H$, there exists, for each tuple $\overline{c}$, a \textit{unique} morphism $\sigma_{2, \overline{c}} \in \mathrm{Aut}_{\mathcal{G}'_2}(S_{2, \overline{c}})$ such that $\sigma_{2, \overline{c}} g_{\overline{c}} h_{\overline{c}} f_{\overline{c}} = g_{\overline{c}} h_{\overline{c}} f_{\overline{c}} \sigma_1|_{S_{\overline{c}}}$. By uniqueness again, the family of maps $(\sigma_{2,a})_{a \in A}$ is coherent. Thus, it defines an automorphism $\sigma_2$ of $C\mathcal{G}_2$ over $\mathbb{U}$, with the equality : $\sigma_2 g_a h_a f_a = g_a h_a f_a \sigma_1$.
\end{proof}

\begin{prop}
The definition of $C$ yields a functor $C : Iso(SFCG_A) \rightarrow AC_A$.
\end{prop}
\begin{proof}
As in the independent case, the action on the identities is easy enough to understand, and composition can be dealt with by picking literal compositions of morphisms when choosing the $h^{13}_a, f^{13}_a$ and $g^{13}_a$. See proposition \ref{C_is_functor_independent_case} for the details.

\end{proof}

\begin{rem}
Under additional hypotheses, such as finitariness of all the automorphism groups in the groupoids, or some stable embeddedness statements, we can define the functor $C$ on morphisms of groupoids, and not just isomorphisms of groupoids.
\end{rem}

\end{subsection}

\begin{subsection}{Comparing simplicial groupoids}
\begin{defi}\label{defi_epsilon_finite_gen}
Let $\mathcal{G}$ be an element of $SFCG_A$. Let $\mathcal{G}'$ be the extension of $\mathcal{G}$ defining the extra structure in the cover $C(\mathcal{G}) = (\mathbb{U}, S)$. Let $\mathcal{G}''(C(\mathcal{G}))$ be the binding groupoid of the cover $C(\mathcal{G}).$ Let $n$ be an integer. We define the set of maps $\varepsilon_{\mathcal{G}, n} := \lbrace g\circ f : O_{1, i_1} \rightarrow O_{2, i_2} \, |\, i_1 \in Ob(\mathcal{G}_n), i_2 \in Ob(GC(\mathcal{G})_n), \exists \overline{c} \in A^{=n}, \, g \in Hom_{\mathcal{G}''(C(\mathcal{G}))}(S_{\overline{c}}, O_{2, i_2}) \wedge f \in Hom_{\mathcal{G}'}(O_{1, i_1}, S_{\overline{c}})\rbrace$.
\end{defi}

\begin{prop}
The sets of maps $\varepsilon_{\mathcal{G}, n}$ define an isomorphism $\varepsilon : \mathcal{G} \rightarrow \mathcal{G}(C(\mathcal{G}))$ in $SFCG_A$.
\end{prop}

\begin{proof}
The same proof as in proposition \ref{epsilon_independent_case} shows that, in each degree, the set of maps $\varepsilon_{\mathcal{G}, n}$ is an isomorphism of groupoids. 
It remains to show that this family of maps is compatible with the inclusions that define the simplicial structures of $\mathcal{G}$ and $\mathcal{G}(C(\mathcal{G}))$. This can be checked using the same ideas as in proposition \ref{G(h)_is_a_simplicial_morphism}, namely, extending the groupoid morphisms to $\mathcal{G}'$ and $\mathcal{G}'(C(\mathcal{G}))$, and checking the conditions there using the distinguished objects $S_{\overline{c}}$.
\end{proof}
\end{subsection}
\end{section}

\begin{subsection}{The equivalence of categories}
\begin{theo}
The functors $G : Iso(AC_A) \rightarrow Iso(SFCG_A)$ and $C : Iso(SFCG_A) \rightarrow Iso(AC_A)$ are equivalences of categories.
\end{theo}

\begin{proof}
The only statement not already proved is the fact that $\varepsilon : id_{SFCG_A} \rightarrow GC$ and $\eta : id_{AC_A} \rightarrow CG$ are natural isomorphisms. Working fiberwise, this is implied by the case of internal covers.
\end{proof}

\begin{conj}
The categories $AC_A$ and $SFCG_A$ are equivalent. Moreover, this equivalence can be proved by extending the functor $C : Iso(SFCG_A) \rightarrow Iso(AC_A)$ to non-isomorphisms using a similar definition.
\end{conj}

\end{subsection}

\end{section}

\begin{section}{Computations}\label{section_computations}

In this section, we intend to understand the general constructions in two examples.

\begin{subsection}{Example : $(\mathbb{Z}/4\mathbb{Z})^{\omega}$}

In this subsection, we will try and apply our general constructions to the example of the two-sorted structure $\mathbb{V}:=((\mathbb{Z}/2\mathbb{Z})^{\omega},(\mathbb{Z}/4\mathbb{Z})^{\omega}, \iota, \pi)$.
Here, the maps $\iota$ and $\pi$ are the group morphisms that appear in the following exact sequence : 

\begin{center}
    
\begin{tikzcd}

0 \arrow[r] & (\mathbb{Z}/2\mathbb{Z})^{\omega}\arrow[r, hook, "\iota"]&(\mathbb{Z}/4\mathbb{Z})^{\omega}\arrow[r, "\pi", hook]&(\mathbb{Z}/2\mathbb{Z})^{\omega} \arrow[r]& 0
\end{tikzcd}

\end{center}

Note that the $\mathbb{F}_2$-vector space $(\mathbb{Z}/2\mathbb{Z})^{\omega}$ does not eliminate imaginaries, so we have to work in $((\mathbb{Z}/2\mathbb{Z})^{\omega})^{eq}$

\begin{prop}
Let $(\sigma_{\overline{c}})_{\overline{c} \in A^{\leq 3}}$ be an element of the projective limit of the restricted system : $(\mathrm{Aut}(S_{\overline{c}}))_{\overline{c}\in A^{\leq 3}}$.

Then this family induces an automorphism of the structure $(\mathbb{U}, S)$ over $\mathbb{U}$.
\end{prop}

\begin{proof}
We define the map $\sigma := \bigcup\limits_{a \in A}\sigma_a$. 

We first compute that, on each fiber $S_a$, the map $\sigma_a$ is of the form $x \mapsto x + x_a$, where $x_a$ is some element of $S_0$. Indeed, if $x, y \in S_a$, then there exists some $z \in (\mathbb{Z}/2\mathbb{Z})^{\omega}$ such that $x-y=\iota(z)$. Since $\sigma_a$ is an automorphism of $S_a$ over $\mathbb{U}$, we compute : $\sigma_a(x) - \sigma_a(y) = \sigma_a(\iota(z))=\iota(\sigma_a(z))=\iota(z) = x-y.$ Thus $\sigma_a(x)-x = \sigma_a(y)-y =:x_a.$

Then, since each $\sigma_a \cup \sigma_b \cup \sigma_c$ is an automorphism of $(\mathbb{U}, S_a \cup S_b \cup S_c)$, the following condition appears : $x_{a+b} = x_a + x_b$ for all $a, b \in A$. In particular, $x_0=0$.
Therefore, the map $\sigma \cup id_{\mathbb{U}}$ is a group automorphism commuting with $\iota$ and $\pi$. It is thus an automorphism of $(\mathbb{U}, S)$.
\end{proof}

\begin{rem}
In this case, the degenerate simplicial groupoid defined with $\mathcal{G}'_n := \mathcal{G}'_3$ if $n\geq 3$, and identities in higher degrees instead of the inclusion morphisms, captures the whole structure of the cover $(\mathbb{U},S)$ over $\mathbb{U}$. 
\end{rem}

\begin{coro}
The sort $(\mathbb{Z}/4\mathbb{Z})^{\omega}$ is not internal to $(\mathbb{Z}/2\mathbb{Z})^{\omega}$.
\end{coro}

\begin{proof}
The automorphism group $\mathrm{Aut}(S / \mathbb{U})$ is isomorphic to the group $\mathrm{End}((\mathbb{Z}/2\mathbb{Z})^{\omega})$. The latter has cardinality $2^{2^{\aleph_0}}$. Thus, it cannot be interpreted in $(\mathbb{Z}/2\mathbb{Z})^{\omega}$.
\end{proof}

\begin{rem}
In this example, the groupoids in degree $1$ can be chosen to be \textit{canonical} ones, i.e. with one object in each isomorphism class. Namely, the object is the kernel $K$ that appears in the short exact sequence, and the bijections $K \simeq S_a$ are the translations. 

However, in higher degrees, the non-existence of a section of the short exact sequence makes canonical groupoids harder, or even impossible to find. For instance, if $a+b=c$ in $(\mathbb{Z}/2\mathbb{Z})^{\omega}$, then you have to create one object in $\mathcal{G}_{a,b,c}$ for each possible value of $(x_a+x_b-x_c) \in K$, if $x_a \in S_a, x_b \in S_b$ and $x_c \in S_c$. These values classify the \textquotedblleft kinds of bijections" $K\sqcup K \sqcup K \simeq S_a \sqcup S_b \sqcup S_c$ mentioned in remark \ref{types_of_bijections}.
And, since sections do not exist, the 0-definable value $0 \in K$ is usually invalid, since no bijection realizes it.

\end{rem}

\end{subsection}

\begin{subsection}{Example : $d(\frac{dx}{x})=0$}

    Let $\mathbb{M}\models DCF_0$. Let $C = \lbrace x \, | \, dx = 0 \rbrace$ be the field of constants.
    Let's consider the following definable subgroup of $\mathbb{M}^{\times}$ : $S:=\lbrace x \in \mathbb{M}^{\times} : d(\frac{dx}{x})=0 \rbrace$.

    We notice that the 2-sorted structure $(C, S)$ comes with a definable short exact sequence $1 \rightarrow C^{\times} \rightarrow S \rightarrow C \rightarrow 0$.

\begin{prop}\label{automorphisms_DCF}
Let $(\sigma_{\overline{c}})_{\overline{c} \in C^{\leq 3}}$ be an element of the projective limit of the restricted system : $(\mathrm{Aut}(S_{\overline{c}}))_{\overline{c}\in C^{\leq 3}}$. In other words, we are given a family of automorphisms $\sigma_a \in \mathrm{Aut}(S_a / \mathbb{U})$, such that, for each triple $(a,b,c)$, the map $\sigma_a \cup \sigma_b \cup \sigma_c$ is an automorphism of $S_a \cup S_b \cup S_c$ over $\mathbb{U}$.

Then this family induces an automorphism of the structure $(C, S)$ over $C$.
\end{prop}

\begin{proof}
We define the map $\sigma := \bigcup_{a \in C}\sigma_a$.
We know that $\sigma_a : x \mapsto x \times \lambda_a$, where $\lambda_a$ is some element of $S_0=C^{\times}$.

When considering $\sigma_{a,b,a+b}$, the following condition appears : $\lambda_{a+b} = \lambda_a \times \lambda_b$ for all $a, b \in C$. In particular, $\lambda_0=1$.

We now have to check that this condition is enough to get an automorphism of the structure induced by the differential field $\mathbb{M}$.

Using quantifier elimination in $DCF_0$, we reduce to the case of formulas of the form  $Q(x_1,...,x_n)=P(x_1,d(x_1), ..., x_n, d(x_n),...,d^m(x_n))=0$, where the coefficients of $P$ are in $C$. We already know that $\sigma$ preserves multiplication and differentiation on $C, S$. 

\begin{claim}
The fibers $S_a$ are linearly independent over $C$.
\end{claim}
\begin{proof}
The differential is $C$-linear, and each fiber $S_a$ is made of eigenvectors associated to the eigenvalue $a$. An argument from linear algebra tells us that eigenvectors associated to distinct eigenvalues have to be linearly independent.
\end{proof}

Let us assume that $x_1,...,x_n \in S$ are roots of the differential polynomial $Q$ above. We write $Q(x_1,...,x_n)$ as a sum of products of elements of $S$. By the claim, all these terms are in the same fiber. Let $t_1,...,t_k \in S_a$ be these terms. Each of them is a product of some of the $d^{j}(x_i)$, with possibly a coefficient in $C$. By hypothesis, on $S_a$, the map $\sigma$ acts as multiplication by $\lambda_a$. Thus, we have $\sigma(t_1)+...+\sigma(t_k) = \lambda_a (t_1+...+t_k) = 0$. Moreover, we recall that $\sigma$ preserves multiplication in $S$, differentiation and product with elements of $C$. So each $\sigma(t_r)$, for $r=1,...,k$, is a product of some of the $d^j(\sigma(x_i))$, with the same coefficient in $C$ as the term $t_r$. Thus, $P(\sigma(x_1),d(\sigma(x_1)),..., d^m(\sigma(x_n)))=0$, as desired.

\end{proof}
\end{subsection}

\begin{prop}\label{exist_sections}
Let $G \leq C$ be a finitely generated subgroup. Let $s : G \rightarrow S$ be a group morphism such that $f \circ s = id_G$, where $f : S \rightarrow C$ is the surjection in the short exact sequence mentioned above.

Then there exists a section of $f$ which is a group morphism that extends $s$.
\end{prop}

\begin{proof}
We will need the following lemma : 

\begin{lemma}
Let $x=x_1$ be an element of $S$. Then, there exists a coherent system of roots of $x$. More precisely, there exists a family $(x_n)_{n \geq 1}$ of elements of $S$ such that, for all $k,l \geq 1$, we have $(x_{kl})^k=x_l$.
\end{lemma}
\begin{proof}
We use $\omega$-saturation of the ambient structure. In fact, the conditions can be expressed straightforwardly in a type, with $\omega$ variables (one for each $x_n$) and $1$ parameter (which is $x$). To prove consistency, we only need to pick a root of $x_1$ of high enough degree, and take its powers.

Finally, to prove that the $x_n$ are in $S$, we use the formula : $n\frac{dx_n}{x_n}=\frac{d((x_n)^n)}{(x_n)^n}= \frac{dx}{x} \in C$. Thus $\frac{dx_n}{x_n} \in C$.
\end{proof}

Now, since $G$ is a finitely generated torsion-free abelian group, we can find a $\mathbb{Z}$-basis $g_1,...,g_k \in G$. This basis is also a $\mathbb{Q}$-basis of the vector space $V$ generated by $G$. Let $W \leq C$ be a $\mathbb{Q}$-vector space such that $V \oplus W = C$. Note that this direct sum is also a direct sum of abelian groups. Thus, it suffices to define a group morphism on $W$ that is section of $f : S \rightarrow C$, and to extend $s : G \rightarrow S$ into a group morphism $s : V \rightarrow S$ such that $f\circ s = id_V$. 

\noindent We will detail the construction for $V$, the one for $W$ being slightly easier. For $i=1,...,k$, let $(h_{i, \frac{1}{n}})_n$ be a coherent system of roots of $s(g_i)$. If $x = \lambda_1 g_1 + ... + \lambda_k g_k \in V$, where the  $\lambda_i=\frac{a_i}{b_i}$ are in $\mathbb{Q}$, we define $s(x):=\prod\limits_{i=1}^k (h_{i, \frac{1}{b_i}})^{a_i}$.
Note that, if $\frac{a}{b} = \frac{c}{d}$, then $h_{i,\frac{1}{d}}=(h_{i,\frac{1}{ad}})^a=(h_{i,\frac{1}{bc}})^a$. Thus, $(h_{i,\frac{1}{d}})^c=(h_{i,\frac{1}{bc}})^{ac} = (h_{i,\frac{1}{b}})^a$. So $s$ is well-defined on $V$.

Let us now check that $s$ is a group morphism. Let $x=\frac{a_1}{b_1} g_1 + ... + \frac{a_k}{b_k}g_k$ and $y=\frac{c_1}{d_1} g_1 + ... + \frac{c_k}{d_k}g_k$. So $x+y = \frac{a_1d_1+c_1b_1}{b_1d_1}g_1+...+\frac{a_kd_k+c_kb_k}{b_kd_k}g_k$.

By definition, we have $s(x+y)=\prod\limits_{i=1}^k (h_{i, \frac{1}{b_id_i}})^{a_id_i+b_ic_i}$.

$=\prod\limits_{i=1}^k (h_{i, \frac{1}{b_id_i}})^{a_id_i} \times \prod\limits_{i=1}^k (h_{i, \frac{1}{b_id_i}})^{b_ic_i}$

$= \prod\limits_{i=1}^k (h_{i, \frac{1}{b_i}})^{a_i} \times \prod\limits_{i=1}^k (h_{i, \frac{1}{d_i}})^{c_i}$

$=s(x) s(y)$.

\noindent Thus, $s$ is a group morphism. We notice that it extends the morphism that was defined on $G$. It remains to check that it is a section of $f$. Let $1\leq i\leq k$, let $b \geq 1$. We want to prove that $fs(\frac{1}{b}g_i)=\frac{1}{b}g_i$. Note that this is enough, since $id_V$ and $fs$ are group morphisms, and the $\frac{1}{b}g_i$ generate the abelian group $V$. 

\noindent We compute : $b f(s(\frac{1}{b}g_i))=f((s(\frac{1}{b}g_i))^b)=fs(g_i)=g_i$, since $s|_G$ was assumed to be a section of $f$. Thus, since this equality holds in $C$, which is torsion-free and divisible, we have $\frac{1}{b} g_i = f(s(\frac{1}{b}g_i))$, as desired.

\end{proof}

\begin{prop}\label{unique_morphism_extension}
Let $B=\lbrace b_1,...,b_k \rbrace$ be a finite subset of $C$. Let $f : B \rightarrow S$ be such that, whenever there is a relation of the form $\sum_i n_i b_i = 0$, then the relation $\prod_i f(b_i)^{n_i}=1$ holds.

Then, $f$ extends uniquely into a group morphism $f : <B> \rightarrow S$.
\end{prop}

\begin{proof}
Uniqueness is clear, for the extension has to be defined by the following formula : $\sum_i n_i b_i \mapsto \prod_i f(b_i)^{n_i}$. For existence, it suffices to check that, if $\sum_i n_i b_i = \sum_i m_i b_i$, then $\prod_i f(b_i)^{n_i}=\prod_i f(b_i)^{m_i}$, for all tuples of integers $(n_1,...,n_k), (m_1,...,m_k)$. This is precisely given by the hypothesis on $f$.
\end{proof}

\begin{prop}\label{explicit_binding_simplicial_groupoid_DCF}
The binding simplicial groupoid of $(\mathbb{U}, S)$ is type-definable, and can be chosen to be canonical.
\end{prop}
\begin{proof}
Let $n$ be an integer. We shall now define the degree n binding groupoid. For each element $c \in A^{=n}$, we define the only object $O_c$ as the set $\bigcup_{a \in c} C^{\times} \times \lbrace a \rbrace$.

Let $X_c$ be the set of tuples of couples of the form $((x_1, a_1), (x_2,a_2),...,(x_n, a_n))$, where the tuple $(a_1,...,a_n)$ defines the same set as $c$, and $x_1,...,x_n$ are elements of $C^{\times}$. 

Let $E_c$ be the quotient of $X_c$ under the relation \textquotedblleft being equal up to permutation".

We now define the set $M_c$ as the set of elements $m$ of $E_c$ such that, for any representative $((x_1, a_1), (x_2,a_2),...,(x_n, a_n))$ of $m$, for any tuple of integers $(k_1,...,k_n)$, if the relation $\sum_i k_i a_i = 0$ holds in the group $C$, then the relation $\prod_i x_i^{k_i}=1$ holds in the group $C^{\times}$.

Note that the set $M_c$ always contains the class of $((1,a_1),..., (1, a_n))$. Moreover, it is type-definable, uniformly over $A^{=n}$. The action of $M_c$ on $O_c$ is given by fiberwise translation.

We thus get a type-definable groupoid in $\mathbb{U}$. Using propositions \ref{exist_sections} and \ref{unique_morphism_extension}, it is not hard to extend this groupoid by adding the object $S_c$ in the component over $c$, for $c \in A^{=n}$. 

Indeed, let us define the morphisms $O_c \rightarrow S_c$ as the fiberwise translations by elements of the $S_a$, for $a \in c$, that satisfy the condition of proposition \ref{unique_morphism_extension}. 
Then, the group of groupoid automorphisms of $S_c$ is the group of fiberwise translations by elements of $C^{\times}$ that also satisfy the condition of proposition \ref{unique_morphism_extension}. Thus, by propositions \ref{exist_sections} and then \ref{automorphisms_DCF}, these bijections are in fact elements of $\mathrm{Aut}(S_c / \mathbb{U})$.

\end{proof}

\begin{rem}

Contrary to the previous example, the groupoids can be chosen to be \textit{canonical} in all degrees, for there exist enough sections $C \rightarrow S$. In other words, the distinguished value $1$ (the identity element in $S$) can be picked. However, these groupoids are \textit{type-definable} without parameters. In fact, the hypothesis of uniform finite generatedness of the languages does not hold, as we shall prove in the following proposition :
\end{rem}

\begin{prop}\label{no_uniform_finite_generatedness_DCF_0}
Let $n \geq 2$. Then, the languages of the covers $(S_{\overline{c}}, \mathbb{U})$ are not uniformly finitely generated over the language of $\mathbb{U}$, for $\overline{c} \in A^{=n}$, even though each language is finitely generated.
\end{prop}

\begin{proof}
We shall prove the fact for $n=2$, the other cases being similar. By contradiction, assume that the languages of the internal covers $(S_a \cup S_b, \mathbb{U})$ are uniformly finitely generated, for $a,b \in A$, $a\neq b$. Then, the automorphism groups $\mathrm{Aut}(S_a\cup S_b / \mathbb{U})$, along with their actions on the $S_a\cup S_b$, are uniformly definable in $(\mathbb{U}, S)$. Let $\phi(a,b,m)$ be a formula defining \textquotedblleft $m$ is the code of an automorphism of $S_a \cup S_b$ over $\mathbb{U}$." Let $\psi(a,b,m,s, s')$ be a formula defining the graphs of the maps defined by such codes $m$.
We may assume that $\psi$ is given as in proposition \ref{explicit_binding_simplicial_groupoid_DCF} above, i.e., the maps are fiberwise translations by elements of the kernel $S_0$.

Moreover, propositions \ref{automorphisms_DCF}, \ref{exist_sections} and \ref{unique_morphism_extension} give us a set of definable relations that generate the languages at stake. Indeed, if $k,l \in \mathbb{Z}$, we define the relation $R_{k,l}(x,y)$ as $x^k = y^l$. Then, by the results proved in the propositions mentioned above, we have :

$\lbrace$  \textquotedblleft $m$ defines, via $\psi(a,b,m,s,s')$, a map on $S_a \cup S_b$ which preserves the relation induced by $R_{k,l}$ on $(S_a\cup S_b)^2$" \, | \, $k,l \in \mathbb{Z}$    $\rbrace \models \phi(a,b,m) $. Here, the variables are both $a,b$ and $m$.

Thus, by $\omega$-saturation and compactness, a finite fragment of the set of relations $R_{k,l}$ suffices to ensure that the map defined by some code $m$ is an automorphism. The crucial point is that this happens \textit{uniformly} in $a,b$. Now, take an integer $N$ that bounds the size of the integers $k,l$ that appear in such a finite fragment. Pick some non-zero element $a \in A$. The relation $R_{2N,1}(x,y)$ has to be preserved by any automorphism of $S_a\cup S_{2N.a}$.
However, there exist tuples $m$ that code permutations of $S_a \cup S_{2N.a}$ that preserve the relations $R_{k,l}$, for $|k|, |l| \leq N$, and which do not preserve the relation $R_{2N,1}$. Simply take any \textquotedblleft generic" pair of elements $\lambda, \mu$ of $S_0$ - in particular, with $\lambda^{2N} \neq \mu$ - and consider the map $m_{\lambda, \mu}$ which acts as translation by $\lambda$ and $\mu$ on $S_a$ and $S_{2N.a}$ respectively.

Now, the relations induced by the $R_{k,l}$, for $|k|, |l| \leq N$, on $(S_a \cup S_{2N.a})^2$ are trivial, and thus automatically preserved by the map $m_{\lambda, \mu}$. However, the relation $R_{2N,1}$ is not preserved by this map, which is a contradiction. So uniform finite generatedness of the languages does not hold.

Now, let us show that each language is finitely generated. For each tuple $\overline{c} \in A^{=n}$, let $G_{\overline{c}}$ be the subgroup of $\mathbb{Z}^n$ which is the kernel of the map $(k_1,...,k_n) \mapsto k_1 c_1 + ... + k_n c_n$. Since $\mathbb{Z}$ is noetherian, we know that $G_{\overline{c}}$ is finitely generated, as an abelian group. Now, any finite system of generators defines a finite family of definable relations that generates the language of $(\mathbb{U}, S_{\overline{c}})$ over that of $\mathbb{U}$. 

\end{proof}

\begin{theo}
The group $S$ is isomorphic, but not definably so, to the direct product $C \times C^{\times}.$ The cover $(C , S)$ is interpretable in the field $C$.
\end{theo}

\begin{proof}
The existence of group-theoretic sections of the short exact sequence $1 \rightarrow C^{\times} \rightarrow S \rightarrow C \rightarrow 0$ implies that the group $S$ is isomorphic to a semidirect product $C^{\times} \rtimes C$. Since the group $S$ is abelian, this semidirect product is in fact a direct product.

Note that the cover is not internal, for there are at least $2^{|C|}$ elements in the automorphism group $\mathrm{Aut}((C, S) / C)$. Thus, there are no \textit{definable} group isomorphisms $S \simeq C^{\times} \times C$.

Interpretability in the field $C$ follows, for the proof of proposition \ref{automorphisms_DCF} shows that the structure induced on $S$ by the differential field is only the group-theoretic structure induced by the short exact sequence mentioned above.
\end{proof}

\begin{rem}
Here, group isomorphisms $S \simeq C^{\times} \times C$ are built using families of morphisms $f_a : C^{\times} \rightarrow S_a$ that are coherent relative to the set-theoretic inclusions on both sides. Since the inductive construction of such coherent families involves many choices, it is not surprising to find that they are not definable. 

Moreover, the theorems proved in the next section should help understand why interpretability of $S$ in $C$, in this example, is not a coincidence. 
\end{rem}

\end{section}

\begin{section}{Type-definable simplicial groupoids}\label{type_def_groupoids}

In this section, we present the type-definable groupoids that appear when the technical hypotheses of uniformly finitely generated languages are not satisfied, just as in the example above. We shall call them \textquotedblleft *-definable" or \textquotedblleft type-definable" indiscriminately, for these are the only kinds of type-definable groupoids we will consider here.

From now on, we study 1-analysable covers that may not satisfy the hypothesis of (Uniform Finite Generatedness) of the languages, and keep only the (Local Stable Embeddedness) hypothesis : For all finite subsets $C, D \subseteq A$, the structure $(\mathbb{U}, \bigcup\limits_{a \in C}S_{a})$, with the structure induced by the $C$-definable sets of $(\mathbb{U}, S)$, is stably embedded in $(\mathbb{U}, \bigcup\limits_{a \in C\cup D}S_{a})$.

\begin{defi}\label{defi_type_def_groupoids}
\begin{enumerate}
    \item A concrete groupoid is \textit{$*$-definable} if its set of objects and its set of morphisms are type-definable without parameters, with the types having possibly infinitely many variables.

We still require the morphisms and objects to be definable maps and definable sets, whose definitions only use finite fragments of the tuples involved. We shall refer to these finite fragments as the \textit{concrete parts} of the infinite tuples. 

We weaken the assumption that two distinct objects (resp. morphisms) define distinct sets (resp. bijections). Instead, we only require that two distinct objects (resp. morphisms) that define the same set (resp. bijection) have the same concrete part. 

So, each object (resp. each morphism) in the \textit{abstract} groupoid is defined by an infinite tuple, and such an infinite tuple contains a finite tuple which is the code, in the model-theoretic sense, of a definable set (resp. a definable bijection). Type-definability means that finite tuples (that code definable bijections) code morphisms between objects if and only if they realize a specific type over the infinite tuples that code the objects.
Note that the infinite tuples defining the morphisms are merely made of a concrete part, and the infinite tuples defining the source and target of the morphism.

We suggest the reader have a look at the proof of theorem \ref{theo_type_def_binding_simplicial_groupoids} to understand why we picked this notion for type-definable groupoids.

    \item A concrete simplicial groupoid is \textit{$*$-definable} if the groupoids in each degree and the inclusion morphisms are $*$-definable.

We will restrict to the case where the *-definable (simplicial) groupoids are naturally projective limits of 0-definable (simplicial) groupoids. More precisely, we ask that there exist coverings of the index sets of variables by finite subsets containing the concrete parts, with the following properties : \begin{itemize}
    \item The objects defined by restriction to these finite sets of coordinates are 0-definable (simplicial) groupoids.
    \item The union of two finite subsets belonging to the covering is also an element of the covering.
\end{itemize}
    
    \item A morphism of *-definable groupoids is a *-definable set of definable maps that satisfies conditions (A) and (B). Again, the types defining the sets of maps may have infinitely many variables, but each map is defined with a finite concrete part.
\end{enumerate}

\end{defi}

\begin{rem}
A morphism of *-definable groupoids yields a collection of morphisms of 0-definable groupoids between the projections of the two groupoids, by projection.
\end{rem}

\begin{subsection}{Type-definable binding simplicial groupoids}

\begin{theo}\label{theo_type_def_binding_simplicial_groupoids}
Let $\mathbb{U}'=(\mathbb{U}, S)$ be a 1-analysable cover of $\mathbb{U}$ over $A$ that satisfies the hypothesis of (Local Stable Embeddedness). Then there exist simplicial groupoids $\mathcal{G}$ and $\mathcal{G}'$, over the set $A$, that are *-definable in $\mathbb{U}$ and $\mathbb{U}'^{eq}$ respectively, with the following properties:

\begin{itemize}
    \item The simplicial groupoid $\mathcal{G}'$ is an extension of $\mathcal{G}$ with the extra object $S_{\overline{c}}$ in the connected component over $\overline{c}$. Moreover, the set-theoretic inclusions between the $S_{\overline{c}}$ belong to the simplicial groupoid.
    \item The automorphism groups $\mathrm{Aut}_{\mathcal{G}'}(S_{\overline{c}})$, along with their action on $S_{\overline{c}}$, are isomorphic to the groups $\mathrm{Aut}(S_{\overline{c}}/\mathbb{U})$.
\end{itemize}
\end{theo}

\begin{proof}
Let $n$ be an integer. We wish to build the degree $n$ groupoid. 

First, by compactness, internality and elimination of imaginaries in $\mathbb{U}$, there exist a 0-definable set $B \subseteq \mathbb{U}$, a 0-definable surjective map $B \rightarrow A^{=n}$ and 0-definable sets $X \subseteq \mathbb{U}, \, F \subseteq B \times S \times X,$ such that, for all $c \in A^{=n}$, for all $b \in B_c$, the set $F_b$ is the graph of an embedding $F_b : S_c \hookrightarrow X$.

Let $R(x)$ be a 0-definable relation in $\mathbb{U}'$. Let $c$ be an element of $A^{=n}$. For the same reasons as in theorem \ref{binding_groupoid_of_internal_cover}, there exists a formula $\phi(z, \alpha_c)$, where $\alpha_c \in \mathbb{U}$ such that $\models \exists z \in B_c, \, \phi(z, \alpha_c)$ and $\models \forall z,t \in B_c, \, (\phi(z, \alpha_c)\wedge \phi(t, \alpha_c))\rightarrow $\textquotedblleft $F_t^{-1} \circ F_z \cup id_{\mathbb{U}}$ preserves the relation induced by $R(x)$ on $(\mathbb{U}, S_c)$ .$"$

Now, by compactness, we can find a formula $\phi_R(z, w_R)$ that satisfies the above statement uniformly for $c \in A^{=n}$.
Then, keeping the same notion of \textquotedblleft kinds of bijections" as in theorem \ref{binding_groupoid_of_internal_cover}, we define the set of codes of the objects of $\mathcal{G}_n$ as the set of infinite tuples $(u, c, (w_R)_R)$ that satisfy $\lbrace c \in A^{=n} \rbrace \cup \lbrace \exists t \in B_c, \, \, \phi_{R_1}(t, w_{R_1}) \wedge ... \wedge \phi_{R_k}(t, w_{R_k}) \wedge $\textquotedblleft $u$ is the code of the set $Im(F_t)" \, | \, k \in \mathbb{N}, \, R_1,...,R_k $: 0-definable relations $ \rbrace $

The set defined by such an infinite tuple $(u,c,(w_R)_R)$ is merely the set coded by $u$. Note that here, two distinct objects in the groupoid may have the same underlying set.

Now, the set of morphisms from $S_c$ to the object corresponding to the infinite tuple $(u, c, (w_R)_R)$ is coded by the set of tuples $(t, c, (w_R)_R)$ realizing the type $\lbrace t \in B_c \rbrace \cup \lbrace \phi_R(t, w_R) \, | \, R $ : 0-definable relation$  \rbrace $. The bijections coded by such tuples $(t, c, (w_R)_R)$ are just the maps $F_t$.

\noindent Note that these sets are non-empty by saturation, and uniformly *-definable. Moreover, the construction gives precisely a projective limit of groupoids, where the covering can be indexed by the finite collections of 0-definable relations.

Now, the rest of the construction follows just as in theorem \ref{theo_defi_simplicial_groupoid}. For instance, the inclusions are conjugates of the set-theoretic ones $S_c \hookrightarrow S_d$ by morphisms in the components over $c$ and $d$. They can  be coded by pairs of codes of morphisms, and are thus *-definable. Again, just as in theorem \ref{theo_defi_simplicial_groupoid}, condition (A) for the inclusions is given by the technical stable embeddedness hypotheses.

\end{proof}

\begin{rem}
The simplicial groupoids built that way are finitely faithful and have the \textit{disjoint union property}, just as those given by theorem \ref{theo_defi_simplicial_groupoid}.
\end{rem}

\end{subsection}

\begin{subsection}{The cover associated to a type-definable simplicial groupoid}

\begin{prop}\label{exists_system_of_inclusions_*}
    Let $\mathcal{G}$ be a *-definable simplicial groupoid over $A$. Then there exists a commuting system of inclusions in $\mathcal{G}$.
\end{prop}

\begin{proof}
The same inductive construction as in proposition \ref{exists_system_of_inclusions} works, for the types involved in the saturation arguments still have few enough variables and parameters. 
\end{proof}
\begin{prop}\label{coherent_families_*}
Let $\mathcal{G}_1$, $\mathcal{G}_2$ be *-definable simplicial groupoids over $A$ with the disjoint union property. Let $H : \mathcal{G}_1 \rightarrow \mathcal{G}_2$ be a morphism of simplicial groupoids.
Let $\iota_1= ( O_{i_{\overline{c}}} \rightarrow O_{i_{\overline{d}}})_{\overline{c} \subset \overline{d} \in A^{<\omega}}$ and $\iota_2 = (O_{j_{\overline{c}}} \rightarrow O_{j_{\overline{d}}})_{\overline{c} \subset \overline{d} \in A^{<\omega}} $ be commuting systems of inclusions in $\mathcal{G}_1$ and $\mathcal{G}_2$ respectively. 

Then there exists a coherent family of maps $(m_a : O_{i_a} \rightarrow O_{j_a})_{a\in A}$ belonging to the morphism of simplicial groupoids $H$.

More generally, for all finite subsets $B \subset A$, any \textit{coherent} family of maps $(m_a : O_{i_a} \rightarrow O_{j_a})_{a\in B}$ can be extended into a coherent family of maps $(m_a : O_{i_a} \rightarrow O_{j_a})_{a\in A}$.
\end{prop}

\begin{proof}
The same inductive construction as in proposition \ref{coherent_families} works. The only detail to take care of is removing the existential quantifiers that express the existence of higher degree \textquotedblleft coherence morphisms", and replacing them with more variables to account for these \textit{uniquely determined} coherence morphisms. Since the types involved still have few enough parameters and variables, they can be realized.
\end{proof}

\begin{theo}\label{cover_from_groupoid_*}
Let $\mathcal{G}$ be a finitely faithful *-definable simplicial groupoid in $\mathbb{U}$ over $A$, with the disjoint union property. Then there exists a cover $(\mathbb{U}, O_*)$ of $\mathbb{U}$, along with an extension $\mathcal{G}'$ of $\mathcal{G}$ that is *-definable in $(\mathbb{U}, O_*)$, such that :

\begin{itemize}

    \item The cover $(\mathbb{U}, O_*)$ is 1-analysable over $A$, and satisfies the technical hypothesis given at the beginning of section \ref{type_def_groupoids}.
    \item The simplicial groupoid $\mathcal{G}'$ is an extension of $\mathcal{G}$ with the extra object $O_{*,\overline{c}}$ in the connected component over $\overline{c}$. Moreover, the set-theoretic inclusions between the $O_{*, \overline{c}}$ belong to the simplicial groupoid.
    \item The automorphism groups $\mathrm{Aut}_{\mathcal{G}'}(O_{*,\overline{c}})$, along with their action on $O_{*,\overline{c}}$, are isomorphic to the groups $\mathrm{Aut}(O_{*,\overline{c}}/\mathbb{U})$.
\end{itemize}

\end{theo}

\begin{proof}
Just as in the case of 0-definable simplicial groupoids, we wish to extend $\mathcal{G}$ by picking a commuting system of inclusions spanning all the connected components, and creating copies of the objects that appear in this system.

We use proposition \ref{exists_system_of_inclusions_*}. We can thus pick a commuting system of inclusions $\iota$ in the simplicial groupoid $\mathcal{G}$.

Then, just as before, we are able to create copies of objects using maps $f_{\overline{c}}$, in a coherent way :

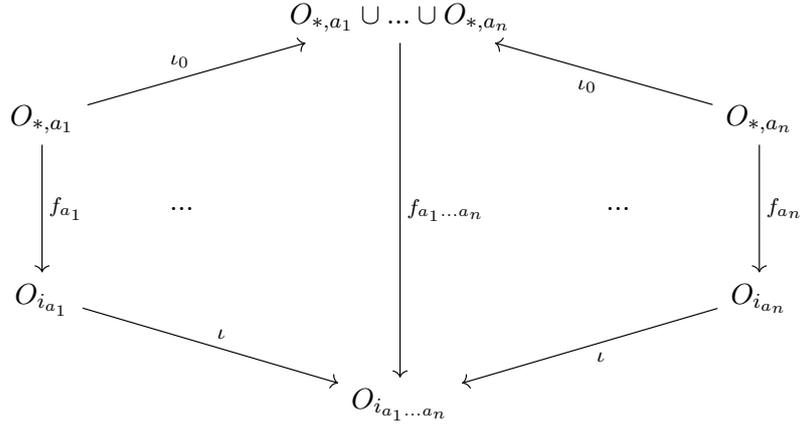
\begin{center}
\begin{tikzcd}
 && O_{*, a_1}\cup ... \cup O_{*, a_n}\arrow[dddd, "f_{a_1...a_n}"] \arrow[from=lld, "\iota_0"] \arrow[from=rrd, "\iota_0"] &

 \\
 O_{*,{a_1}}\arrow[dd, "f_{a_1}"]&&&&  O_{*,{a_n}}\arrow[dd, "f_{a_n}"]   
 \\
 & ... && ...
 \\
 O_{i_{a_1}}& && &  O_{i_{a_n}}
 \\
  && O_{i_{a_1...a_n}}  \arrow[from=llu, "\iota"] \arrow[from=rru, "\iota"] &
 \end{tikzcd} 
\end{center}
Here, the $\iota_0$ are the set-theoretic inclusions, and the $\iota$ are in the commuting system picked earlier. Note that we only consider the case of pairwise distinct $a_i$.

Again, apart from the new object sort $O_*$, we create new morphism sorts $M_*^{(n)}$ and $N_*^{(i,j)}$, by creating copies of the \textit{concrete parts} of the tuples defining the morphisms and inclusion morphisms. 

Now, we use the projective limit property of $\mathcal{G}$. We pick an appropriate covering of the sets of coordinates by finite subsets.
To fix notations, let $I_n$ be the set of variables in the types defining the degree $n$ groupoid. We cover $I_n$ with finite subsets $F_{i,n}$, for $i$ belonging to some set $J$.
We let $\mathcal{G}_i$ denote the 0-definable simplicial groupoid given by the restrictions to the $F_{i,n}$. 
If $i,j \in J$, we let $i \vee j \in J$ denote the index corresponding to the unions $F_{i,n} \cup F_{j,n}$.

We thus get back to the case of theorem \ref{cover_from_simplicial_groupoid}, by looking at each $i \in J$. This way, we add definable relations to the structure $(\mathbb{U}, O_*)$ that still yield a cover of $\mathbb{U}$, and build the only 0-definable extension $\mathcal{G}'_i$ of $\mathcal{G}_i$ for which the copying maps $f_c : O_{*,c} \rightarrow O_{i_c}$ we already picked are isomorphisms. To be clear, the new relations define which \textquotedblleft concrete parts" of morphisms, i.e. elements of the $M_*^{(n)}$ or $N_*^{(i,j)}$, are morphisms between given objects of $\mathcal{G}'_i$.

Then, we combine all these structures together. We still have a 1-analysable cover of $\mathbb{U}$. In fact, because of the finite union property, any relation definable in the whole structure is definable using only the structure of one of the $\mathcal{G}'_i$.

We show that the projective limit of the $\mathcal{G}'_i$ is an extension $\mathcal{G}'$ of $\mathcal{G}$ that is *-definable and satisfies the conditions of the theorem. 

The inclusion $\mathrm{Aut}_{}(O_{*, c} / \mathbb{U})\subseteq \mathrm{Aut}_{\mathcal{G}'}(O_{*, c})$ is proved the same way as in theorem \ref{cover_from_simplicial_groupoid}.

Let $\sigma \in \mathrm{Aut}_{\mathcal{G}'}(O_{*, c})$. We wish to extend $\sigma$ into an automorphism of the full structure $(\mathbb{U}, O_*)$. This amounts to extending $\sigma$ into an automorphism of the simplicial groupoid $\mathcal{G}'$. Applying proposition \ref{coherent_families_*}, we find a coherent family of automorphisms, in the groupoid sense, of the $O_{*, d}$, for $d \in A^{< \omega}$, that extends $\sigma$. Such a family defines an automorphism of the full structure $(\mathbb{U}, O_*)$.

Finally, let us show that the required stable embeddedness conditions hold. Let $c \subset d \in A^{< \omega}$. Let $\sigma$ be an automorphism of $(\mathbb{U}, O_{*,c})$, which we may assume to fix $\mathbb{U}$ pointwise. Thus, we have $\sigma \in \mathrm{Aut}_{\mathcal{G}'}(O_{*, c})$, and the previous paragraph shows that $\sigma$ can be extended into an automorphism of $(\mathbb{U}, O_{*,d})$.
\end{proof}

\begin{prop}
The cover built from the simplicial groupoid $\mathcal{G}$, up to isomorphism, does not depend on the choices made.
\end{prop}

\begin{proof}
The only fact not shown before is that two choices of coverings define the same structure. Let $J_1, J_2$ be index sets for coverings of the set of variables that satisfy the two conditions given in definition \ref{defi_type_def_groupoids}. Since the underlying sets of the new sorts $O_*, M_*^{(n)}$ can be assumed to be fixed, we only have to show that the definable relations are the same.

\noindent Let $i \in J_1$. We want to show that the extended simplicial groupoid $\mathcal{G}'_i$ is definable with the structure induced by some $\mathcal{G}'_j$, for some $j \in J_2$. Since the set defined by $i$ is finite, and $J_2$ defines a covering, it can be \textquotedblleft covered" with finitely many elements $j_1,...,j_k \in J_2$ of the covering associated to $J_2$. So, let $j=j_1\vee j_2,...,\vee j_k \in J_2$ be the index of a subset covering that of $i$. It remains to notice that the simplicial groupoid $\mathcal{G}'_i$ is the projection of the simplicial groupoid $\mathcal{G}'_j$ on the coordinates coreesponding to $i$.

By symmetry, we conclude that the structure thus defined does not depend on the covering.

\end{proof}

\end{subsection}

\begin{subsection}{Comparing covers}

From now on, we let $\mathbb{V} = (\mathbb{U}, S, p : S\rightarrow A)$ be a 1-analysable cover satisfying the  conditions given at the beginning of section \ref{type_def_groupoids}. The simplicial groupoids $\mathcal{G}(\mathbb{V})$ and $\mathcal{G}'(\mathbb{V})$ are the binding simplicial groupoids that are *-definable in $\mathbb{U}$ and $\mathbb{V}^{eq}$ respectively.
As in the case of 0-definable simplicial groupoids, we shall also work with the cover $C\mathcal{G}(\mathbb{V}) = (\mathbb{U}, O_*)$ and the simplicial groupoids $\mathcal{G}(C\mathcal{G}(\mathbb{V}))$ and $\mathcal{G}''(C\mathcal{G}(\mathbb{V}))$.

\begin{theo}\label{theo_comparing_covers_*}
There exists an isomorphism of covers $\eta_{\mathbb{V}} : \mathbb{V} \simeq C\mathcal{G}(\mathbb{V}).$
\end{theo}

\begin{proof}
Pick a coherent systems of inclusions $\iota = (O_{i_{\overline{c}}} \rightarrow O_{i_{\overline{d}}})$ in $\mathcal{G}(\mathbb{V})$. Then pick families of morphisms $(g_a : S_a \rightarrow O_{i_a})_a$ and $(h_a : O_{i_a} \rightarrow O_{*,a})_a$ that belong to $\mathcal{G}'(\mathbb{V})$ and $\mathcal{G}''(C\mathcal{G}(\mathbb{V}))$ respectively, and which are \textit{coherent}, with respect to $\iota$ and the set-theoretic inclusions. By proposition \ref{coherent_families_*}, such families exist. Then, just as in theorem \ref{theo_comparing_covers}, the structure $(\mathbb{V}\cup C\mathcal{G}(\mathbb{V}), \bigcup\limits_{a \in A} h_a g_a)$ defines an isomorphism of covers.

\end{proof}

\begin{rem}
This theorem shows that the *-definable groupoids built in theorem \ref{theo_type_def_binding_simplicial_groupoids} indeed capture the whole structure of the 1-analysable covers.
\end{rem}

\begin{coro}
Assume that the binding groupoid $\mathcal{G}$ of the cover $\mathbb{V}$ over $\mathbb{U}$ is canonical, and that there exists a commuting system of inclusions which is $0$-definable. Then the cover $\mathbb{V}$ is interpretable in $\mathbb{U}$.
\end{coro}
\begin{proof}
Here, the reconstruction of the cover $\mathbb{V}$ from its binding groupoid, as described in theorem \ref{cover_from_groupoid_*}, can be done inside $\mathbb{U}$ itself. Indeed, adding a new object in each connected component is the same as duplicating the whole type-definable simplicial groupoid, since the latter is canonical. Note that the copying maps in all degrees can be chosen to be 0-definable (over $A$, of course), precisely because we can choose the commuting system of inclusions in $\mathcal{G}$ to be 0-definable. 

Now, for each 0-definable projection of the type-definable simplicial groupoid $\mathcal{G}$, the corresponding extended simplicial groupoid is 0-definable, using conjugation by the 0-definable copying maps. Thus, the structure of the cover $C\mathcal{G}(\mathbb{V})$ is interpretable in $\mathbb{U}$. Since isomorphisms of covers are bi-interpreations, the structure $\mathbb{V}$ is therefore interpretable in $\mathbb{U}$.
\end{proof}

\end{subsection}

\begin{subsection}{Functoriality}

From now on, we let $AC^*_A$ denote the category of 1-analysable covers over $A$ that satisfy the (Local Stable Embeddedness) condition, and $SFCG^*_A$ denote the category of the *-definable finitely faithful simplicial groupoids over $A$ that satisfy the disjoint union property.

\begin{theo}
The categories $Iso(AC^*_A)$ and $Iso(SFCG^*_A)$ are equivalent.
\end{theo}
\begin{proof}
The definitions of the functors $C : Iso(SFCG^*_A) \rightarrow Iso(AC^*_A)$ and $G : AC^*_A \rightarrow SFCG^*_A$ are the same as in section \ref{section_functoriality_finite_generatedness}. The natural isomorphism $\varepsilon : GC \simeq id_{Iso(SFCG^*_A)}$ is also built as in definition \ref{defi_epsilon_finite_gen}.

\end{proof}
\end{subsection}

\end{section}

\printbibliography[
heading=bibintoc,
title={References}
]


Paul Z. WANG, Ecole Normale Superieure, Département de Mathématiques et Applications, 45 rue d'Ulm, 75005 Paris

\end{document}